\documentclass[reqno]{amsart}
\usepackage{amsmath, amsthm, amssymb, amstext}
\usepackage{hyperref}
\usepackage[dvipsnames]{xcolor}
\usepackage{orcidlink}
\hypersetup{pdfborder={0 0 0},colorlinks}
\allowdisplaybreaks
\usepackage{tikz}
\usetikzlibrary{patterns,calc}

\newcommand{\A}{{\mathcal A}}
\newcommand{\As}[1][]{A_\sigma #1}
\newcommand{\abs}[1]{\left|#1\right|}
\newcommand{\Bs}[1][]{B_\sigma #1}

\newcommand{\bgset}[1]{\big\{#1\big\}}
\newcommand{\D}{{\mathcal D}}

\newcommand{\ds}[1]{\displaystyle #1}
\newcommand{\eps}{\varepsilon}
\newcommand{\F}{{\mathcal F}}
\newcommand{\G}{{\mathcal G}}
\renewcommand{\H}{{\mathcal H}}
\newcommand{\loc}{{\rm loc}}
\newcommand{\M}{{\mathcal M}}

\newcommand{\norm}[2][]{\left\|#2\right\|_{#1}}
\renewcommand{\O}{\text{O}}
\renewcommand{\o}{\text{o}}
\newcommand{\PS}[1]{\textnormal{(\text{PS})$_{#1}$}}
\newcommand{\pnorm}[2][]{\if #1'' \left|#2\right|_p \else \left|#2\right|_{#1} \fi}
\newcommand{\R}{\mathbb R}
\newcommand{\rad}{\text{rad}}
\newcommand{\restr}[2]{\left.#1\right|_{#2}}
\newcommand{\seq}[1]{\left(#1\right)}
\newcommand{\set}[1]{\left\{#1\right\}}
\newcommand{\T}{{\mathcal T}}
\newcommand{\Tc}{{\mathcal T}_{\rm c}}
\newcommand{\Te}{{\mathcal T}_{\rm e}}
\newcommand{\vol}[1]{\left|#1\right|}
\newcommand{\wto}{\rightharpoonup}

\DeclareMathOperator{\supp}{supp}

\newtheorem{theorem}{Theorem}
\newtheorem{remark}[theorem]{Remark}
\newtheorem{lemma}[theorem]{Lemma}
\newtheorem{proposition}[theorem]{Proposition}
\newtheorem{corollary}[theorem]{Corollary}
\newtheorem{definition}[theorem]{Definition}

\newcommand{\dist}[2]{\text{dist}\, (#1,#2)}

\numberwithin{theorem}{section}
\numberwithin{equation}{section}

\title[A scaling approach to double phase equations in $\mathbb{R}^N$]{A scaling approach to\\double phase equations in $\mathbb{R}^N$}

\author[K. Perera]{Kanishka Perera\orcidlink{0000-0001-6168-247X}}
\address[K. Perera]{Department of Mathematics, Florida Institute of Technology, 150 W University Blvd, Melbourne, FL 32901-6975, USA}
\email{kperera@fit.edu}

\author[P. Winkert]{Patrick Winkert\orcidlink{0000-0003-0320-7026}}
\address[P. Winkert]{Technische Universit\"{a}t Berlin, Institut f\"{u}r Mathematik, Stra\ss{}e des 17.\,Juni 136, 10623 Berlin, Germany}
\email{winkert@math.tu-berlin.de}


\begin{document}

\begin{abstract}
We study radial solutions of the double phase equation $-\Delta_p u - \operatorname{div}(|x|^{-b}\,|\nabla u|^{q-2}\nabla u) = f(|x|,u)$ in $\mathbb{R}^N$, whose operator admits a scaling that acts on both phases simultaneously. Our first main result is a sharp description of the continuous and compact embeddings of the radial energy space into weighted Lebesgue spaces $L^s(\mathbb{R}^N;|x|^{-c})$ with an explicit characterization of their admissible regions in the $(s,c)$-plane. The two critical edges reflect the two phases: on the edge attached to the $p$-phase compactness is lost by spreading to infinity, while on the edge attached to the weighted $q$-phase it is lost by concentration at the origin. Our second main result concerns equations with one critical term on each edge. We show that bounded Palais-Smale sequences either converge or lose at least the energy $\beta^\ast = \min\{\beta_1,\beta_2\}$, where $\beta_1$ and $\beta_2$ are the explicit single-phase thresholds, and that this level is optimal when the critical terms are unweighted. Thus, the two mechanisms cannot combine to lower the compactness threshold. Combining these results with a Poho\v{z}aev identity, variational methods for scaled functionals and the Fadell-Rabinowitz cohomological index, we prove nonexistence results, obtain an unbounded sequence of eigenvalues on the scaling line and infinitely many solutions for pure power equations, and show that large superscaled and small subscaled perturbations of the doubly critical equation yield multiple solutions.
\end{abstract}

\maketitle

\section{Introduction}

Consider the double phase problem
\begin{equation} \label{1}
	- \operatorname{div} \left(|\nabla u|^{p-2}\, \nabla u + |x|^{-b}\, |\nabla u|^{q-2}\, \nabla u\right) = f(|x|,u) \quad \text{in } \mathbb{R}^N,
\end{equation}
where $1 < p < q < N$, $0 \le b < N - q$, and $f$ is a Carath\'{e}odory function on $[0,\infty) \times \R$ satisfying a suitable growth condition. An important feature of the operator
\begin{align}\label{special-double-phase}
	\A(u) := - \operatorname{div} \left(|\nabla u|^{p-2}\, \nabla u + |x|^{-b}\, |\nabla u|^{q-2}\,\nabla u\right)
\end{align}
is its scaling property. Let
\begin{equation} \label{100}
	\alpha = \frac{b}{q - p}
\end{equation}
and consider
\begin{align*}
	u_t(x) = t^{\alpha + 1}\, u(x/t), \quad t > 0.
\end{align*}
Then
\begin{equation} \label{3}
	\A(u_t)(x) = t^{\beta - 1}\, \A(u)(x/t),
\end{equation}
where
\begin{align*}
	\beta = \alpha\, (p - 1),
\end{align*}
since $\alpha\, (q - 1) - b = \alpha\, (p - 1)$. Thus, both gradient terms transform according to the same scaling law, a property that plays a central role in our approach.

The operator $\A$ contains no lower order terms, so the variational problem has no natural length scale. It is governed by the $p$-phase at infinity and by the weighted $q$-phase near the origin. This paper makes the two-phase structure precise and exploits it. Our two main contributions are:
\begin{enumerate}
	\item[\textnormal{(a)}]
		\emph{The sharp embedding triangle.} We determine precisely which weighted Lebesgue spaces $L^s(\R^N;|x|^{-c})$ contain the radial energy space, continuously and compactly (Theorem \ref{thm:intro-embedding}). The admissible set is a triangle in the $(s,c)$-plane with two critical edges. Each critical edge belongs to one phase and comes with its own mechanism of noncompactness: spreading to infinity for the $p$-phase and concentration at the origin for the $q$-phase.
	\item[\textnormal{(b)}]
		\emph{The threshold analysis.} For equations with a critical term on each edge, we identify the energy level below which Palais-Smale sequences are compact (Theorem \ref{thm:intro-threshold}). This level is the smaller of the two explicit single-phase thresholds, and we prove its optimality when the critical terms are unweighted. A noncompact sequence may a priori lose mass through both mechanisms at once. We show that any such loss costs at least as much as the cheaper pure loss, so simultaneous losses cannot lower the compactness threshold.
\end{enumerate}
We then use these results, together with the variational methods for scaled functionals recently developed by Mercuri and Perera \cite{Mercuri-Perera-2026} and Perera \cite{Perera-2026}, to establish nonexistence, existence, and multiplicity results for radial weak solutions of equation \eqref{1}.

\subsection{Background and related work}

The operator $\A$ in \eqref{special-double-phase} belongs to the class of double phase operators, which are generally of the form
\begin{align*}
	\A_\mu(u) = -\operatorname{div}\left( |\nabla u|^{p-2}\,\nabla u +\mu(x)\,|\nabla u|^{q-2}\,\nabla u \right),
\end{align*}
where $1<p<q$ and $\mu(\cdot)$ is a nonnegative coefficient. The corresponding energy density is characterized by the presence of two different growth exponents,
\begin{align*}
	\H(x,t)=t^p+\mu(x)t^q.
\end{align*}
The associated energy functional is given by
\begin{align*}
	I_\mu(u) = \int_{\R^N} \left( \frac{1}{p}\,|\nabla u|^p +\frac{\mu(x)}{q}\,|\nabla u|^q \right)\,\mathrm{d}x.
\end{align*}
Several familiar operators arise as special cases. If $\mu\equiv0$, one recovers the classical $p$-Laplacian, whereas a positive constant coefficient leads to a $(p,q)$-Laplacian. More generally, a spatially varying coefficient allows the two growth regimes to interact in a way that depends on both the position and the size of the gradient. At points where $\mu(\cdot)$ vanishes, the energy exhibits purely $p$-growth, while the $q$-term becomes increasingly important for large gradients wherever $\mu(\cdot)$ is positive.

In the present paper, we consider the particular choice $\mu(x)=|x|^{-b}$. For this choice of coefficient, the relative strength of the two phases is determined by
\begin{align*}
	\frac{|x|^{-b}\,t^q}{t^p} = |x|^{-b}\,t^{q-p}, \quad x\neq 0,\quad t>0.
\end{align*}
Thus, the two contributions are comparable when $t=|x|^\alpha$, where $\alpha$ is defined in \eqref{100}. For $t<|x|^\alpha$, the $p$-phase dominates, whereas for $t>|x|^\alpha$, the $q$-phase is dominant. This position-dependent transition between the two growth regimes reflects the particular structure of our operator and is closely related to its scaling property. For $b=0$, our operator reduces to the standard $(p,q)$-Laplacian. When $b>0$, however, the coefficient becomes singular at the origin and tends to zero as $|x|\to\infty$. This distinguishes our setting from the more commonly studied double phase models with bounded coefficients.

The origins of double phase functionals can be traced to the pioneering work of Zhikov \cite{Zhikov-1986}, particularly in connection with homogenization, elasticity, and the Lavrentiev phenomenon (see also \cite{Zhikov-1995,Zhikov-2011}). Such energies provide a mathematical description of heterogeneous materials whose constituents exhibit different growth behavior. From a broader perspective, double phase problems belong to the theory of variational integrals with nonstandard growth, to which Marcellini \cite{Marcellini-1989,Marcellini-1991} made fundamental contributions. Their regularity properties have since been investigated extensively. In particular, Colombo and Mingione \cite{Colombo-Mingione-2015a,Colombo-Mingione-2015b} and Baroni, Colombo, and Mingione \cite{Baroni-Colombo-Mingione-2015,Baroni-Colombo-Mingione-2016,Baroni-Colombo-Mingione-2018} developed important regularity results for double phase functionals and more general variational problems with nonuniform ellipticity. These results also reveal the significance of the regularity of the coefficient and the relationship between the growth exponents $p$ and $q$. A broader account of the theory of nonstandard growth and nonuniform ellipticity is given by Mingione and R\u{a}dulescu \cite{Mingione-Radulescu-2021}. Further contributions to the theory of double phase problems on bounded domains, including results on eigenvalue problems, existence and multiplicity, critical growth, and the underlying function spaces and operators, can be found in the works of Colasuonno and Perera \cite{Colasuonno-Perera-2025}, Colasuonno and Squassina \cite{Colasuonno-Squassina-2016}, Crespo-Blanco, Gasi\'{n}ski, Harjulehto, and Winkert \cite{Crespo-Blanco-Gasinski-Harjulehto-Winkert-2022}, Ho and Winkert \cite{Ho-Winkert-2023}, Liu and Dai \cite{Liu-Dai-2018}, and Perera and Squassina \cite{Perera-Squassina-2018}.

We next turn to double phase problems in the entire space $\R^N$. In contrast to problems on bounded domains, the lack of compactness of Sobolev embeddings introduces additional difficulties in the application of variational methods. Consequently, suitable assumptions on the coefficients, potential terms, or nonlinearities, as well as symmetry conditions, play an important role in obtaining existence and multiplicity results. An early contribution in this direction was made by Liu and Dai \cite{Liu-Dai-2020}, who studied the problem
\begin{align*}
	-\operatorname{div}\left( |\nabla u|^{p-2}\,\nabla u +a(x)\,|\nabla u|^{q-2}\,\nabla u \right)+|u|^{p-2}\,u+a(x)\,|u|^{q-2}\,u =f(x,u) \quad\text{in }\R^N,
\end{align*}
where $a(\cdot)$ is a nonnegative H\"{o}lder continuous function and the reaction term satisfies suitable subcritical growth conditions. Combining variational methods with the Nehari manifold approach, topological degree theory, and deformation arguments, they established the existence of positive and negative ground state solutions, together with a least energy sign-changing solution. Under additional symmetry assumptions, they also obtained infinitely many pairs of solutions. Of particular relevance to our work is their treatment of radially symmetric coefficients and nonlinearities, which leads to corresponding existence and multiplicity results for radial solutions.

Another class of double phase problems in $\R^N$ was investigated by Liu and Winkert \cite{Liu-Winkert-2022}, who considered the combined effects of singular and superlinear nonlinearities in equations of the form
\begin{align*}
	&-\operatorname{div}\left( |\nabla u|^{p-2}\,\nabla u +\mu(x)\,|\nabla u|^{q-2}\,\nabla u \right)+u^{p-1}+\mu(x)\,u^{q-1}\\
	&=\xi(x)\,u^{-\delta}+\lambda u^{\nu-1} \quad\text{in }\R^N,
\end{align*}
where $0<\delta<1$, $q<\nu<p^\ast$, and the coefficient $\mu(\cdot)$ is nonnegative and bounded. Using a splitting of the Nehari manifold, they proved the existence of at least two positive weak solutions for sufficiently small positive values of $\lambda$, one with negative energy and the other with positive energy. Here, the singularity appears in the reaction term, whereas in our problem the singular behavior for $b>0$ is caused by the coefficient of the $q$-gradient term.

Radial solutions have also been investigated by Lian, Ge, Cao, and Zhou \cite{Lian-Ge-Cao-Zhou-2026}. They considered
\begin{align*}
	-\operatorname{div}\left( |\nabla u|^{p-2}\,\nabla u +\mu(x)\,|\nabla u|^{q-2}\,\nabla u \right)+|u|^{p-2}\,u+\mu(x)|u|^{q-2}\,u =f(x,u) \quad\text{in }\R^N,
\end{align*}
with a nonnegative, integrable, radially symmetric coefficient $\mu(\cdot)$ and a radial nonlinearity $f$. By working on a suitable nodal Nehari set, they established the existence of a least energy sign-changing radial solution. Under an additional condition on the nonlinearity, this solution has precisely two nodal domains. Their analysis extends the radial sign-changing results of Liu and Dai \cite{Liu-Dai-2020} by allowing weaker monotonicity assumptions on the reaction term.

Further existence and multiplicity results for double phase equations on unbounded domains have been obtained through various variational approaches. Steg\-li\'{n}ski \cite{Steglinski-2022} investigated problems with an unbounded potential and nonlinearities that need not satisfy the Ambrosetti-Rabinowitz condition. Using the Fountain and Dual Fountain Theorems, together with a suitable compactness condition, he proved the existence of infinitely many solutions. Ge and Pucci \cite{Ge-Pucci-2022} employed a perturbation method to obtain nontrivial solutions for double phase equations with superlinear reactions, again without requiring the Ambrosetti-Rabinowitz condition. Ambrosio and Essebei \cite{Ambrosio-Essebei-2023} established multiplicity results for problems involving two weighted gradient terms by combining Ricceri's variational principle with truncation techniques and Moser iteration. In addition, Ge and Yuan \cite{Ge-Yuan-2024} studied parameter-dependent double phase problems with sign-changing potentials and obtained the existence of at least two nontrivial solutions under appropriate assumptions.

The influence of the spatial coefficient on the number and behavior of solutions was investigated more recently by Liu and Winkert \cite{Liu-Winkert-2026}. They considered
\begin{align*}
	-\operatorname{div}\left( |\nabla u|^{p-2}\,\nabla u +a(\varepsilon x)\,|\nabla u|^{q-2}\,\nabla u \right)+|u|^{p-2}\,u +a(\varepsilon x)\,|u|^{q-2}\,u =f(u) \quad\text{in }\R^N,
\end{align*}
where $\varepsilon>0$ and $a(\cdot)$ is a nonnegative, bounded coefficient. By applying genus theory, the Nehari manifold method, and the Lusternik-Schnirelmann category, they derived multiplicity results for constant-sign and sign-changing solutions. In particular, the number of solutions is related to topological properties of the zero set of the coefficient, and concentration phenomena are obtained as $\varepsilon\to0$. This demonstrates how the spatial dependence of the double phase coefficient can substantially influence the solution structure.

Critical growth in double phase equations on $\R^N$ was studied by Ha and Ho \cite{Ha-Ho-2025}. They investigated equations of Schr\"odinger type with variable exponents, including problems of the form
\begin{align*}
	&-\operatorname{div}\left( |\nabla u|^{p(x)-2}\,\nabla u +a(x)\,|\nabla u|^{q(x)-2}\,\nabla u \right)\\
	&+\lambda V(x)\left( |u|^{p(x)-2}\,u +a(x)|u|^{q(x)-2}\,u \right)\\
	&=f(x,u)+\theta B(x,u) \quad\text{in }\R^N,
\end{align*}
where $B$ represents a nonlinearity with generalized critical growth. They developed a concentration-compactness principle for the associated Musielak-Orlicz-Sobolev spaces, including a version that accounts for concentration at infinity. These results were then applied to establish existence and concentration properties of solutions under suitable assumptions on the potential and nonlinearities.

Other recent contributions address more general growth structures and nonlinearities. Bahrouni, Bahrouni, and Winkert \cite{Bahrouni-Bahrouni-Winkert-2025} investigated double phase problems with variable exponents depending on the solution and its gradient. They established continuous and compact embeddings for the corresponding function spaces and applied abstract critical point theory to obtain existence and multiplicity results in $\R^N$. Bahrouni, Fiscella, and Winkert \cite{Bahrouni-Fiscella-Winkert-2026} considered critical double phase equations involving logarithmic nonlinearities and proved an existence result by variational methods. Furthermore, Ha, Ho, Quan, and Sim \cite{Ha-Ho-Quan-Sim-2026} established global boundedness results for generalized Schr\"{o}dinger-type double phase problems with variable exponents, using De Giorgi iteration techniques. Their regularity estimates also yield applications to problems with supercritical growth.

The preceding studies illustrate the variety of analytical approaches to double phase equations in $\R^N$. Many of these works rely on additional lower-order terms, suitable potential functions, or regularity assumptions on the spatial coefficients to overcome the lack of compactness. In contrast, our problem \eqref{1} involves no lower-order terms on the left-hand side and is driven by the particular coefficient $|x|^{-b}$, which becomes singular at the origin when $b>0$. This setting leads to a homogeneous variational framework in which both gradient terms admit a common scaling. Together with radial symmetry, this structure motivates the study of sharp weighted embeddings and provides the basis for our existence and multiplicity results, including problems with critical nonlinearities.

\subsection{Functional setting and the embedding triangle}\label{intro-functional-setting}

The natural space to look for solutions to equation \eqref{1} is the intersection space
\begin{align*}
	W = \D^{1,p}(\R^N) \cap \D^{1,q}(\R^N;|x|^{-b}),
\end{align*}
defined as the completion of $C^\infty_c(\R^N)$ with respect to the norm
\begin{align*}
	\norm{u} = \norm[\D^{1,p}(\R^N)]{u} + \norm[\D^{1,q}(\R^N;|x|^{-b})]{u}.
\end{align*}
Here, $\D^{1,q}(\R^N;|x|^{-b})$ denotes the weighted homogeneous Sobolev space equipped with the norm
\begin{align*}
	\norm[\D^{1,q}(\R^N;|x|^{-b})]{u} := \left(\int_{\R^N} |x|^{-b}\, |\nabla u|^q\, \mathrm{d}x\right)^{1/q}.
\end{align*}
Throughout the paper, we work in the subspace $W_\rad$ of radial functions in $W$.

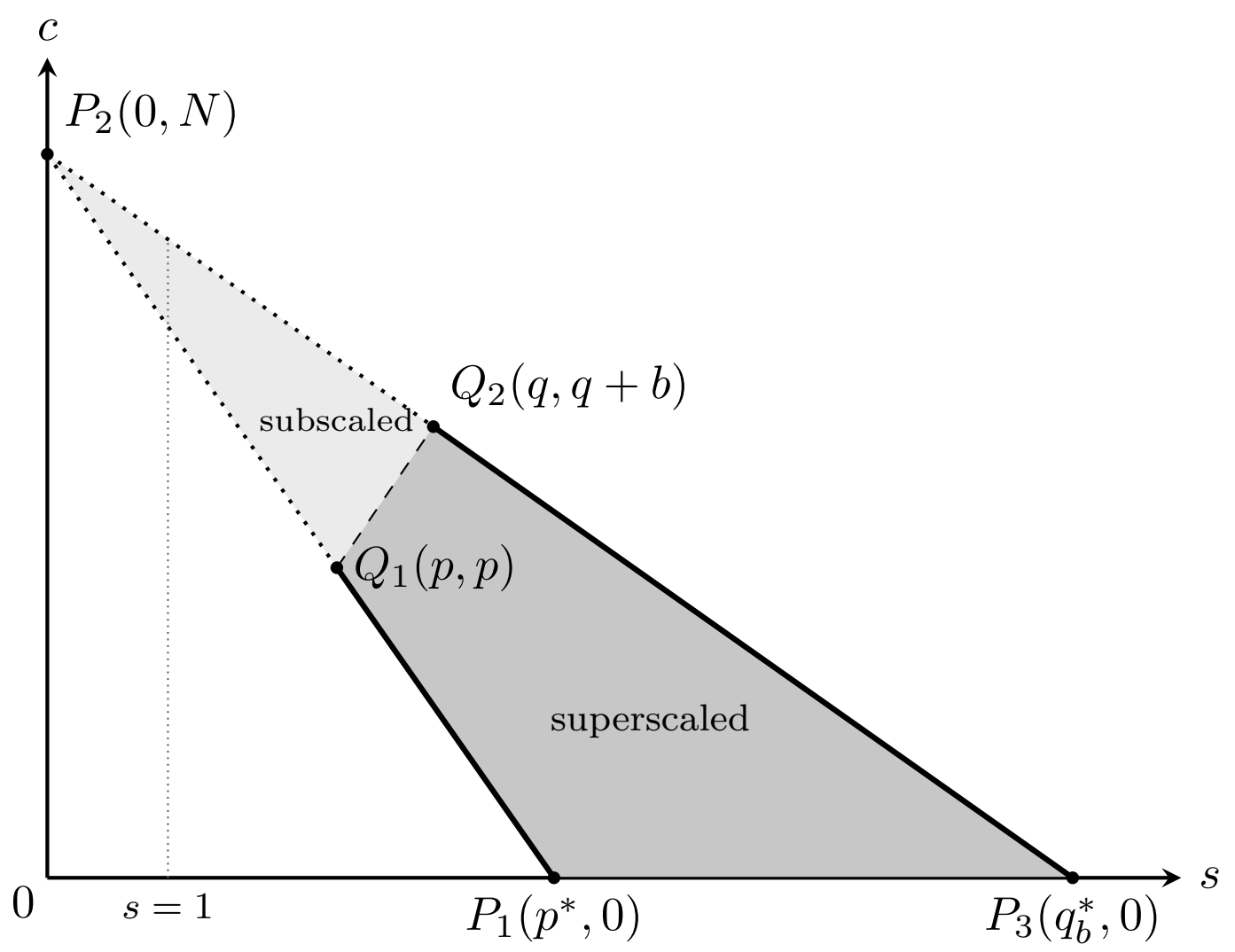
\begin{figure}[ht]
	\centering
	\begin{tikzpicture}[scale=0.9]
		\draw[-stealth, thick] (0,0) -- (9.4,0) node[right] {$s$};
		\draw[-stealth, thick] (0,0) -- (0,6.8) node[above] {$c$};
		\node at (-0.2,-0.2) {0};
		\coordinate (P1) at (4.2,0);
		\coordinate (P2) at (0,6);
		\coordinate (P3) at (8.5,0);
		\coordinate (Q1) at (2.4,2.571);
		\coordinate (Q2) at (3.2,3.741);
		\fill[black!8] (P2) -- (Q1) -- (Q2) -- cycle;
		\fill[black!22] (Q1) -- (P1) -- (P3) -- (Q2) -- cycle;
		\draw[thick, dotted] (P2) -- (Q1);
		\draw[thick, dotted] (P2) -- (Q2);
		\draw[very thick] (Q1) -- (P1);
		\draw[very thick] (Q2) -- (P3);
		\draw (P1) -- (P3);
		\draw[dashed] (Q1) -- (Q2);
		\draw[densely dotted, black!60] (1,0) -- (1,5.29);
		\node[below] at (1,0) {\footnotesize $s = 1$};
		\fill (P1) circle (1.5pt) node[below] {$P_1(p^\ast,0)$};
		\fill (P2) circle (1.5pt) node[above right] {$P_2(0,N)$};
		\fill (P3) circle (1.5pt) node[below] {$P_3(q^\ast_b,0)$};
		\fill (Q1) circle (1.5pt) node[right] {$Q_1(p,p)$};
		\fill (Q2) circle (1.5pt) node[above right] {$Q_2(q,q + b)$};
		\node at (2.4,3.8) {\scriptsize subscaled};
		\node at (5.0,1.3) {\footnotesize superscaled};
	\end{tikzpicture}
	\caption{The triangle $\T$ in the $(s,c)$-plane. The dashed segment $Q_1Q_2$ is the scaling line separating the subscaled region (light) from the superscaled region (dark). The thick segments $[Q_1,P_1]$ and $[Q_2,P_3]$ carry the critical exponents. The dotted segments $[P_2,Q_1)$ and $[P_2,Q_2)$ are excluded from the continuous embedding range, and compactness holds precisely for $(s,c)\in\mathcal{T}_{\mathrm{c}}$.}
	\label{fig:triangle}
\end{figure}

For $s \ge 1$ and $c \ge 0$, let $L^s(\R^N;|x|^{-c})$ be the weighted Lebesgue space with norm
\begin{align*}
	\norm[L^s(\R^N;|x|^{-c})]{u} := \left(\int_{\R^N} |x|^{-c} \,|u|^s\, \mathrm{d}x\right)^{1/s}.
\end{align*}
Let
\begin{align*}
	p^\ast = \frac{Np}{N - p}, \quad q^\ast_b = \frac{Nq}{N - q - b}.
\end{align*}
Radial functions in $W$ satisfy the two Strauss-type pointwise estimates
\begin{align*}
	|u(x)| \le C\, |x|^{-N/p^\ast} \norm[\D^{1,p}(\R^N)]{u}, \quad |u(x)| \le C \,|x|^{-N/q^\ast_b} \norm[\D^{1,q}(\R^N;|x|^{-b})]{u},
\end{align*}
for $x\neq 0$, with a suitable radial representative of $u$ (see Lemma \ref{Lemma Strauss}). The first estimate controls radial functions at infinity, while the second provides control near the origin. Together, they allow us to characterize the range of weighted Lebesgue spaces into which $W_\rad$ embeds. More precisely, this range is described by the triangle $\T$ in the $(s,c)$-plane with the vertices
\begin{align*}
	P_1(p^\ast,0), \quad P_2(0,N), \quad P_3(q^\ast_b,0)
\end{align*}
(see Figure \ref{fig:triangle}). The edges $P_1 P_2$ and $P_2 P_3$ lie on the lines
\begin{align*}
	c = N \left(1 - \frac{s}{p^\ast}\right) \quad \text{and} \quad c = N \left(1 - \frac{s}{q^\ast_b}\right),
\end{align*}
respectively, and they contain the points
\begin{align*}
	Q_1(p,p) \in P_1 P_2, \qquad Q_2(q,q + b) \in P_2 P_3,
\end{align*}
which correspond to the Hardy inequality in $\D^{1,p}(\R^N)$ and the weighted Hardy inequality in $\D^{1,q}(\R^N;|x|^{-b})$, respectively. For two points $A$ and $B$, we write $[A,B]$ for the closed segment joining them and $[A,B)$ for the same segment with the endpoint $B$ removed. Define
\begin{align*}
	\Tc & = \set{(s,c) \colon s \ge 1,\ c \ge 0,\ N\, (1 - s/p^\ast) < c < N\, (1 - s/q^\ast_b)},\\
	\Te & = \Tc \cup [Q_1,P_1] \cup [Q_2,P_3].
\end{align*}
Thus, $\Tc$ consists of those points of $\T$ with $s \ge 1$ that do not lie on the edges $P_1 P_2$ and $P_2 P_3$. The set $\Te$ is obtained from $\set{(s,c) \in \T \colon s \ge 1}$ by removing the half-open segments $[P_2,Q_1)$ and $[P_2,Q_2)$.

Our first main result is the following sharp embedding theorem.

\begin{theorem} \label{thm:intro-embedding}
	Let $s \ge 1$ and $c \ge 0$.
	\begin{enumerate}
		\item[\textnormal{(i)}]
			$W_\rad$ is continuously embedded in $L^s(\R^N;|x|^{-c})$ if and only if $(s,c) \in \Te$, and $W_\rad \not\subset L^s(\R^N;|x|^{-c})$ if $(s,c) \notin \Te$.
		\item[\textnormal{(ii)}]
			$W_\rad$ is compactly embedded in $L^s(\R^N;|x|^{-c})$ if and only if $(s,c) \in \Tc$.
		\item[\textnormal{(iii)}]
			On the critical edges, compactness fails through the two phases separately. If $(s,c) \in [Q_1,P_1]$, the dilations $u_\eps(x) = \eps^{N/p^\ast} u(\eps x)$, $\eps \searrow 0$, of a fixed radial $u \in C^\infty_c(\R^N \setminus \set{0}) \setminus \set{0}$ spread to infinity, converge weakly to $0$ in $W_\rad$, and preserve the $L^s(\R^N;|x|^{-c})$-norm, while their $\D^{1,q}(\R^N;|x|^{-b})$-norms tend to zero. If $(s,c) \in [Q_2,P_3]$, the dilations $u_\eps(x) = \eps^{-N/q^\ast_b} u(x/\eps)$, $\eps \searrow 0$, concentrate at the origin and have the same properties with the roles of the two phases exchanged.
	\end{enumerate}
\end{theorem}

This is Theorem \ref{Theorem embedding}. The two Strauss-type estimates explain the shape of $\T$. The edge $P_1P_2$ is where the decay of the $p$-phase at infinity becomes borderline, and the edge $P_2P_3$ is where the singularity allowed by the $q$-phase at the origin becomes borderline. The points $Q_1$ and $Q_2$ given by the two Hardy inequalities cut off the parts of these edges where integrability fails at the opposite end. In particular, $W_\rad$ is embedded in $L^s(\R^N)$ continuously for all $s \in [p^\ast,q^\ast_b]$ and compactly for $s \in (p^\ast,q^\ast_b)$. The restriction to radial functions is essential here when $b > 0$. The whole space $W$ is not embedded in $L^{q^\ast_b}(\R^N)$, and the Caffarelli-Kohn-Nirenberg inequalities give embeddings of $W$ only on a smaller quadrilateral (see Remark \ref{Remark 1}).

Theorem \ref{thm:intro-embedding} is a double phase counterpart of the classical radial compactness results. These go back to Strauss \cite{Strauss-1977}, who obtained compactness of radial $H^1(\R^N)$ in $L^s(\R^N)$ for $2 < s < 2^\ast$. In single phase problems, the radial lemma gives decay at infinity, and the lack of compactness at the Sobolev exponent comes from concentration alone. For radial weighted spaces with potentials, Su, Wang, and Willem \cite{Su-Wang-Willem-2007a,Su-Wang-Willem-2007b} described the compact embeddings in terms of the behavior of the weights at the origin and at infinity. In the present homogeneous double phase setting, the two phases play the roles of these two ends. The critical exponents $p^\ast$ and $q^\ast_b$ are distinct, and the corresponding losses of compactness occur at opposite ends of $(0,\infty)$. This is what makes the threshold analysis of Theorem \ref{thm:intro-threshold} below possible.

In view of these embeddings, a natural growth condition on the nonlinearity $f$ in equation \eqref{1} is
\begin{equation} \label{6}
	|f(|x|,t)| \le \sum_{i=1}^n a_i\, |x|^{-c_i}\, |t|^{s_i-1} \quad \text{for a.e.\,}x \in \R^N \text{ and all }t\in \R,
\end{equation}
where $(s_i,c_i) \in \Te$ and $a_i > 0$ are constants. Under this growth condition, a weak radial solution of equation \eqref{1} is a function $u \in W_\rad$ satisfying
\begin{align*}
	\int_{\R^N} \left(|\nabla u|^{p-2} + |x|^{-b}\, |\nabla u|^{q-2}\right) \nabla u \cdot \nabla v\, \mathrm{d}x = \int_{\R^N} f(|x|,u) \,v \,\mathrm{d}x \quad \text{for all }v \in W_\rad.
\end{align*}
Such solutions are precisely the critical points of the $C^1$-functional
\begin{align*}
	E(u) = \frac{1}{p} \int_{\R^N} |\nabla u|^p\, \mathrm{d}x + \frac{1}{q} \int_{\R^N} |x|^{-b}\, |\nabla u|^q\, \mathrm{d}x - \int_{\R^N} F(|x|,u)\, \mathrm{d}x, \quad u \in W_\rad,
\end{align*}
where $F(|x|,t) = \int_0^t f(|x|,\tau)\,\mathrm{d}\tau$ is a primitive of $f$ with respect to its second variable. By a symmetric criticality argument, every weak radial solution also satisfies the weak formulation of equation \eqref{1} for all test functions in $W$ (see Proposition \ref{Proposition symmetric criticality}), although the energy functional $E$ need not be defined on the whole space $W$.

\subsection{Scaling and nonlinear regimes}

As in Mercuri and Perera \cite{Mercuri-Perera-2026}, the scaling property \eqref{3} of the operator $\A$ leads to a natural classification of nonlinearities $f$ according to their behavior under scaling. More precisely, we distinguish the following three regimes:
\begin{enumerate}
	\item[\textnormal{(i)}]
		{\bf Subscaled:} The nonlinearity $f$ is subscaled if
		\begin{align*}
			\lim_{r \to \infty}\,\frac{f(r,\pm r^{\alpha + 1})}{\pm r^{\beta - 1}} = 0.
		\end{align*}
	\item[\textnormal{(ii)}]
		{\bf Asymptotically scaled:} The nonlinearity $f$ is asymptotically scaled if
		\begin{align*}
			\lim_{r \to \infty} \,\frac{f(r,\pm r^{\alpha + 1})}{\pm r^{\beta - 1}} = \lambda \in (0,\infty).
		\end{align*}
	\item[\textnormal{(iii)}]
		{\bf Superscaled:} The nonlinearity $f$ is superscaled if
		\begin{align*}
			\lim_{r \to \infty} \,\frac{f(r,\pm r^{\alpha + 1})}{\pm r^{\beta - 1}} = \infty.
		\end{align*}
\end{enumerate}
Here, $r = |x|$ represents the radius, while $\pm r^{\alpha + 1}$ corresponds to the amplitude of the scaled function $u_r$. In each case, the limit is required to hold for both choices of sign. To illustrate this classification, consider the weighted power nonlinearity
\begin{align*}
	f(|x|,t) = \lambda\, |x|^{-c}\, |t|^{s-2}\, t, \quad\lambda>0, \quad s>1.
\end{align*}
A direct computation shows that
\begin{align*}
	\frac{f(r,\pm r^{\alpha + 1})}{\pm r^{\beta - 1}} = \lambda\, r^{(\alpha + 1) \,s - \alpha p - c} = \lambda\, r^{\rho - \sigma},
\end{align*}
where
\begin{equation} \label{rho sigma}
	\sigma = \alpha p + N = \alpha q + N - b, \quad \rho = \rho(s,c) = (\alpha + 1) \,s - c + N.
\end{equation}
Thus, the three regimes are determined by the sign of $\rho-\sigma$. The relation $\rho=\sigma$, or equivalently, $c = (\alpha + 1)\, s - \alpha p$, defines a line passing through $Q_1$ and $Q_2$. We call the segment $Q_1 Q_2$ the scaling line. It divides $\T$ into the subscaled region above it, which is the triangle with the vertices $P_2$, $Q_1$, and $Q_2$, and the superscaled region below it, which is the quadrilateral with the vertices $Q_1$, $P_1$, $P_3$, and $Q_2$ (see Figure \ref{fig:triangle}). The scaling also explains the role of the critical segments. In fact,
\begin{equation} \label{crit-segments}
	\rho = \frac{\sigma s}{p} \quad \text{on } [Q_1,P_1], \quad \rho = \frac{\sigma s}{q} \quad \text{on } [Q_2,P_3].
\end{equation}
These relations will be useful in distinguishing the various existence and multiplicity results below.

\subsection{Main contributions and results}

We now present several representative results that illustrate the main contributions of this paper. They are special cases or direct consequences of the general results proved in Section \ref{Section results}. Throughout, a weak radial solution is understood to be a critical point in $W_\rad$ of the corresponding energy functional $E$, as described in Subsection \ref{intro-functional-setting}.

Our approach combines the sharp weighted embeddings of radial functions with the scaling properties of the operator and the variational theory of scaled functionals. This allows us to treat eigenvalue problems on the scaling line as well as equations with nonlinearities lying above or below it. For critical problems, the loss of compactness requires additional arguments based on a Poho\v{z}aev identity and suitable energy thresholds.

\subsubsection*{Scaled eigenvalue problem}

Let $(s,c)$ lie on the open scaling line $(Q_1,Q_2)$, i.e., $p < s < q$ and $c = (\alpha + 1) \,s - \alpha p$, and consider the eigenvalue problem
\begin{equation} \label{intro:eigenvalue}
	- \operatorname{div} \left(|\nabla u|^{p-2}\, \nabla u + |x|^{-b}\, |\nabla u|^{q-2}\, \nabla u\right) = \lambda\, |x|^{-c} \,|u|^{s-2} u \quad \text{in } \R^N.
\end{equation}
Since $\rho(s,c) = \sigma$, this problem is invariant under the scaling $u \mapsto u_t$. In other words, whenever $u$ is a solution, so is $u_t$ for every $t > 0$. We call $\lambda \in \R$ an eigenvalue if problem \eqref{intro:eigenvalue} admits a nontrivial weak radial solution.

\begin{theorem} \label{thm:intro-eigenvalues}
	Assume that $p < s < q$ and $c = (\alpha + 1)\, s - \alpha p$. Then problem \eqref{intro:eigenvalue} admits an unbounded sequence of positive eigenvalues $\lambda_k(s) \nearrow \infty$ with radial eigenfunctions. Moreover, the first eigenvalue is characterized by
	\begin{align*}
		\lambda_1(s) = \min_{u \in W_\rad \setminus \set{0}}\, \frac{\ds{\frac{1}{p} \int_{\R^N} |\nabla u|^p\, \mathrm{d}x + \frac{1}{q} \int_{\R^N} |x|^{-b}\, |\nabla u|^q\, \mathrm{d}x}}{\ds{\frac{1}{s} \int_{\R^N} |x|^{-c}\, |u|^s\, \mathrm{d}x}}.
	\end{align*}
\end{theorem}

The existence of the eigenvalue sequence follows from Theorem \ref{Theorem 1}. To obtain the formula for $\lambda_1(s)$, we apply Theorem \ref{Theorem 3} \textnormal{(\hyperlink{Theorem 3.i}{i})}. Indeed, $I_\sigma$ and $J_\sigma$, introduced in Section \ref{Section abstract}, have the same scaling behavior, and hence
\begin{align*}
	\widetilde{\Psi}(\pi(u)) = \frac{I_\sigma(u)}{J_\sigma(u)}\quad \text{for all } u \ne 0.
\end{align*}

\subsubsection*{Superscaled pure power equation}

Next, we consider the unweighted pure power equation
\begin{equation} \label{intro:pure-power}
	- \operatorname{div} \left(|\nabla u|^{p-2}\, \nabla u + |x|^{-b}\, |\nabla u|^{q-2}\, \nabla u\right) = \lambda |u|^{s-2}\, u \quad \text{in } \R^N,
\end{equation}
with associated energy functional
\begin{align*}
	E(u) = \frac{1}{p} \int_{\R^N} |\nabla u|^p\, \mathrm{d}x + \frac{1}{q} \int_{\R^N} |x|^{-b}\, |\nabla u|^q\, \mathrm{d}x - \frac{\lambda}{s} \int_{\R^N} |u|^s\, \mathrm{d}x, \quad u \in W_\rad.
\end{align*}
Observe that $p^\ast < q^\ast_b$, so the interval $(p^\ast,q^\ast_b)$ is nonempty. Moreover, for every $s$ in this interval, the point $(s,0)$ lies in the superscaled region.

\begin{theorem} \label{thm:intro-pure-power}
	Let $\lambda > 0$ and $p^\ast < s < q^\ast_b$. Then equation \eqref{intro:pure-power} has an infinite sequence of nontrivial weak radial solutions $\seq{u_k}$ such that $E(u_k) > 0$ and $E(u_k) \nearrow \infty$.
\end{theorem}

This is Corollary \ref{cor:pure-power}, which follows as a special case of Theorem \ref{Theorem 2}.

\subsubsection*{Critical equations}

We first consider the pure critical power equations
\begin{equation} \label{intro:p-ast}
	- \operatorname{div} \left(|\nabla u|^{p-2}\, \nabla u + |x|^{-b}\, |\nabla u|^{q-2} \,\nabla u\right) = |u|^{p^\ast - 2} \,u \quad \text{in } \R^N
\end{equation}
and
\begin{equation} \label{intro:q-ast-b}
	- \operatorname{div} \left(|\nabla u|^{p-2}\, \nabla u + |x|^{-b}\, |\nabla u|^{q-2} \,\nabla u\right) = |u|^{q^\ast_b - 2}\, u \quad \text{in } \R^N.
\end{equation}
In both cases, a Poho\v{z}aev identity rules out the existence of nontrivial weak radial solutions.

\begin{theorem} \label{thm:intro-nonexistence}
	Equations \eqref{intro:p-ast} and \eqref{intro:q-ast-b} have no nontrivial weak radial solutions.
\end{theorem}

This is Corollary \ref{cor:nonexistence}, which follows from the more general nonexistence result in Theorem \ref{Theorem nonexistence}.

\begin{remark} 
	When both critical terms are present, we arrive at the equation
	\begin{equation} \label{intro:both-critical}
		- \operatorname{div} \left(|\nabla u|^{p-2} \,\nabla u + |x|^{-b}\, |\nabla u|^{q-2} \,\nabla u\right) = |u|^{p^\ast - 2} \,u + |u|^{q^\ast_b - 2}\,u \quad \text{in } \R^N.
	\end{equation}
	In this situation, the Poho\v{z}aev identity only implies
	\begin{align*}
		\int_{\R^N} |\nabla u|^p\, \mathrm{d}x = \int_{\R^N} |u|^{p^\ast}\, \mathrm{d}x, \quad \int_{\R^N} |x|^{-b} \,|\nabla u|^q\, \mathrm{d}x = \int_{\R^N} |u|^{q^\ast_b}\, \mathrm{d}x
	\end{align*}
	for every weak radial solution $u$ (see Lemma \ref{Lemma Pohozaev critical}). These identities do not, however, exclude nontrivial solutions. Whether equation \eqref{intro:both-critical} admits a nontrivial weak radial solution remains an open problem (see also Remark \ref{Remark lambda zero}). This corresponds to the special case $\lambda = 0$, $\mu_1 = \mu_2 = 1$, and $c_1 = c_2 = 0$ of the general critical equation \eqref{19} considered in Section \ref{Section results}.
\end{remark}

We next consider perturbations of equation \eqref{intro:both-critical} of the form
\begin{equation} \label{intro:critical}
	\begin{aligned}
		&- \operatorname{div} \left(|\nabla u|^{p-2} \,\nabla u + |x|^{-b}\, |\nabla u|^{q-2} \,\nabla u\right)\\
		&= \lambda\, |x|^{-c}\, |u|^{s-2} u + |u|^{p^\ast - 2}\, u + |u|^{q^\ast_b - 2}\, u \quad \text{in } \R^N,
	\end{aligned}
\end{equation}
where $\lambda > 0$. The corresponding energy functional is
\begin{align*}
	E(u)
	&= \frac{1}{p} \int_{\R^N} |\nabla u|^p\, \mathrm{d}x + \frac{1}{q} \int_{\R^N} |x|^{-b}\, |\nabla u|^q\, \mathrm{d}x - \frac{\lambda}{s} \int_{\R^N} |x|^{-c}\, |u|^s\, \mathrm{d}x\\
	&\quad - \frac{1}{p^\ast} \int_{\R^N} |u|^{p^\ast}\, \mathrm{d}x - \frac{1}{q^\ast_b} \int_{\R^N} |u|^{q^\ast_b}\, \mathrm{d}x, \quad u \in W_\rad.
\end{align*}
Equation \eqref{intro:critical} has a critical term on each critical edge of $\T$, namely at the vertices $P_1$ and $P_3$. Let
\begin{align*}
	S = \inf_{u \in W_\rad \setminus \set{0}}\, \frac{\ds{\int_{\R^N} |\nabla u|^p\, \mathrm{d}x}}{\ds{\left(\int_{\R^N} |u|^{p^\ast}\, \mathrm{d}x\right)^{p/p^\ast}}}, \quad S_b = \inf_{u \in W_\rad \setminus \set{0}}\, \frac{\ds{\int_{\R^N} |x|^{-b}\, |\nabla u|^q\, \mathrm{d}x}}{\ds{\left(\int_{\R^N} |u|^{q^\ast_b}\, \mathrm{d}x\right)^{q/q^\ast_b}}}.
\end{align*}
Here $S$ is the Sobolev constant of $\D^{1,p}(\R^N)$, attained by the Aubin-Talenti functions (see Aubin \cite{Aubin-1976} and Talenti \cite{Talenti-1976}), while $S_b$ is the sharp radial weighted Sobolev constant, whose infimum over $W_\rad$ can be approached by truncations of a Bliss-type extremal (see Bliss \cite{Bliss-1930} and Proposition \ref{Proposition sharp}). The single-phase thresholds are
\begin{align*}
	\beta_1 = \frac{1}{N}\, S^{N/p}, \quad \beta_2 = \frac{q + b}{Nq}\, S_b^{N/(q + b)}.
\end{align*}
The first is the familiar compactness level for the $p$-Laplacian with critical growth, and the second is its counterpart for the weighted $q$-phase. Our second main result shows that the compactness threshold for the doubly critical problem is the smaller of the two.

\begin{theorem} \label{thm:intro-threshold}
	Assume $q < p^\ast$, $s > 1$, $(s,c) \in \Tc$, and $\lambda \in \R$. Let
	\begin{align*}
		\beta^\ast = \min \set{\beta_1,\beta_2}.
	\end{align*}
	\begin{enumerate}
		\item[\textnormal{(i)}]
			If $\seq{u_j}$ is a bounded \PS{c} sequence of $E$, then a subsequence converges weakly to a weak radial solution $u$ of equation \eqref{intro:critical}, and either it converges strongly or $c \ge E(u) + \beta^\ast$.
		\item[\textnormal{(ii)}]
			$E$ has bounded \PS{\beta_1} and \PS{\beta_2} sequences that converge weakly to $0$ and have no convergent subsequences. The first loses compactness by spreading to infinity in the $p$-phase, and the second by concentration at the origin in the $q$-phase.
		\item[\textnormal{(iii)}]
			If $\lambda > 0$ and $\rho(s,c) \ge \sigma$, then $E$ satisfies the \PS{c} condition for every $c < \beta^\ast$ but not for $c = \beta^\ast$.
	\end{enumerate}
\end{theorem}

Part \textnormal{(i)} is the special case $\mu_1 = \mu_2 = 1$, $c_1 = c_2 = 0$ of Proposition \ref{Proposition PS}, and part \textnormal{(ii)} is Proposition \ref{Proposition sharp}. For part \textnormal{(iii)}, every \PS{c} sequence is bounded by Lemma \ref{Lemma bounded}, and $E(u) \ge 0$ for every weak radial solution $u$ by \eqref{energy critical}.

The structure of Theorem \ref{thm:intro-threshold} is worth noting. A noncompact Palais-Smale sequence could, a priori, lose mass through both mechanisms at the same time, partly spreading to infinity in the $p$-phase and partly concentrating at the origin in the $q$-phase. One might expect such mixed losses to be cheaper than either pure loss, or the threshold to depend on how the two critical terms interact. Neither happens. The key algebraic estimate (Lemma \ref{Lemma threshold}) shows that every mixed profile costs at least $\min\set{\beta_1,\beta_2}$, and that this bound is attained by the pure profiles. Thus, the interaction between the two critical terms does not lower the compactness threshold. In contrast, a nontrivial solution of the unperturbed equation \eqref{intro:both-critical} would have energy at least $\beta_1 + \beta_2$ (Remark \ref{Remark lambda zero}). The threshold is used in the minimax arguments behind Theorems \ref{thm:intro-large-lambda} and \ref{thm:intro-concave-convex}.

The existence and multiplicity of solutions depend on the scaling regime of the perturbation. We first consider the superscaled case, where a sufficiently large parameter $\lambda$ leads to multiple solutions with positive energy.

\begin{theorem} \label{thm:intro-large-lambda}
	Assume $q < p^\ast$, $c = 0$, and $p^\ast < s < q^\ast_b$. Then, for any $m \ge 1$, there exists $\lambda^\ast_m > 0$ such that equation \eqref{intro:critical} has $m$ distinct pairs of nontrivial weak radial solutions at positive energy levels whenever $\lambda > \lambda^\ast_m$. In particular, the number of solutions tends to infinity as $\lambda \to \infty$.
\end{theorem}

This is the special case $\mu_1 = \mu_2 = 1$ and $c_1 = c_2 = c = 0$ of Theorem \ref{Theorem large lambda}. Indeed, assumptions \textnormal{(\hyperlink{H1}{H$_1$})} and \textnormal{(\hyperlink{H2}{H$_2$})}, together with $\rho > \sigma$, are satisfied because $(s,0)$ lies on the open base of $\mathcal{T}$, while $P_1\in[Q_1,P_1]$ and $P_3\in[Q_2,P_3]$. Moreover, \textnormal{(\hyperlink{H4}{H$_4$})} follows from $q < p^\ast < q^\ast_b$ and condition \textnormal{(\hyperlink{B}{B})} holds as $s > p^\ast > q$.

We turn to the subscaled case, for which a small positive parameter $\lambda$ leads to a different multiplicity result of concave and convex type.

\begin{theorem} \label{thm:intro-concave-convex}
	Assume $q < p^\ast$, $s > 1$, and
	\begin{align*}
		\max \set{N \left(1 - \frac{s}{p^\ast}\right),(\alpha + 1)\, s - \alpha p} < c < N \left(1 - \frac{s}{q^\ast_b}\right).
	\end{align*}
	Equivalently, $(s,c)$ lies in the open subscaled triangle with vertices $P_2$, $Q_1$, and $Q_2$. Then there exists $\lambda_\ast > 0$ such that, for all $\lambda \in (0,\lambda_\ast)$, equation \eqref{intro:critical} admits a sequence of weak radial solutions $\seq{u_k}$ such that $E(u_k) < 0$ and $E(u_k) \nearrow 0$.
\end{theorem}

This follows from Theorem \ref{Theorem concave convex} by taking $\mu_1 = \mu_2 = 1$ and $c_1 = c_2 = 0$. Indeed, \textnormal{(\hyperlink{H1}{H$_1$})} holds and $\rho < \sigma$ because $(s,c)$ belongs to $\Tc$ above the scaling line. Assumption \textnormal{(\hyperlink{H2}{H$_2$})} is satisfied as $P_1 \in [Q_1,P_1]$ and $P_3 \in [Q_2,P_3]$, while \textnormal{(\hyperlink{H4}{H$_4$})} follows from $q < p^\ast < q^\ast_b$.

\subsection{Organization of the paper}

The paper is organized as follows. In Section \ref{Section functional setting}, we study the functional framework, prove the sharp embedding theorem (Theorem \ref{Theorem embedding}, which contains Theorem \ref{thm:intro-embedding}) and a symmetric criticality result, examine the scaling and the associated scaled operator, and derive a Poho\v{z}aev identity for weak radial solutions. Section \ref{Section abstract} recalls the abstract variational theory of scaled functionals from Mercuri and Perera \cite{Mercuri-Perera-2026} and Perera \cite{Perera-2026}, which provides the main tools for our existence and multiplicity arguments. Section \ref{Section results} contains the general versions of Theorems \ref{thm:intro-eigenvalues}--\ref{thm:intro-concave-convex} and deals with three classes of problems. First, we study the scaled eigenvalue problem on the scaling line $Q_1Q_2$. The abstract theory leads to an unbounded sequence of eigenvalues $\lambda_k(s)$ constructed using the $\mathbb{Z}_2$ cohomological index of Fadell and Rabinowitz \cite{Fadell-Rabinowitz-1978} (Theorem \ref{Theorem 1}). Next, we treat pure power equations in the subscaled and superscaled regions. These equations admit infinitely many solutions with negative and positive energies, respectively (Theorem \ref{Theorem 2}). Finally, we investigate critical equations involving one or two critical terms on the segments $[Q_1,P_1]$ and $[Q_2,P_3]$. The main results include nonexistence statements based on the Poho\v{z}aev identity (Theorem \ref{Theorem nonexistence}), an explicit compactness threshold (Proposition \ref{Proposition PS}) based on the sharp algebraic estimate of Lemma \ref{Lemma threshold}, the optimality of this threshold for unweighted critical terms (Proposition \ref{Proposition sharp}), a multiplicity theorem for superscaled perturbations (Theorem \ref{Theorem large lambda}), and a multiplicity result of concave-convex type for subscaled perturbations (Theorem \ref{Theorem concave convex}). Some of these results rely on abstract critical point theorems involving the piercing property of the cohomological index. This is why we use eigenvalues defined through the cohomological index rather than the standard genus construction. We conclude the paper with a discussion of several open problems in Subsection \ref{Section open problems}.

\section{Functional setting} \label{Section functional setting}

\subsection{Density and spherical averages}

We begin with a useful criterion for membership in $W$.

\begin{lemma} \label{Lemma density}
	Let $u \in W^{1,1}_\loc(\R^N \setminus \set{0})$ satisfy
	\begin{equation} \label{eq:density}
		\begin{aligned}
			\int_{\R^N} \left(|\nabla u|^p + |x|^{-p}\, |u|^p\right)\,\mathrm{d}x &< \infty, \\
			\int_{\R^N} |x|^{-b} \left(|\nabla u|^q + |x|^{-q}\,  |u|^q\right)\, \mathrm{d}x &< \infty.
		\end{aligned}
	\end{equation}
	Then $u \in W$, and its gradient as an element of $W$ coincides with its weak gradient on $\R^N \setminus \set{0}$. Consequently, the quantities
	\begin{align*}
		\norm[\D^{1,p}(\R^N)]{u}^p\quad \text{and}\quad \norm[\D^{1,q}(\R^N;|x|^{-b})]{u}^q
	\end{align*}
	are given by the integrals of $|\nabla u|^p$ and $|x|^{-b}\, |\nabla u|^q$, respectively. Moreover, there are functions $\varphi_n \in C^\infty_c(\R^N \setminus \set{0})$ such that $\varphi_n \to u$ in $W$. If $u$ is radial, the functions $\varphi_n$ may also be chosen radial.
\end{lemma}

\begin{proof}
	Choose $\zeta \in C^\infty([0,\infty))$ such that
	\begin{align*}
		0 \le \zeta \le 1,\quad \zeta = 1\quad \text{on }[0,1],\quad \zeta = 0\quad \text{on }[2,\infty),\quad\text{and}\quad  |\zeta'| \le 2.
	\end{align*}
	For $0 < 2 \delta < m$, set
	\begin{align*}
		u_{\delta,m}(x) = \big(1 - \zeta(|x|/\delta)\big)\, \zeta(|x|/m)\, u(x).
	\end{align*}
	Then $u_{\delta,m}$ is supported in $\{x \in \R^N \colon\delta \le |x| \le 2m\}$ and $u - u_{\delta,m} = 0$ on $\{x \in \R^N \colon2\delta \le |x| \le m\}$. Since $2/\delta \le 4/|x|$ for $|x| \le 2\delta$ and $2/m \le 4/|x|$ for $|x| \le 2m$, we find
	\begin{align*}
		|\nabla (u - u_{\delta,m})| \le {\bf1}_{B_{2\delta} \cup (\R^N \setminus B_m)} \left(|\nabla u| + 4\, |x|^{-1}\, |u|\right) \quad \text{a.e.\,in } \R^N \setminus \set{0}.
	\end{align*}
	Thus, \eqref{eq:density} and the dominated convergence theorem imply
	\begin{align*}
		\int_{\R^N} |\nabla (u - u_{\delta,m})|^p\, \mathrm{d}x + \int_{\R^N} |x|^{-b}\, |\nabla (u - u_{\delta,m})|^q\, \mathrm{d}x \to 0 \quad \text{as } \delta \to 0,\, m \to \infty.
	\end{align*}

	Now fix $\delta$ and $m$, and let $w = u_{\delta,m}$. Since $|x|^{-b}$ is bounded from above and below by positive constants on $\{x \in \R^N \colon\delta/2 \le |x| \le 2m + 1\}$, \eqref{eq:density} gives $w \in W^{1,1}(\R^N)$ and $\nabla w \in L^p(\R^N) \cap L^q(\R^N)$. Let $\eta_\eps$ be a standard mollifier, which we take to be radial, and set $w_\eps = \eta_\eps * w$ for $0 < \eps < \min \set{\delta/2,1}$. Then $w_\eps \in C^\infty_c(\R^N \setminus \set{0})$ is supported in $\{x \in \R^N \colon\delta/2 \le |x| \le 2m + 1\}$, $w_\eps \to w$ in $L^1(\R^N)$, and $\nabla w_\eps = \eta_\eps * \nabla w \to \nabla w$ in $L^p(\R^N) \cap L^q(\R^N)$. Consequently,
	\begin{align*}
		\int_{\R^N} |\nabla (w_\eps - w)|^p\, \mathrm{d}x + \int_{\R^N} |x|^{-b}\,|\nabla (w_\eps - w)|^q\, \mathrm{d}x \to 0 \quad \text{as } \eps \to 0.
	\end{align*}
	A diagonal argument now provides functions $\varphi_n \in C^\infty_c(\R^N \setminus \set{0})$ such that $\varphi_n \to u$ in $L^1_\loc(\R^N \setminus \set{0})$ and
	\begin{align*}
		\int_{\R^N} |\nabla (\varphi_n - u)|^p\, \mathrm{d}x + \int_{\R^N} |x|^{-b}\,|\nabla (\varphi_n - u)|^q\, \mathrm{d}x \to 0.
	\end{align*}
	In particular, $\seq{\varphi_n}$ is a Cauchy sequence with respect to $\norm{\cdot}$ and hence converges in $W$ to some $\tilde u$. Since $W \hookrightarrow L^{p^\ast}(\R^N)$, $\varphi_n \to \tilde u$ in $L^{p^\ast}(\R^N)$, so $\tilde u = u$ a.e., and the gradient of $\tilde u$ is the $L^p$-limit of $\nabla \varphi_n$, namely, $\nabla u$. If $u$ is radial, then all the functions above are radial since $\zeta(|\cdot|/\delta)$, $\zeta(|\cdot|/m)$, and $\eta_\eps$ are radial.
\end{proof}

\begin{remark} 
	Conversely, every $u \in W$ satisfies \eqref{eq:density} by the Hardy inequality in $\D^{1,p}(\R^N)$ and the weighted Hardy inequality in $\D^{1,q}(\R^N;|x|^{-b})$. Both inequalities apply because $N > p$ and $N - b > q$. Hence Lemma \ref{Lemma density} characterizes $W$ as the space of functions $u \in W^{1,1}_\loc(\R^N \setminus \set{0})$ satisfying \eqref{eq:density}. In particular, $W = \D^{1,p}(\R^N) \cap \D^{1,q}(\R^N;|x|^{-b})$ when both spaces are realized as function spaces via these inequalities.
\end{remark}

For $v \in L^1_\loc(\R^N)$, let
\begin{align*}
	\bar v(x) = \frac{1}{\omega_{N-1}} \int_{S^{N-1}}\, v(|x| \xi)\,\mathrm{d}S(\xi),
\end{align*}
where $\omega_{N-1}$ denotes the surface area of $S^{N-1}$. Thus, $\bar{v}(x)$ is the average of $v$ over the sphere of radius $|x|$ centered at the origin.

We first record some basic properties of spherical averages.

\begin{lemma} \label{Lemma spherical}
	Let $v \in L^1_\loc(\R^N)$. The following statements hold:
	\begin{enumerate}
		\item[\textnormal{(i)}]\hypertarget{Lemma spherical.i}{}
			The function $\bar{v}$ is defined a.e.\,in $\R^N$, belongs to $L^1_\loc(\R^N)$, and is radial. Moreover, $\bar v = v$ a.e.\,whenever $v$ is radial, and for every radial measurable function $g$,
			\begin{align} \label{eq:average}
				\int_{\R^N} gv\, \mathrm{d}x = \int_{\R^N} g \bar v\, \mathrm{d}x
			\end{align}
			provided $gv \in L^1(\R^N)$, or $g \ge 0$ and $v \ge 0$.
		\item[\textnormal{(ii)}]\hypertarget{Lemma spherical.ii}{}
			For $1 \le m < \infty$, $|\bar v|^m \le \overline{|v|^m}$ a.e.\,in $\R^N$, and if $v \in L^m(\R^N)$, then $\norm[L^m(\R^N)]{\bar v} \le \norm[L^m(\R^N)]{v}$.
		\item[\textnormal{(iii)}]\hypertarget{Lemma spherical.iii}{}
			If $v$ is Lipschitz continuous with compact support, then $v \in W$ and its a.e.\,gradient coincides with its weak gradient. If $v$ is also radial, it can be approximated in $W$ by radial functions in $C^\infty_c(\R^N)$.
	\end{enumerate}
\end{lemma}

\begin{proof}
	\textnormal{(\hyperlink{Lemma spherical.i}{i})}
	Integration in polar coordinates gives
	\begin{align*}
		\int_0^R r^{N-1} \int_{S^{N-1}} |v(r \xi)|\,\mathrm{d}S(\xi)\,\mathrm{d}r = \int_{B_R} |v|\, \mathrm{d}x < \infty \quad \text{for all }R > 0.
	\end{align*}
	Thus, the inner integral is finite for a.e.\,$r > 0$, and $\bar v$ is defined a.e. Clearly, $\bar v$ is radial and $\bar v = v$ a.e.\,if $v$ is radial. Since $|\bar v| \le \overline{|v|}$, the same calculation shows that
	\begin{align*}
		\int_{B_R} |\bar v|\, \mathrm{d}x \le \int_{B_R} \overline{|v|}\, \mathrm{d}x = \int_{B_R} |v|\, \mathrm{d}x.
	\end{align*}
	We conclude that $\bar v \in L^1_\loc(\R^N)$.

	Write $g(x) = \gamma(|x|)$. If $g \ge 0$ and $v \ge 0$, Tonelli's theorem implies
	\begin{align*}
		\int_{\R^N} gv\, \mathrm{d}x
		& = \int_0^\infty r^{N-1}\, \gamma(r) \int_{S^{N-1}} v(r \xi)\,\mathrm{d}S(\xi)\,\mathrm{d}r\\
		&= \omega_{N-1} \int_0^\infty r^{N-1}\, \gamma(r) \bar v(r)\,\mathrm{d}r \\
		& = \int_{\R^N} g \bar v\, \mathrm{d}x.
	\end{align*}
	If $gv \in L^1(\R^N)$, the same calculation is justified by Fubini's theorem, since its absolute integrability follows by applying the preceding argument to $|g|$ and $|v|$.

	\textnormal{(\hyperlink{Lemma spherical.ii}{ii})}
	Since $\mathrm{d}S/\omega_{N-1}$ is a probability measure on $S^{N-1}$, Jensen's inequality gives $|\bar v|^m \le \overline{|v|^m}$ wherever the averages are defined. If $v \in L^m(\R^N)$, then \eqref{eq:average} with $g = 1$ applied to $|v|^m$ gives
	\begin{align*}
		\int_{\R^N} |\bar v|^m\, \mathrm{d}x \le \int_{\R^N} \overline{|v|^m}\, \mathrm{d}x = \int_{\R^N} |v|^m\, \mathrm{d}x.
	\end{align*}

	\textnormal{(\hyperlink{Lemma spherical.iii}{iii})}
	By Rademacher's theorem, $v$ is differentiable a.e., its a.e.\,gradient is its weak gradient, and $\nabla v$ is bounded. Since $v$ is also bounded and has compact support, it satisfies \eqref{eq:density}. Indeed, the weights $|x|^{-p}$ and $|x|^{-q-b}$ are locally integrable because $p < N$ and $q + b < N$. The assertions therefore follow from Lemma \ref{Lemma density}.
\end{proof}

We next show that spherical averaging preserves the relevant gradient bounds and allows radial approximation in $W$.

\begin{lemma} \label{Lemma averaging}
	The following assertions hold:
	\begin{enumerate}
		\item[\textnormal{(i)}]\hypertarget{Lemma averaging.i}{}
			If $v$ is Lipschitz continuous with compact support, then so is $\bar v$, and $\bar v \in W_\rad$. Furthermore,
			\begin{align*}
				\norm[\D^{1,p}(\R^N)]{\bar v} \le \norm[\D^{1,p}(\R^N)]{v}
			\end{align*}
			and
			\begin{align*}
				\norm[\D^{1,q}(\R^N;|x|^{-b})]{\bar v} \le \norm[\D^{1,q}(\R^N;|x|^{-b})]{v}.
			\end{align*}
			In particular, $\norm{\bar v} \le \norm{v}$.
		\item[\textnormal{(ii)}]\hypertarget{Lemma averaging.ii}{}
			If $u \in W_\rad$ and $v$ is Lipschitz continuous with compact support, then
			\begin{align*}
				& \int_{\R^N} \left(|\nabla u|^{p-2} + |x|^{-b}\,  |\nabla u|^{q-2}\right) \nabla u \cdot \nabla v\, \mathrm{d}x \\
				& = \int_{\R^N} \left(|\nabla u|^{p-2} + |x|^{-b}\,  |\nabla u|^{q-2}\right) \nabla u \cdot \nabla \bar v\, \mathrm{d}x.
			\end{align*}
		\item[\textnormal{(iii)}]\hypertarget{Lemma averaging.iii}{}
			Radial functions in $C^\infty_c(\R^N)$ are dense in $W_\rad$.
	\end{enumerate}
\end{lemma}

\begin{proof}
	\textnormal{(\hyperlink{Lemma averaging.i}{i})}
	Write $\bar v(x) = \psi(|x|)$, where
	\begin{align*}
		\psi(r) = \frac{1}{\omega_{N-1}} \int_{S^{N-1}} v(r \xi)\,\mathrm{d}S(\xi).
	\end{align*}
	Suppose that $v$ has Lipschitz constant $L$ and $\supp v \subset B_R$. Then
	\begin{align*}
		|\psi(r) - \psi(r')| \le L |r - r'|
	\end{align*}
	and $\psi(r) = 0$ for $r \ge R$. Thus, $\bar v$ is Lipschitz continuous with compact support, and Lemma \ref{Lemma spherical}  \textnormal{(\hyperlink{Lemma spherical.i}{i})}, \textnormal{(\hyperlink{Lemma spherical.iii}{iii})} shows that $\bar v \in W_\rad$. For almost every $\xi\in S^{N-1}$, the function $r \mapsto v(r \xi)$ is Lipschitz continuous with derivative $\nabla v(r \xi) \cdot \xi$ for a.e.\,$r$. By Fubini's theorem,
	\begin{equation} \label{eq:radial derivative}
		\psi'(r) = \frac{1}{\omega_{N-1}} \int_{S^{N-1}} \nabla v(r \xi) \cdot \xi\,\mathrm{d}S(\xi) \quad \text{for a.e.\,} r > 0.
	\end{equation}
	Thus, spherical averaging commutes with radial differentiation. It follows that $|\nabla \bar v(x)| = |\psi'(|x|)| \le \overline{|\nabla v|}(x)$, and Lemma \ref{Lemma spherical} \textnormal{(\hyperlink{Lemma spherical.ii}{ii})} gives $|\nabla \bar v|^m \le \overline{|\nabla v|^m}$ a.e.\,for $m \ge 1$. Taking $m = p$, or $m = q$ and multiplying by the radial weight $|x|^{-b}$, we may apply \eqref{eq:average} to conclude that
	\begin{align*}
		\norm[\D^{1,p}(\R^N)]{\bar v} \le \norm[\D^{1,p}(\R^N)]{v}\quad \text{and} \quad\norm[\D^{1,q}(\R^N;|x|^{-b})]{\bar v} \le \norm[\D^{1,q}(\R^N;|x|^{-b})]{v}.
	\end{align*}

	\textnormal{(\hyperlink{Lemma averaging.ii}{ii})}
	Write $u(x) = u(r)$, where $r = |x|$. By Lemma \ref{Lemma Strauss} \textnormal{(\hyperlink{Lemma Strauss.iii}{iii})} below,
	\begin{align*}
		\nabla u(x) = u'(r)\, x/r\quad \text{and}\quad \nabla u \cdot \nabla v = u'(r)\, \nabla v(x) \cdot x/r\quad \text{a.e.}
	\end{align*}
	Set $A(r) = |u'(r)|^{p-2} + r^{-b} |u'(r)|^{q-2}$. Using polar coordinates and \eqref{eq:radial derivative}, we arrive at
	\begin{align*}
		\int_{\R^N} A(|x|)\, \nabla u \cdot \nabla v\, \mathrm{d}x & = \omega_{N-1} \int_0^\infty r^{N-1}\, A(r)\, u'(r) \psi'(r)\,\mathrm{d}r \\
		& = \int_{\R^N} A(|x|)\, \nabla u \cdot \nabla \bar v\, \mathrm{d}x.
	\end{align*}

	\textnormal{(\hyperlink{Lemma averaging.iii}{iii})}
	Let $u \in W_\rad$ and take $\varphi_n \in C^\infty_c(\R^N)$ with $\varphi_n \to u$ in $W$. By part \textnormal{(\hyperlink{Lemma averaging.i}{i})}, $\bar \varphi_n \in W_\rad$, and
	\begin{align*}
		\norm{\bar \varphi_n - \bar \varphi_m} = \norm{\overline{\varphi_n - \varphi_m}} \le \norm{\varphi_n - \varphi_m}.
	\end{align*}
	Hence $\seq{\bar \varphi_n}$ is a Cauchy sequence in $W_\rad$ and converges to some $\tilde u \in W_\rad$. Since $u$ is radial, Lemma \ref{Lemma spherical} \textnormal{(\hyperlink{Lemma spherical.i}{i})}, \textnormal{(\hyperlink{Lemma spherical.ii}{ii})} gives
	\begin{align*}
		\norm[L^{p^\ast}(\R^N)]{\bar \varphi_n - u} = \norm[L^{p^\ast}(\R^N)]{\bar \varphi_n - \bar u} = \norm[L^{p^\ast}(\R^N)]{\overline{\varphi_n - u}} \le \norm[L^{p^\ast}(\R^N)]{\varphi_n - u} \to 0.
	\end{align*}
	Thus, $\tilde u = u$. Finally, each $\bar \varphi_n$ is radial, Lipschitz continuous, and compactly supported. By Lemma \ref{Lemma spherical} \textnormal{(\hyperlink{Lemma spherical.iii}{iii})}, it can be approximated in $W$ by radial functions in $C^\infty_c(\R^N)$. The required density follows.
\end{proof}

\begin{remark}
	The proof of Lemma \ref{Lemma averaging} \textnormal{(\hyperlink{Lemma averaging.ii}{ii})} uses Lemma \ref{Lemma Strauss} \textnormal{(\hyperlink{Lemma Strauss.iii}{iii})} below, whose proof only uses part \textnormal{(\hyperlink{Lemma averaging.iii}{iii})} of Lemma \ref{Lemma averaging}, and the proof of part \textnormal{(\hyperlink{Lemma averaging.iii}{iii})} only uses part \textnormal{(\hyperlink{Lemma averaging.i}{i})}. So there is no circularity.
\end{remark}

\subsection{Radial estimates and the sharp embedding theorem}

We begin with two pointwise estimates for radial functions, which will be essential for the embedding results below.

\begin{lemma} \label{Lemma Strauss}
	Every $u \in W_\rad$ has a radial representative $u(x) = u(r)$, where $r = |x|$, such that $u$ is continuous on $(0,\infty)$ and the following hold for some constant $C > 0$, independent of $u$:
	\begin{enumerate}
		\item[\textnormal{(i)}]\hypertarget{Lemma Strauss.i}{}
			$|u(r)| \le Cr^{-N/p^\ast} \norm[\D^{1,p}(\R^N)]{u}$ for all $r > 0$;
		\item[\textnormal{(ii)}]\hypertarget{Lemma Strauss.ii}{}
			$|u(r)| \le Cr^{-N/q^\ast_b} \norm[\D^{1,q}(\R^N;|x|^{-b})]{u}$ for all $r > 0$;
		\item[\textnormal{(iii)}]\hypertarget{Lemma Strauss.iii}{}
			The function $u$ is locally absolutely continuous on $(0,\infty)$,
			\begin{align*}
				\nabla u(x) = u'(|x|)\,x/|x| \quad\text{for a.e.\,}x \in \R^N,
			\end{align*}
			and
			\begin{align*}
				\norm[\D^{1,p}(\R^N)]{u}^p& = \omega_{N-1} \int_0^\infty r^{N-1}\, |u'(r)|^p\,\mathrm{d}r, \\
				\norm[\D^{1,q}(\R^N;|x|^{-b})]{u}^q & = \omega_{N-1} \int_0^\infty r^{N-b-1}\, |u'(r)|^q\,\mathrm{d}r.
			\end{align*}
	\end{enumerate}
\end{lemma}

\begin{proof}
	We first assume that $u\in C^\infty_c(\R^N)$ is radial. For $r > 0$, H\"{o}lder's inequality implies
	\begin{align*}
		|u(r)|
		& = \abs{\int_r^\infty u'(\tau)\, \tau^{(N-1)/p}\, \tau^{-(N-1)/p}\,\mathrm{d}\tau} \\
		& \le \left(\int_r^\infty |u'(\tau)|^p\, \tau^{N-1}\,\mathrm{d}\tau\right)^{1/p} \left(\int_r^\infty \tau^{-(N-1)/(p-1)}\,\mathrm{d}\tau\right)^{(p-1)/p} \\
		& \le Cr^{-N/p^\ast} \norm[\D^{1,p}(\R^N)]{u},
	\end{align*}
	as $N > p$. Similarly, writing
	\begin{align*}
		u'(\tau) = u'(\tau)\, \tau^{(N-b-1)/q}\, \tau^{-(N-b-1)/q}
	\end{align*}
	and using $N - b > q$ we derive the estimate in \textnormal{(\hyperlink{Lemma Strauss.ii}{ii})}.

	For general $u \in W_\rad$, Lemma \ref{Lemma averaging} \textnormal{(\hyperlink{Lemma averaging.iii}{iii})} provides radial functions $u_n \in C^\infty_c(\R^N)$ such that $u_n \to u$ in $W$. Applying \textnormal{(\hyperlink{Lemma Strauss.i}{i})} to $u_n - u_m$, we see that $\seq{u_n}$ is uniformly Cauchy on $\{x \in \R^N \colon|x| \ge \delta\}$ for every $\delta > 0$. It therefore converges locally uniformly on $\R^N \setminus \set{0}$ to a continuous radial function. Since $u_n \to u$ in $L^{p^\ast}(\R^N)$, this limit coincides with $u$ a.e. Passing to the limit gives \textnormal{(\hyperlink{Lemma Strauss.i}{i})} and \textnormal{(\hyperlink{Lemma Strauss.ii}{ii})}.

	 For $0 < \delta < R$, we also have
	\begin{align*}
		\int_\delta^R |u_n'(r) - u_m'(r)|^p\, r^{N-1}\,\mathrm{d}r \le \frac{\norm[\D^{1,p}(\R^N)]{u_n - u_m}^p}{\omega_{N-1}} \to 0.
	\end{align*}
	Consequently, $u_n' \to w$ in $L^p(\delta,R)$ for some $w \in L^p_\loc(0,\infty)$. Passing to the limit in $u_n(r') - u_n(r) = \int_r^{r'} u_n'(\tau)\,\mathrm{d}\tau$ shows that $u(r') - u(r) = \int_r^{r'} w(\tau)\,\mathrm{d}\tau$. Thus, $u$ is locally absolutely continuous on $(0,\infty)$ with $u' = w$ a.e. Since $\nabla u_n(x) = u_n'(|x|)\, x/|x|$ converges to $\nabla u$ in $L^p(\R^N)$ and to $w(|x|)\, x/|x|$ in $L^p_\loc(\R^N \setminus \set{0})$, $\nabla u(x) = u'(|x|)\, x/|x|$ a.e., and the formulas for $\norm[\D^{1,p}(\R^N)]{u}$ and  $\norm[\D^{1,q}(\R^N;|x|^{-b})]{u}$ follow by integration in polar coordinates.
\end{proof}

The pointwise estimates lead to the following embedding result for points in the interior of the triangle.

\begin{lemma} \label{Lemma wedge}
	Assume that $s > 0$ and
	\begin{align}\label{lemma wedge assumption}
		N\,\left (1 - \frac{s}{p^\ast}\right) < c < N\,\left (1 - \frac{s}{q^\ast_b}\right).
	\end{align}
	Then there exists $C>0$ such that
	\begin{align*}
		\int_{\R^N} |x|^{-c}\, |u|^s\, \mathrm{d}x \le C\, \big(\norm[\D^{1,p}(\R^N)]{u}^s + \norm[\D^{1,q}(\R^N;|x|^{-b})]{u}^s\big)
	\end{align*}
	for all $u \in W_\rad$.
\end{lemma}

\begin{proof}
	By Lemma \ref{Lemma Strauss} \textnormal{(\hyperlink{Lemma Strauss.i}{i})} and \textnormal{(\hyperlink{Lemma Strauss.ii}{ii})},
	\begin{align*}
		\int_{B_1^c} |x|^{-c}\, |u|^s\, \mathrm{d}x &\le C \norm[\D^{1,p}(\R^N)]{u}^s \int_1^\infty r^{N\, (1 - s/p^\ast) - c - 1}\,\mathrm{d}r,\\
		\int_{B_1} |x|^{-c}\, |u|^s\, \mathrm{d}x &\le C \norm[\D^{1,q}(\R^N;|x|^{-b})]{u}^s \int_0^1 r^{N\, (1 - s/q^\ast_b) - c - 1}\,\mathrm{d}r
	\end{align*}
	and both integrals on the right-hand sides are finite because \eqref{lemma wedge assumption}.
\end{proof}

We now consider the two critical boundary segments of the triangle.

\begin{lemma} \label{Lemma critical segments}
	The following statements hold.
	\begin{enumerate}
		\item[\textnormal{(i)}]\hypertarget{Lemma critical segments.i}{}
			If $(s,c) \in [Q_1,P_1]$, then $\norm[L^s(\R^N;|x|^{-c})]{u} \le C \norm[\D^{1,p}(\R^N)]{u}$ for all $u \in W$.
		\item[\textnormal{(ii)}]\hypertarget{Lemma critical segments.ii}{}
			If $(s,c) \in [Q_2,P_3]$, then $\norm[L^s(\R^N;|x|^{-c})]{u} \le C \norm[\D^{1,q}(\R^N;|x|^{-b})]{u}$ for all $u \in W_\rad$. In particular, $W_\rad \hookrightarrow L^{q^\ast_b}(\R^N)$.
	\end{enumerate}
\end{lemma}

\begin{proof}
	\textnormal{(\hyperlink{Lemma critical segments.i}{i})} Write $(s,c) = \theta P_1 + (1 - \theta)\, Q_1$ with $\theta \in [0,1]$. Then
	\begin{align*}
		s = \theta p^\ast + (1 - \theta)\, p\quad \text{and}\quad c = (1 - \theta)\, p.
	\end{align*}
	H\"{o}lder's inequality gives
	\begin{align*}
		\int_{\R^N} |x|^{-c}\, |u|^s\, \mathrm{d}x
		& = \int_{\R^N} |u|^{\theta p^\ast}\, (|x|^{-p} |u|^p)^{1 - \theta}\, \mathrm{d}x \\
		& \le \norm[L^{p^\ast}(\R^N)]{u}^{\theta p^\ast} \left(\int_{\R^N} |x|^{-p}\, |u|^p\, \mathrm{d}x\right)^{1 - \theta}.
	\end{align*}
	The Sobolev and Hardy inequalities now imply \textnormal{(\hyperlink{Lemma critical segments.i}{i})}.

	\textnormal{(\hyperlink{Lemma critical segments.ii}{ii})}
	For $(s,c)\in [Q_2,P_3]$, we have $q \leq s \leq q^\ast_b$ and $c = N\, (1 - s/q^\ast_b)$. Therefore,
	\begin{align*}
		c + (s - q)\, \frac{N}{q^\ast_b} = N - \frac{Nq}{q^\ast_b} = q + b.
	\end{align*}
	By Lemma \ref{Lemma Strauss} \textnormal{(\hyperlink{Lemma Strauss.ii}{ii})},
	\begin{align*}
		|x|^{-c}\, |u|^s = |x|^{-c}\, |u|^{s-q}\, |u|^q \le C \norm[\D^{1,q}(\R^N;|x|^{-b})]{u}^{s-q}\, |x|^{-q-b}\, |u|^q,
	\end{align*}
	and the conclusion follows from the weighted Hardy inequality.
\end{proof}

Combining the preceding estimates, we obtain a sharp characterization of the continuous and compact embeddings of $W_\rad$.

\begin{theorem} \label{Theorem embedding}
	Let $s \ge 1$ and $c \ge 0$. Then the following hold:
	\begin{enumerate}
		\item[\textnormal{(i)}]\hypertarget{Theorem embedding.i}{}
			If $(s,c) \in \Te$, then the embedding $W_\rad \hookrightarrow L^s(\R^N;|x|^{-c})$ is continuous.
		\item[\textnormal{(ii)}]\hypertarget{Theorem embedding.ii}{}
			If $(s,c) \in \Tc$, then the embedding $W_\rad \hookrightarrow L^s(\R^N;|x|^{-c})$ is compact.
		\item[\textnormal{(iii)}]\hypertarget{Theorem embedding.iii}{}
			If $(s,c) \in [Q_1,P_1] \cup [Q_2,P_3]$, then the embedding $W_\rad \hookrightarrow L^s(\R^N;|x|^{-c})$ is not compact.
		\item[\textnormal{(iv)}]\hypertarget{Theorem embedding.iv}{}
			If $(s,c) \notin \Te$, then $W_\rad$ is not contained in $L^s(\R^N;|x|^{-c})$.
	\end{enumerate}
\end{theorem}

\begin{proof}
	\textnormal{(\hyperlink{Theorem embedding.i}{i})} The assertion follows from Lemma \ref{Lemma wedge} when $(s,c) \in \Tc$ and from Lemma \ref{Lemma critical segments} when $(s,c) \in [Q_1,P_1] \cup [Q_2,P_3]$.

	\textnormal{(\hyperlink{Theorem embedding.ii}{ii})} Let $\seq{u_j}$ be a bounded sequence in $W_\rad$. Passing to a subsequence, we may assume that $u_j \wto u$ in $W_\rad$. Set $v_j = u_j - u$. Then $v_j \wto 0$ and $\norm{v_j} \le M$ for some $M > 0$. For $R > \rho > 0$,
	\begin{equation} \label{12}
		\begin{aligned}
			&\int_{\R^N} |x|^{-c}\, |v_j|^s\, \mathrm{d}x\\
			&= \int_{B_\rho} |x|^{-c}\, |v_j|^s\, \mathrm{d}x + \int_{A_{\rho,R}} |x|^{-c}\, |v_j|^s\, \mathrm{d}x + \int_{B_R^c} |x|^{-c}\, |v_j|^s\, \mathrm{d}x,
		\end{aligned}
	\end{equation}
	where $A_{\rho,R} = B_R \setminus B_\rho$. By Lemma \ref{Lemma Strauss} \textnormal{(\hyperlink{Lemma Strauss.ii}{ii})},
	\begin{align*}
		\int_{B_\rho} |x|^{-c}\, |v_j|^s\, \mathrm{d}x = \omega_{N-1} \int_0^\rho r^{N-c-1}\, |v_j(r)|^s\,\mathrm{d}r \le CM^s \int_0^\rho r^{N\,(1 - s/q^\ast_b) - c - 1}\,\mathrm{d}r.
	\end{align*}
	The last integral tends to zero as $\rho \searrow 0$, because $c < N\, (1 - s/q^\ast_b)$. Similarly, Lemma \ref{Lemma Strauss} \textnormal{(\hyperlink{Lemma Strauss.i}{i})} implies
	\begin{align*}
		\int_{B_R^c} |x|^{-c}\, |v_j|^s\, \mathrm{d}x = \omega_{N-1} \int_R^\infty r^{N-c-1}\, |v_j(r)|^s\,\mathrm{d}r \le C M^s \int_R^\infty r^{N\,(1 - s/p^\ast) - c - 1}\,\mathrm{d}r,
	\end{align*}
	and the last integral goes to zero as $R \nearrow \infty$ since $c > N\, (1 - s/p^\ast)$. So the first and the last integrals in \eqref{12} go to zero as $\rho \searrow 0$ and $R \nearrow \infty$, uniformly in $j$.

	For fixed $R > \rho > 0$, we claim that $v_j \to 0$ uniformly on $A_{\rho,R}$. By Lemma \ref{Lemma Strauss} \textnormal{(\hyperlink{Lemma Strauss.i}{i})}, $|v_j(r)| \le CM \rho^{-N/p^\ast}$ for $r \in [\rho,R]$, and for $\rho \le r < r' \le R$, Lemma \ref{Lemma Strauss} \textnormal{(\hyperlink{Lemma Strauss.iii}{iii})} and H\"{o}lder's inequality give
	\begin{align*}
		|v_j(r') - v_j(r)|
		& \le \int_r^{r'} |v_j'(\tau)|\,\mathrm{d}\tau \\
		& \le \left(\int_\rho^R \tau^{N-1}\, |v_j'(\tau)|^p\,\mathrm{d}\tau\right)^{1/p} \rho^{-(N-1)/p}\, (r' - r)^{(p-1)/p} \\
		& \le C_\rho\,  M\, (r' - r)^{(p-1)/p}.
	\end{align*}
	So $\seq{v_j}$ is bounded and equicontinuous on $[\rho,R]$. By the Arzel\`{a}-Ascoli theorem, every subsequence has a further subsequence that converges uniformly on $[\rho,R]$. Since $v_j \wto 0$ in $W$ and $W \hookrightarrow L^{p^\ast}\!(\R^N)$, we also have $v_j \wto 0$ in $L^{p^\ast}\!(A_{\rho,R})$, so the uniform limit is $0$, and hence the whole sequence $v_j \to 0$ uniformly on $A_{\rho,R}$. Consequently,
	\begin{align*}
		\int_{A_{\rho,R}} |x|^{-c}\, |v_j|^s\, \mathrm{d}x \le \rho^{-c} \vol{A_{\rho,R}}\, \sup_{A_{\rho,R}} |v_j|^s \to 0.
	\end{align*}
	Returning to \eqref{12}, we obtain $v_j \to 0$ in $L^s(\R^N;|x|^{-c})$, as required.

	\textnormal{(\hyperlink{Theorem embedding.iii}{iii})} Suppose first that $(s,c) \in [Q_1,P_1]$, so $c = N\, (1 - s/p^\ast)$. Let $u \in C^\infty_c(\R^N) \setminus \set{0}$ be a radial function that vanishes in a neighborhood of the origin, let $\eps_j \searrow 0$, and let
	\begin{align*}
		u_j(x) = \eps_j^{N/p^\ast} u(\eps_j\, x), \quad x \in \R^N.
	\end{align*}
	Then
	\begin{align*}
		\norm{u_j} = \norm[\D^{1,p}(\R^N)]{u} + \eps_j^{N/p - (N-b)/q} \norm[\D^{1,q}(\R^N;|x|^{-b})]{u}
	\end{align*}
	and $N/p - (N - b)/q > 0$, so $\seq{u_j}$ is bounded in $W_\rad$. Passing to a subsequence, we may assume that $u_j \wto \tilde u$ for some $\tilde u \in W_\rad$. Since $W \hookrightarrow L^{p^\ast}\!(\R^N)$, $u_j \wto \tilde u$ in $L^{p^\ast}\!(\R^N)$. Since $\dist{0}{\supp u_j} \to \infty$ and hence $u_j \to 0$ pointwise, it follows that $\tilde u = 0$. However, $\norm[L^s(\R^N;|x|^{-c})]{u_j} = \norm[L^s(\R^N;|x|^{-c})]{u}$ and hence $u_j \nrightarrow 0$ in $L^s(\R^N;|x|^{-c})$. The case $(s,c) \in [Q_2,P_3]$ is similar, using $u_j(x) = \eps_j^{-N/q^\ast_b}\, u(x/\eps_j)$, for which
	\begin{align*}
		\norm{u_j} = \eps_j^{N/p - (N-b)/q} \norm[\D^{1,p}(\R^N)]{u} + \norm[\D^{1,q}(\R^N;|x|^{-b})]{u}
	\end{align*}
	and $\supp u_j$ shrinks to the origin.

	\textnormal{(\hyperlink{Theorem embedding.iv}{iv})} Let $(s,c) \notin \Te$ with $s \ge 1$ and $c \ge 0$. Then one of the following holds:
	\begin{enumerate}
		\item[\textnormal{(1)}]\hypertarget{Theorem embedding.iv.1}{} $c < N\, (1 - s/p^\ast)$;
		\item[\textnormal{(2)}]\hypertarget{Theorem embedding.iv.2}{} $c > N\, (1 - s/q^\ast_b)$;
		\item[\textnormal{(3)}]\hypertarget{Theorem embedding.iv.3}{} $(s,c) \in [P_2,Q_1)$, i.e., $c = N\, (1 - s/p^\ast)$ and $s < p$;
		\item[\textnormal{(4)}]\hypertarget{Theorem embedding.iv.4}{} $(s,c) \in [P_2,Q_2)$, i.e., $c = N\,(1 - s/q^\ast_b)$ and $s < q$.
	\end{enumerate}
	For cases \textnormal{(\hyperlink{Theorem embedding.iv.1}{1})} and \textnormal{(\hyperlink{Theorem embedding.iv.2}{2})}, let $u \in C^\infty_c(\R^N) \setminus \set{0}$ be a radial function that vanishes in a neighborhood of the origin and set $u_\eps(x) = u(x/\eps)$. A change of variables gives
	\begin{align*}
		\norm[\D^{1,p}(\R^N)]{u_\eps} &= \eps^{(N-p)/p} \norm[\D^{1,p}(\R^N)]{u},\\
		\norm[\D^{1,q}(\R^N;|x|^{-b})]{u_\eps} &= \eps^{(N-q-b)/q} \norm[\D^{1,q}(\R^N;|x|^{-b})]{u},\\
		\norm[L^s(\R^N;|x|^{-c})]{u_\eps} &= \eps^{(N-c)/s} \norm[L^s(\R^N;|x|^{-c})]{u}.
	\end{align*}
	Note that
	\begin{align*}
		\frac{N - q - b}{q} < \frac{N - p}{p}.
	\end{align*}
	In case \textnormal{(\hyperlink{Theorem embedding.iv.1}{1})}, $(N - c)/s > (N - p)/p$ and hence
	\begin{align*}
		\frac{\norm[L^s(\R^N;|x|^{-c})]{u_\eps}}{\norm{u_\eps}} \to \infty \quad \text{as } \eps \to \infty.
	\end{align*}
	In case \textnormal{(\hyperlink{Theorem embedding.iv.2}{2})}, $(N - c)/s < (N - q - b)/q$ and hence this ratio goes to infinity as $\eps \to 0$. So there is no continuous embedding, and since the embedding would be continuous by the closed graph theorem if $W_\rad \subset L^s(\R^N;|x|^{-c})$, $W_\rad$ is not contained in $L^s(\R^N;|x|^{-c})$.

	In case \textnormal{(\hyperlink{Theorem embedding.iv.3}{3})}, fix $\gamma \in (1/p,1/s)$ and let $u(x) = h(|x|)$, where $h \in C^\infty([0,\infty))$ is constant near zero and $h(r) = r^{-(N-p)/p}\, (\log r)^{-\gamma}$ for $r \ge 3$. Then $|h'(r)| \le C r^{-N/p}\, (\log r)^{-\gamma}$ for $r \ge 3$, so
	\begin{align*}
		\int_3^\infty r^{N-1}\, |h'(r)|^p\,\mathrm{d}r \le C \int_3^\infty \frac{\mathrm{d}r}{r\,(\log r)^{\gamma p}} < \infty
	\end{align*}
	since $\gamma p > 1$. Moreover,
	\begin{align*}
		\int_3^\infty r^{N-b-1}\, |h'(r)|^q\,\mathrm{d}r \le C \int_3^\infty r^{N-b-Nq/p-1}\,\mathrm{d}r < \infty
	\end{align*}
	as $Nq/p > N - b$, while
	\begin{align*}
		\int_3^\infty r^{N-c-1}\, |h(r)|^s\,\mathrm{d}r = \int_3^\infty \frac{\mathrm{d}r}{r\,(\log r)^{\gamma s}} = \infty
	\end{align*}
	since $N - c = (N - p)\, s/p$ and $\gamma s < 1$. Furthermore,
	\begin{align*}
		\int_3^\infty r^{N-p-1}\, |h(r)|^p\,\mathrm{d}r = \int_3^\infty \frac{\mathrm{d}r}{r\, (\log r)^{\gamma p}} < \infty
	\end{align*}
	and
	\begin{align*}
		\int_3^\infty r^{N-q-b-1}\, |h(r)|^q\,\mathrm{d}r \le \int_3^\infty r^{N-b-Nq/p-1}\,\mathrm{d}r < \infty.
	\end{align*}
	Since $h$ is smooth on $[0,\infty)$, it follows that $u$ satisfies \eqref{eq:density}, so $u \in W_\rad$ by Lemma \ref{Lemma density}. Hence $u \in W_\rad \setminus L^s(\R^N;|x|^{-c})$.

	In case \textnormal{(\hyperlink{Theorem embedding.iv.4}{4})}, choose $\gamma \in (1/q,1/s)$ and let $u(x) = h(|x|)$, where $h \in C^\infty(0,\infty)$ vanishes on $[1,\infty)$ and $h(r) = r^{-(N-q-b)/q} \abs{\log r}^{-\gamma}$ for $0 < r \le 1/3$. Then $|h'(r)| \le Cr^{-(N-b)/q} \abs{\log r}^{-\gamma}$ for $0 < r \le 1/3$, so
	\begin{align*}
		\int_0^{1/3} r^{N-1}\, |h'(r)|^p\,\mathrm{d}r \le C \int_0^{1/3} r^{N-p\,(N-b)/q-1}\,\mathrm{d}r < \infty
	\end{align*}
	since $p(N - b)/q < N$, and
	\begin{align*}
		\int_0^{1/3} r^{N-b-1}\, |h'(r)|^q\,\mathrm{d}r \le C \int_0^{1/3} \frac{\mathrm{d}r}{r\, \abs{\log r}^{\gamma q}} < \infty
	\end{align*}
	since $\gamma q > 1$. However,
	\begin{align*}
		\int_0^{1/3} r^{N-c-1}\, |h(r)|^s\,\mathrm{d}r = \int_0^{1/3} \frac{\mathrm{d}r}{r\, \abs{\log r}^{\gamma s}} = \infty
	\end{align*}
	since $N - c = s\, (N - q - b)/q$ and $\gamma s < 1$. Finally,
	\begin{align*}
		\int_0^{1/3} r^{N-p-1}\, |h(r)|^p\,\mathrm{d}r \le C \int_0^{1/3} r^{N-p\,(N-b)/q-1}\,\mathrm{d}r < \infty
	\end{align*}
	and
	\begin{align*}
		\int_0^{1/3} r^{N-q-b-1}\, |h(r)|^q\,\mathrm{d}r = \int_0^{1/3} \frac{\mathrm{d}r}{r\, \abs{\log r}^{\gamma q}} < \infty.
	\end{align*}
	As $h$ is smooth on $(0,\infty)$ and vanishes on $[1,\infty)$, it follows that $u$ satisfies \eqref{eq:density}, so $u \in W_\rad$ by Lemma \ref{Lemma density}. Hence $u \in W_\rad \setminus L^s(\R^N;|x|^{-c})$.
\end{proof}

\begin{remark} \label{Remark 1}
	When $b > 0$, Theorem \ref{Theorem embedding} does not extend to the whole space $W$ near the vertex $P_3$. To see this, let $\varphi \in C^\infty_c(\R^N) \setminus \set{0}$, choose $\eps \in (0,1)$, and set
	\begin{align*}
		u_\eps(x) = \eps^{-(N-p)/p}\, \varphi((x - x_\eps)/\eps), \quad |x_\eps| = \eps^{-N\,(q-p)/pb}.
	\end{align*}
	Then $\norm[\D^{1,p}(\R^N)]{u_\eps}$ is independent of $\eps$, while $\norm[\D^{1,q}(\R^N;|x|^{-b})]{u_\eps}$ remains bounded as $\eps \to 0$. On the other hand,
	\begin{align*}
		\norm[L^{q^\ast_b}(\R^N)]{u_\eps}^{q^\ast_b} = \eps^{N\,(1-q^\ast_b/p^\ast)} \norm[L^{q^\ast_b}(\R^N)]{\varphi}^{q^\ast_b} \to \infty.
	\end{align*}
	The Caffarelli-Kohn-Nirenberg inequality \cite{Caffarelli-Kohn-Nirenberg-1984}, however, gives
	\begin{align*}
		\D^{1,q}(\R^N;|x|^{-b}) \hookrightarrow L^s(\R^N;|x|^{-c})
	\end{align*}
	for $(s,c)$ on the segment of $[Q_2,P_3]$ joining $Q_2$ to
	\begin{align*}
		P_3' = \left(q^\ast,\frac{Nb}{N - q}\right), \quad q^\ast = \frac{Nq}{N - q}.
	\end{align*}
	H\"{o}lder's inequality therefore implies $W \hookrightarrow L^s(\R^N;|x|^{-c})$ for all $(s,c)$ in the convex hull of $P_1$, $Q_1$, $Q_2$, and $P_3'$. By restricting ourselves to $W_\rad$, we can work in the considerably larger region $\Te$.
\end{remark}

\subsection{Weak radial solutions and symmetric criticality}

Although the energy functional $E$ need not be defined on the whole $W$, radial critical points satisfy the weak formulation of equation \eqref{1} for all test functions in $W$.

\begin{proposition} \label{Proposition symmetric criticality}
	Assume the growth condition \eqref{6}. If $u \in W_\rad$ is a weak radial solution of equation \eqref{1}, then $f(|x|,u) v \in L^1(\R^N)$ for all $v \in W$, and
	\begin{align*}
		\int_{\R^N} \left(|\nabla u|^{p-2} + |x|^{-b}\, |\nabla u|^{q-2}\right) \nabla u \cdot \nabla v\, \mathrm{d}x = \int_{\R^N} f(|x|,u)\, v\, \mathrm{d}x
	\end{align*}
	for all $v \in W$.
\end{proposition}

\begin{proof}
	Set $g = \sum_{i=1}^n a_i\, |x|^{-c_i}\, |u|^{s_i-1}$. Then $g$ is radial and nonnegative, and \eqref{6} implies $|f(|x|,u)| \le g$. We first show that
	\begin{equation} \label{eq:g bound}
		\int_{\R^N} g\, |v|\, \mathrm{d}x \le C \norm{v} \quad \text{for all }v \in W,
	\end{equation}
	where $C>0$ may depend on $u$. Let $v \in C^\infty_c(\R^N)$. Then the function $|v|$ is Lipschitz continuous with compact support, so $|v| \in W$ and $\norm{|v|} = \norm{v}$ by Lemma \ref{Lemma spherical} \textnormal{(\hyperlink{Lemma spherical.iii}{iii})} since $|\nabla |v|| = |\nabla v|$ a.e. Using \eqref{eq:average}, H\"{o}lder's inequality, Theorem \ref{Theorem embedding} \textnormal{(\hyperlink{Theorem embedding.i}{i})}, and Lemma \ref{Lemma averaging} \textnormal{(\hyperlink{Lemma averaging.i}{i})}, we find
	\begin{align*}
		\int_{\R^N} g\, |v|\, \mathrm{d}x
		& = \int_{\R^N} g\, \overline{|v|}\, \mathrm{d}x \\
		& \le \sum_{i=1}^n a_i \norm[L^{s_i}(\R^N;|x|^{-c_i})]{u}^{s_i-1} \big\|\overline{|v|}\big\|_{L^{s_i}(\R^N;|x|^{-c_i})} \\
		& \le C\, \big\|\overline{|v|}\big\| \le C \norm{v}.
	\end{align*}
	For general $v \in W$, choose $\varphi_n \in C^\infty_c(\R^N)$ with $\varphi_n \to v$ in $W$. Since $W \hookrightarrow L^{p^\ast}(\R^N)$, a subsequence also converges to $v$ a.e.\,in $\R^N$. Fatou's lemma extends \eqref{eq:g bound} to all $v\in W$. It follows that $f(|x|,u)\,v \in L^1(\R^N)$ for all $v \in W$, and the functional
	\begin{align*}
		L(v) = \int_{\R^N} \left(|\nabla u|^{p-2} + |x|^{-b}\, |\nabla u|^{q-2}\right) \nabla u \cdot \nabla v\, \mathrm{d}x - \int_{\R^N} f(|x|,u)\, v\, \mathrm{d}x
	\end{align*}
	belongs to $W^\ast$. For $v \in C^\infty_c(\R^N)$, Lemma \ref{Lemma averaging} \textnormal{(\hyperlink{Lemma averaging.ii}{ii})} and \eqref{eq:average}, applied to the radial function $f(|x|,u)$, show that
	\begin{align*}
		L(v) = L(\bar v) = 0,
	\end{align*}
	because $u$ is a weak radial solution and $\bar v \in W_\rad$ by Lemma \ref{Lemma averaging} \textnormal{(\hyperlink{Lemma averaging.i}{i})}. Since $C^\infty_c(\R^N)$ is dense in $W$, we conclude that $L = 0$.
\end{proof}

\subsection{The scaling and the scaled operator} 

Set
\begin{equation} \label{scaling}
	u_t(x) =
	\begin{cases}
		t^{\alpha + 1}\,u(x/t) & \text{if } t > 0\\
		0 & \text{if } t = 0,
	\end{cases}
\end{equation}
where $\alpha \ge 0$ is given by \eqref{100}. A change of variables leads to
\begin{equation} \label{scaling relations}
	\begin{aligned}
		\int_{\R^N} |\nabla u_t|^p\, \mathrm{d}x & = t^\sigma \int_{\R^N} |\nabla u|^p\, \mathrm{d}x, \\
		\int_{\R^N} |x|^{-b}\, |\nabla u_t|^q\, \mathrm{d}x & = t^\sigma \int_{\R^N} |x|^{-b}\, |\nabla u|^q\, \mathrm{d}x, \\
		\int_{\R^N} |x|^{-c}\, |u_t|^s\, \mathrm{d}x & = t^{\rho(s,c)} \int_{\R^N} |x|^{-c}\, |u|^s\, \mathrm{d}x,
	\end{aligned}
\end{equation}
where $\sigma$ and $\rho(s,c)$ are defined in \eqref{rho sigma}. In particular,
\begin{equation} \label{eq.1.prop1}
	\norm{u_t} = t^{\sigma/p} \norm[\D^{1,p}(\R^N)]{u} + t^{\sigma/q} \norm[\D^{1,q}(\R^N;|x|^{-b})]{u} \le \max \bgset{t^{\sigma/p},t^{\sigma/q}} \norm{u}
\end{equation}
for all $u \in W$ and $t \ge 0$.

We verify that this scaling satisfies the assumptions of the abstract framework introduced in Section \ref{Section abstract}.

\begin{proposition} 
	The mapping
	\begin{align*} 
		W_\rad \times [0,\infty) \to W_\rad, \quad (u,t) \mapsto u_t,
	\end{align*}
	is a scaling in $W_\rad$ in the sense of Definition \ref{Definition scaling} with $\sigma$ as in \eqref{rho sigma}.
\end{proposition}

\begin{proof}
	Assumptions \textnormal{(\hyperlink{A1}{A$_1$})}--\textnormal{(\hyperlink{A3}{A$_3$})} are clearly satisfied, while \textnormal{(\hyperlink{A4}{A$_4$})} and \textnormal{(\hyperlink{A5}{A$_5$})} follow from \eqref{eq.1.prop1}. It only remains to verify continuity. Note that $(u,t) \mapsto u_t$ is linear in $u$. Let $(u_j,t_j) \to (u,t)$ in $W_\rad \times [0,\infty)$. If $t = 0$, then $u_t = 0$ and $(u_j)_{t_j} \to 0$ by \eqref{eq.1.prop1} since $\seq{u_j}$ is bounded, so suppose $t > 0$. Then
	\begin{align*}
		\norm{(u_j)_{t_j} - u_t} \le \norm{(u_j)_{t_j} - u_{t_j}} + \norm{u_{t_j} - u_t},
	\end{align*}
	and \eqref{eq.1.prop1} together with \textnormal{(\hyperlink{A1}{A$_1$})} gives
	\begin{align*}
		\norm{(u_j)_{t_j} - u_{t_j}} = \norm{(u_j - u)_{t_j}} \le \max \bgset{t_j^{\sigma/p},t_j^{\sigma/q}} \norm{u_j - u} \to 0
	\end{align*}
	and
	\begin{align*}
		\norm{u_{t_j} - u_t} = \norm{(u_{t_j/t})_t - u_t} = \norm{(u_{t_j/t} - u)_t} \le \max \set{t^{\sigma/p},t^{\sigma/q}} \norm{u_{t_j/t} - u}.
	\end{align*}
	Thus, it suffices to show that $u_{\tau_j} \to u$ in $W$ whenever $\tau_j \to 1$. Moreover, we may assume that $u \in C^\infty_c(\R^N)$ since $W$ is the completion of $C^\infty_c(\R^N)$. Indeed, for $\varphi \in C^\infty_c(\R^N)$, \eqref{eq.1.prop1} gives
	\begin{align*}
		\norm{u_{\tau_j} - u}
		&\le \norm{(u - \varphi)_{\tau_j}} + \norm{\varphi_{\tau_j} - \varphi} + \norm{\varphi - u}\\
		&\le \big(\max \bgset{\tau_j^{\sigma/p},\tau_j^{\sigma/q}} + 1\big) \norm{u - \varphi} + \norm{\varphi_{\tau_j} - \varphi}.
	\end{align*}
	The first term on the right-hand side can be made arbitrarily small, uniformly for sufficiently large $j$, by choosing $\varphi$ sufficiently close to $u$ in $W$.

	Now let $u \in C^\infty_c(\R^N)$ with $\supp u \subset B_R$ and let $\tau_j \to 1$. For all sufficiently large $j$, $\tau_j \in [1/2,2]$. Then $\supp u_{\tau_j} \subset B_{2R}$ and
	\begin{align*}
		\nabla u_{\tau_j}(x) = \tau_j^\alpha\, \nabla u(x/\tau_j) \to \nabla u(x)
	\end{align*}
	uniformly on $\R^N$, by the uniform continuity of $\nabla u$. Since $b < N$ and hence $\int_{B_{2R}} |x|^{-b}\, \mathrm{d}x < \infty$, it follows that
	\begin{align*}
		\int_{\R^N} |\nabla (u_{\tau_j} - u)|^p\, \mathrm{d}x \to 0, \quad \int_{\R^N} |x|^{-b}\, |\nabla (u_{\tau_j} - u)|^q\, \mathrm{d}x \to 0.
	\end{align*}
	Hence, $u_{\tau_j} \to u$ in $W$.
\end{proof}

The scaled operator associated with the left-hand side of equation \eqref{1} and its potential are given by
\begin{equation} \label{As}
	\As(u) v = \int_{\R^N} \left(|\nabla u|^{p-2} + |x|^{-b}\, |\nabla u|^{q-2}\right) \nabla u \cdot \nabla v\, \mathrm{d}x, \quad u, v \in W_\rad
\end{equation}
and
\begin{equation} \label{I_sigma}
	I_\sigma(u) = \frac{1}{p} \int_{\R^N} |\nabla u|^p\, \mathrm{d}x + \frac{1}{q} \int_{\R^N} |x|^{-b}\, |\nabla u|^q\, \mathrm{d}x, \quad u \in W_\rad,
\end{equation}
respectively. A change of variables shows that $\As(u_t) v_t = t^\sigma \As(u) v$ for all $u, v \in W_\rad$ and $t \ge 0$. Thus, $\As$ is a scaled operator in the sense of Definition \ref{Definition scaled operator}. Assumptions \textnormal{(\hyperlink{A6}{A$_6$})} and \textnormal{(\hyperlink{A10}{A$_{10}$})} of Section \ref{Scaled eigenvalue problems} are immediate. We verify \textnormal{(\hyperlink{A7}{A$_7$})} in the next proposition.

\begin{proposition} 
	Every sequence $\seq{u_j}$ in $W_\rad$ such that $u_j \wto u$ and
	\begin{equation} \label{eq:S+hp}
		\limsup_{j \to \infty} \int_{\R^N} \left(|\nabla u_j|^{p-2} + |x|^{-b}\, |\nabla u_j|^{q-2}\right) \nabla u_j \cdot \nabla (u_j - u)\, \mathrm{d}x \le 0
	\end{equation}
	converges strongly to $u$. In particular, \textnormal{(\hyperlink{A7}{A$_7$})} holds.
\end{proposition}

\begin{proof}
	We adapt the argument of Perera, Agarwal, and O'Regan \cite[Proposition 1.3]{Perera-Agarwal-ORegan-2010} to each of the two terms of $\As$. Write $\As = T_1 + T_2$, where
	\begin{align*}
		T_1(u)\, v = \int_{\R^N} |\nabla u|^{p-2}\, \nabla u \cdot \nabla v\, \mathrm{d}x, \quad T_2(u)\, v = \int_{\R^N} |x|^{-b}\, |\nabla u|^{q-2}\, \nabla u \cdot \nabla v\, \mathrm{d}x.
	\end{align*}
	Let $\ell_1 = \norm[\D^{1,p}(\R^N)]{\cdot}$ and $\ell_2 = \norm[\D^{1,q}(\R^N;|x|^{-b})]{\cdot}$ and set $r_1 = p$ and $r_2 = q$. H\"{o}lder's inequality implies
	\begin{align*}
		T_i(u)\, v \le \ell_i(u)^{r_i - 1}\, \ell_i(v), \quad T_i(u)\, u = \ell_i(u)^{r_i}
	\end{align*}
	for all $u, v \in W_\rad$ and $i = 1, 2$. Therefore,
	\begin{align*}
		D_{i,j}
		&:= T_i(u_j)(u_j - u) - T_i(u)(u_j - u)\\
		&= T_i(u_j)\, u_j - T_i(u_j)\, u - T_i(u)\, u_j + T_i(u)\, u\\
		&\ge \left(\ell_i(u_j)^{r_i - 1} - \ell_i(u)^{r_i - 1}\right)\!\big(\ell_i(u_j) - \ell_i(u)\big) \ge 0.
	\end{align*}
	Since $u_j \wto u$ and $\As(u) \in W_\rad^\ast$, $\As(u)(u_j - u) \to 0$. Combining this with \eqref{eq:S+hp}, we obtain
	\begin{align*}
		\limsup_{j \to \infty}\, (D_{1,j} + D_{2,j}) = \limsup_{j \to \infty}\, \As(u_j)(u_j - u) \le 0.
	\end{align*}
	As both terms are nonnegative, $D_{1,j} \to 0$ and $D_{2,j} \to 0$, so $\ell_1(u_j) \to \ell_1(u)$ and $\ell_2(u_j) \to \ell_2(u)$.

	The linear maps $W_\rad \to L^p(\R^N;\R^N)$, $u \mapsto \nabla u$ and $W_\rad \to L^q(\R^N;\R^N)$, $u \mapsto |x|^{-b/q}\, \nabla u$ are continuous. Thus, $u_j \wto u$ implies $\nabla u_j \wto \nabla u$ in $L^p(\R^N;\R^N)$ and $|x|^{-b/q}\, \nabla u_j \wto |x|^{-b/q}\, \nabla u$ in $L^q(\R^N;\R^N)$. Since the corresponding norms converge and both Lebesgue spaces are uniformly convex, the convergence is strong in each space. Consequently, $u_j \to u$ in $W_\rad$.
\end{proof}

\subsection{Poho\v{z}aev identity}

In this subsection and the next, we consider nonlinearities of the form
\begin{equation} \label{power f}
	f(r,t) = \sum_{i=1}^n \lambda_i\, r^{-c_i}\, |t|^{s_i-2}\, t, \quad r > 0,\, t \in \R,
\end{equation}
where $\lambda_i \in \R$, $(s_i,c_i) \in \Te$, and $s_i > 1$. For such $f$, the energy functional takes the form
\begin{align*}
	E(u) = I_\sigma(u) - \sum_{i=1}^n \frac{\lambda_i}{s_i}\, J_i(u), \quad J_i(u) = \int_{\R^N} |x|^{-c_i}\, |u|^{s_i}\, \mathrm{d}x, \quad u \in W_\rad,
\end{align*}
which is of class $C^1$ by Theorem \ref{Theorem embedding} \textnormal{(\hyperlink{Theorem embedding.i}{i})}. We write $\rho_i = \rho(s_i,c_i)$, so that $J_i(u_t) = t^{\rho_i} J_i(u)$ by \eqref{scaling relations}.

In the proof of the Poho\v{z}aev identity, we will use a regularization of the radial integrand of $I_\sigma$ in the spirit of Guedda and V\'{e}ron \cite[Proof of Theorem 1.1]{Guedda-Veron-1989}. For $\eps \in [0,1]$, $r > 0$, and $\xi \in \R$, let
\begin{align*}
	a_\eps(r,\xi) = (\eps + \xi^2)^{(p-2)/2}\, \xi + r^{-b}\, (\eps + \xi^2)^{(q-2)/2}\, \xi
\end{align*}
and
\begin{align*}
	A_\eps(r,\xi) = \frac{(\eps + \xi^2)^{p/2} - \eps^{p/2}}{p} + r^{-b}\, \frac{(\eps + \xi^2)^{q/2} - \eps^{q/2}}{q}.
\end{align*}
Then $\partial_\xi\,  A_\eps = a_\eps$. For $\varepsilon=0$, we write
\begin{align*}
	a(r,\xi) &= a_0(r,\xi) = |\xi|^{p-2}\, \xi + r^{-b}\, |\xi|^{q-2}\, \xi,\\
	A(r,\xi) &= A_0(r,\xi) = \frac{|\xi|^p}{p} + r^{-b}\, \frac{|\xi|^q}{q}.
\end{align*}
The functions $a_\eps(r,\xi)$, $A_\eps(r,\xi)$, and
\begin{align*}
	\partial_r\,  A_\eps(r,\xi) = - \frac{b}{q}\, r^{-b-1} \left[(\eps + \xi^2)^{q/2} - \eps^{q/2}\right]
\end{align*}
are continuous in $(\eps,r,\xi) \in [0,1] \times (0,\infty) \times \R$. For each $\eps \in [0,1]$ and $r > 0$, the function $a_\eps(r,\cdot)$ is odd, continuous, and strictly increasing from $\R$ onto $\R$. We denote its inverse by $\widetilde{a}_\eps(r,\cdot)$ and set $\widetilde{a} = \widetilde{a}_0$.

The following lemma provides the regularity properties of the inverse functions needed in the proof.

\begin{lemma} \label{Lemma regularization}
	The following assertions hold.
	\begin{enumerate}
		\item[\textnormal{(i)}]\hypertarget{Lemma regularization.i}{}
			The mapping $(\eps,r,\eta) \mapsto \widetilde{a}_\eps(r,\eta)$ is continuous on $[0,1] \times (0,\infty) \times \R$.
		\item[\textnormal{(ii)}]\hypertarget{Lemma regularization.ii}{}
			For every $\eps \in (0,1]$, $\widetilde{a}_\eps \in C^1((0,\infty) \times \R)$.
	\end{enumerate}
\end{lemma}

\begin{proof}
	\textnormal{(\hyperlink{Lemma regularization.i}{i})}
	For $\eps \in [0,1]$, $r > 0$, and $\xi \in \R$, we have
	\begin{equation} \label{eq:a_eps-lower-bd}
		|a_\eps(r,\xi)| \ge (\eps + \xi^2)^{(p-2)/2}\, |\xi| \ge \theta(|\xi|),
	\end{equation}
	where
	\begin{equation} \label{eq:theta-coercive}
		\theta(t) = \min\, \{t^{p-1},(1 + t^2)^{(p-2)/2}\, t\} \to \infty \quad \text{as } t \to \infty.
	\end{equation}
	Let $(\eps_n,r_n,\eta_n) \to (\eps,r,\eta)$ in $[0,1] \times (0,\infty) \times \R$ and set $\xi_n = \widetilde{a}_{\eps_n}(r_n,\eta_n)$. Then $a_{\eps_n}(r_n,\xi_n) = \eta_n$ and hence $\theta(|\xi_n|) \le |\eta_n|$ by \eqref{eq:a_eps-lower-bd}, so $\seq{\xi_n}$ is bounded by \eqref{eq:theta-coercive}. If $\xi'$ is a limit point of $\seq{\xi_n}$, then $a_\eps(r,\xi') = \lim a_{\eps_n}(r_n,\xi_n) = \eta$ by the continuity of $(\eps,r,\xi) \mapsto a_\eps(r,\xi)$, so $\xi' = \widetilde{a}_\eps(r,\eta)$. Hence $\xi_n \to \widetilde{a}_\eps(r,\eta)$.

	\textnormal{(\hyperlink{Lemma regularization.ii}{ii})}
	For $\eps > 0$, the function $a_\eps$ belongs to $C^\infty((0,\infty) \times \R)$. Moreover,
	\begin{align*}
		\partial_\xi\,  a_\eps(r,\xi)
		& = (\eps + \xi^2)^{(p-4)/2} \left[\eps + (p - 1)\, \xi^2\right] \\
		& \quad + r^{-b}\, (\eps + \xi^2)^{(q-4)/2} \left[\eps + (q - 1)\, \xi^2\right] > 0.
	\end{align*}
	The conclusion follows from the inverse function theorem.
\end{proof}

We are now ready to prove the following Poho\v{z}aev identity.

\begin{theorem} \label{thm:Pohozaev}
	Let $f$ be as in \eqref{power f}. Then every weak radial solution $u \in W_\rad$ of equation \eqref{1} satisfies the Poho\v{z}aev identity
	\begin{equation} \label{Pohozaev}
		\sigma I_\sigma(u) = \sum_{i=1}^n \frac{\rho_i\,  \lambda_i}{s_i}\, J_i(u).
	\end{equation}
\end{theorem}

\begin{proof}
	Let $u(x) = u(r)$, where $r = |x|$, be the representative given by Lemma \ref{Lemma Strauss}. Set
	\begin{align*}
		F_i(r,t) = \frac{\lambda_i}{s_i}\, r^{-c_i}\, |t|^{s_i}, \quad F(r,t) = \sum_{i=1}^n F_i(r,t).
	\end{align*}
	Since $s_i > 1$ for each $i$, both $f$ and $F$ are continuous on $(0,\infty) \times \R$. Moreover,
	\begin{equation} \label{eq:tf}
		t\,  f(r,t) = \sum_{i=1}^n s_i\, F_i(r,t).
	\end{equation}
	For $v \in C^\infty_c(0,\infty)$, the radial function $v(|x|)$ belongs to $W_\rad$. By Lemma \ref{Lemma Strauss} \textnormal{(\hyperlink{Lemma Strauss.iii}{iii})}, the weak formulation of equation \eqref{1}, tested with $v(|x|)$, becomes
	\begin{equation} \label{radial weak}
		\int_0^\infty \Phi(r)\, v'(r)\,\mathrm{d}r = \int_0^\infty r^{N-1}\, f(r,u(r))\, v(r)\,\mathrm{d}r,
	\end{equation}
	where $\Phi(r) = r^{N-1}\, a(r,u'(r))$. We have $\Phi \in L^1_\loc(0,\infty)$, as $u' \in L^p_\loc(0,\infty) \cap L^q_\loc(0,\infty)$.

	The continuity of $u$ on $(0,\infty)$ implies that $r \mapsto r^{N-1}\, f(r,u(r))$ is continuous. Hence the function
	\begin{align*}
		\Lambda(r) = - \int_1^r \tau^{N-1}\, f(\tau,u(\tau))\,\mathrm{d}\tau
	\end{align*}
	belongs to $C^1(0,\infty)$ with $\Lambda'(r) = - r^{N-1}\, f(r,u(r))$. Integrating by parts on the right-hand side of \eqref{radial weak} gives
	\begin{align*}
		\int_0^\infty (\Phi(r) - \Lambda(r))\, v'(r)\,\mathrm{d}r = 0 \quad \text{for all }v \in C^\infty_c(0,\infty),
	\end{align*}
	so $\Phi - \Lambda$ is equal to a constant a.e.\,by the du Bois-Reymond lemma. Hence, after modification on a null set, $\Phi \in C^1(0,\infty)$ and
	\begin{equation} \label{Phi prime}
		\Phi'(r) = - r^{N-1}\, f(r,u(r)) \quad \text{for all }r > 0.
	\end{equation}
	Then $u'(r) = \widetilde{a}(r,r^{1-N}\, \Phi(r))$ for a.e.\,$r > 0$ and the right-hand side is continuous by Lemma \ref{Lemma regularization} \textnormal{(\hyperlink{Lemma regularization.i}{i})}, so $u \in C^1(0,\infty)$ and $a(r,u'(r)) = r^{1-N}\, \Phi(r)$ for all $r > 0$.

	Define
	\begin{align*}
		P(r) = r^N\, \big[a(r,u'(r))\, u'(r) - A(r,u'(r))\big] + r^N\, F(r,u(r)) - (\alpha + 1)\, \Phi(r)\, u(r),
	\end{align*}
	for $r>0$. Clearly, $P$ is continuous on $(0,\infty)$. We claim that, for all $0 < \delta < R$,
	\begin{equation} \label{P difference}
		P(R) - P(\delta) = \int_\delta^R r^{N-1} \left[\sum_{i=1}^n \rho_i\, F_i(r,u(r)) - \sigma A(r,u'(r))\right] \,\mathrm{d}r.
	\end{equation}
	If $u$ were of class $C^2$, this identity would follow directly by differentiating $P$. Since we only know that $u \in C^1(0,\infty)$, we use the regularization as follows. For $\eps \in (0,1]$, set
	\begin{align*}
		v_\eps(r) = \widetilde{a}_\eps(r,r^{1-N}\, \Phi(r)), \quad r > 0.
	\end{align*}
	By Lemma \ref{Lemma regularization} \textnormal{(\hyperlink{Lemma regularization.ii}{ii})}, $v_\eps \in C^1(0,\infty)$, and its definition implies
	\begin{equation} \label{frozen flux}
		r^{N-1}\, a_\eps(r,v_\eps(r)) = \Phi(r), \quad r > 0.
	\end{equation}
	Now define
	\begin{align*}
		P_\eps(r)
		& = r^N\, \big[a_\eps(r,v_\eps(r))\, v_\eps(r) - A_\eps(r,v_\eps(r))\big] \\
		& \quad + r^N\, F(r,u(r)) - (\alpha + 1)\, \Phi(r)\, u(r), \quad r > 0.
	\end{align*}
	By \eqref{frozen flux}, the first term can be written as $r\, \Phi(r)\, v_\eps(r) - r^N\, A_\eps(r,v_\eps(r))$, and since $\partial_\xi\,  A_\eps = a_\eps$, \eqref{frozen flux} also gives
	\begin{align*}
		\frac{\mathrm{d}}{\mathrm{d}r} \left[r^N\, A_\eps(r,v_\eps)\right] = Nr^{N-1}\, A_\eps(r,v_\eps(r)) + r^N\, \partial_r\,  A_\eps(r,v_\eps(r)) + r\, \Phi(r)\, v_\eps'(r).
	\end{align*}
	So the terms containing $v_\eps'(r)$ cancel. Using \eqref{frozen flux} and \eqref{Phi prime}, we find
	\begin{equation} \label{DP_eps-1}
		\begin{aligned}
			&\frac{\mathrm{d}}{\mathrm{d}r} \left[r\, \Phi(r)\, v_\eps(r) - r^N\, A_\eps(r,v_\eps(r))\right]\\
			&= r^{N-1}\, \big[a_\eps(r,v_\eps(r))\, v_\eps(r) - r\,  f(r,u(r))\, v_\eps(r)\\
			&\quad - N A_\eps(r,v_\eps(r)) - r\, \partial_r\,  A_\eps(r,v_\eps(r))\big].
		\end{aligned}
	\end{equation}
	Next, the chain rule shows that
	\begin{equation} \label{Dp_eps-2}
		\frac{\mathrm{d}}{\mathrm{d}r} \left[r^N\, F(r,u(r))\right] = r^{N-1} \left[\sum_{i=1}^n (N - c_i)\, F_i(r,u(r)) + r\,  f(r,u(r))\, u'(r)\right],
	\end{equation}
	because $r\, \partial_r\,  F_i(r,t) = - c_i\, F_i(r,t)$ and $\partial_t\,  F(r,t) = f(r,t)$. Furthermore, \eqref{Phi prime} and \eqref{eq:tf} imply
	\begin{equation} \label{Dp_eps-3}
		\begin{aligned}
			&\frac{\mathrm{d}}{\mathrm{d}r} \left[- (\alpha + 1)\, \Phi(r)\, u(r)\right]\\
			&= (\alpha + 1)\, r^{N-1} \left[\sum_{i=1}^n s_i\, F_i(r,u(r)) - a(r,u'(r))\, u'(r)\right].
		\end{aligned}
	\end{equation}
	Adding \eqref{DP_eps-1}, \eqref{Dp_eps-2}, and \eqref{Dp_eps-3}, and recalling that $\rho_i = (\alpha + 1)\, s_i - c_i + N$, we conclude that $P_\eps \in C^1(0,\infty)$ with
	\begin{align*}
		P_\eps'(r) = r^{N-1}\, \Bigg[&\sum_{i=1}^n \rho_i\, F_i(r,u(r)) + a_\eps(r,v_\eps(r))\, v_\eps(r)\\
		&\quad- N A_\eps(r,v_\eps(r)) - r\, \partial_r\,  A_\eps(r,v_\eps(r))\\
		&\quad- (\alpha + 1)\, a(r,u'(r))\, u'(r) + r\,  f(r,u(r))\, (u'(r) - v_\eps(r))\Bigg].
	\end{align*}
	Now fix $0 < \delta < R$. By Lemma \ref{Lemma regularization} \textnormal{(\hyperlink{Lemma regularization.i}{i})}, the function $(\eps,r) \mapsto \widetilde{a}_\eps(r,r^{1-N}\, \Phi(r))$ is uniformly continuous on $[0,1] \times [\delta,R]$. Therefore, $v_\eps(r) \to \widetilde{a}(r,r^{1-N}\, \Phi(r)) = u'(r)$ uniformly on $[\delta,R]$ as $\eps \to 0$. Since $a_\eps(r,\xi)$, $A_\eps(r,\xi)$, and $\partial_r\,  A_\eps(r,\xi)$ are uniformly continuous on compact subsets of $[0,1] \times (0,\infty) \times \R$, it follows that $P_\eps(r) \to P(r)$ and
	\begin{align*}
		P_\eps'(r) \to r^{N-1} \bigg[
		&\sum_{i=1}^n \rho_i\, F_i(r,u(r)) - \alpha\,  a(r,u'(r))\, u'(r)\\
		&\quad - N A(r,u'(r)) - r\, \partial_r\,  A(r,u'(r))\bigg]
	\end{align*}
	uniformly on $[\delta,R]$. Since
	\begin{align*}
		&\alpha\,  a(r,u'(r))\, u'(r) + N A(r,u'(r)) + r\, \partial_r\,  A(r,u'(r))\\
		&= (\alpha p + N)\, \frac{|u'(r)|^p}{p} + (\alpha q + N - b)\, r^{-b}\, \frac{|u'(r)|^q}{q} = \sigma A(r,u'(r))
	\end{align*}
	by \eqref{rho sigma}, now passing to the limit in $P_\eps(R) - P_\eps(\delta) = \int_\delta^R P_\eps'(r)\,\mathrm{d}r$ proves \eqref{P difference}.

	We have
	\begin{align*}
		0 \le A(r,u'(r)) \le a(r,u'(r))\, u'(r) = |u'(r)|^p + r^{-b}\, |u'(r)|^q
	\end{align*}
	and
	\begin{align*}
		\frac{|\Phi(r)\, u(r)|}{r^N}
		& \le |u'(r)|^{p-1}\, \frac{|u(r)|}{r} + r^{-b}\, |u'(r)|^{q-1}\, \frac{|u(r)|}{r} \\
		& \le |u'(r)|^p + r^{-p}\, |u(r)|^p + r^{-b}\, |u'(r)|^q+ r^{-q-b}\, |u(r)|^q
	\end{align*}
	by Young's inequality. Consequently,
	\begin{align*}
		\frac{|P(r)|}{r} \le Cr^{N-1} \bigg[&|u'(r)|^p + r^{-b}\, |u'(r)|^q + \sum_{i=1}^n r^{-c_i}\, |u(r)|^{s_i}\\
		&\quad + r^{-p}\, |u(r)|^p + r^{-q-b}\, |u(r)|^q\bigg]
	\end{align*}
	for some constant $C > 0$. Hence $\int_0^\infty |P(r)|\,\mathrm{d}r/r < \infty$ by Lemma \ref{Lemma Strauss} \textnormal{(\hyperlink{Lemma Strauss.iii}{iii})}, Theorem \ref{Theorem embedding} \textnormal{(\hyperlink{Theorem embedding.i}{i})}, and the Hardy and weighted Hardy inequalities. So there are sequences $\delta_k \searrow 0$ and $R_k \nearrow \infty$ such that $P(\delta_k) \to 0$ and $P(R_k) \to 0$. Since the integrand in \eqref{P difference} is integrable on $(0,\infty)$ for the same reasons, taking $\delta = \delta_k$ and $R = R_k$ there and letting $k \to \infty$, we obtain
	\begin{align*}
		\sigma \int_0^\infty r^{N-1}\, A(r,u'(r))\,\mathrm{d}r = \sum_{i=1}^n \rho_i \int_0^\infty r^{N-1}\, F_i(r,u(r))\,\mathrm{d}r,
	\end{align*}
	which is \eqref{Pohozaev} after multiplying by $\omega_{N-1}$.
\end{proof}

The Poho\v{z}aev identity has the following immediate consequences.

\begin{corollary} 
	Let $f$ be as in \eqref{power f}. Then every weak radial solution $u \in W_\rad$ of equation \eqref{1} satisfies the Nehari identity
	\begin{equation} \label{Nehari}
		\int_{\R^N} |\nabla u|^p\, \mathrm{d}x + \int_{\R^N} |x|^{-b}\, |\nabla u|^q\, \mathrm{d}x = \sum_{i=1}^n \lambda_i \int_{\R^N} |x|^{-c_i}\, |u|^{s_i}\, \mathrm{d}x,
	\end{equation}
	the dilation identity
	\begin{equation} \label{dilation}
		\begin{aligned}
			&\frac{N - p}{p} \int_{\R^N} |\nabla u|^p\, \mathrm{d}x + \frac{N - q - b}{q} \int_{\R^N} |x|^{-b}\, |\nabla u|^q\, \mathrm{d}x\\
			&= \sum_{i=1}^n \frac{(N - c_i)\, \lambda_i}{s_i} \int_{\R^N} |x|^{-c_i}\, |u|^{s_i}\, \mathrm{d}x,
		\end{aligned}
	\end{equation}
	and the energy identity
	\begin{equation} \label{energy identity}
		E(u) = \sum_{i=1}^n \left(\frac{\rho_i}{\sigma} - 1\right) \frac{\lambda_i}{s_i} \int_{\R^N} |x|^{-c_i}\, |u|^{s_i}\, \mathrm{d}x.
	\end{equation}
\end{corollary}

\begin{proof}
	The Nehari identity \eqref{Nehari} is $E'(u)\, u = 0$. Since $\sigma/p - (\alpha + 1) = (N - p)/p$ and $\sigma/q - (\alpha + 1) = (N - q - b)/q$ by \eqref{rho sigma} and $\rho_i - (\alpha + 1)\, s_i = N - c_i$, subtracting $\alpha + 1$ times \eqref{Nehari} from \eqref{Pohozaev} gives \eqref{dilation}. Dividing \eqref{Pohozaev} by $\sigma$ and subtracting from the definition of $E(u)$ gives \eqref{energy identity}.
\end{proof}

\subsection{Weak limits of Palais-Smale sequences}

Recall that a sequence $\seq{u_j} \subset W_\rad$ is a \PS{c} sequence for $E$ if $E(u_j) \to c$ and $E'(u_j) \to 0$ in $W_\rad^\ast$. The functional $E$ is said to satisfy the \PS{c} condition if every \PS{c} sequence has a strongly convergent subsequence.

We conclude this section by identifying the weak limits of bounded sequences of approximate critical points.

\begin{lemma} \label{Lemma weak limit}
	Let $f$ be as in \eqref{power f}. Suppose that $\seq{u_j}$ is a bounded sequence in $W_\rad$ such that $u_j \wto u$ and $E'(u_j) \to 0$. Then, for a renamed subsequence, $u_j \to u$ uniformly on every annulus $A_{\rho,R} = \set{x \in \R^N \colon\rho \le |x| \le R}$, $\nabla u_j \to \nabla u$ a.e.\,in $\R^N$, and $u$ is a weak radial solution of equation \eqref{1}.
\end{lemma}

\begin{proof}
	Set $w_j = u_j - u$. By the proof of Theorem \ref{Theorem embedding} \textnormal{(\hyperlink{Theorem embedding.ii}{ii})}, $w_j \to 0$ uniformly on every annulus $A_{\rho,R}$. Let $\varphi \in C^\infty_c(\R^N \setminus \set{0})$ be radial with $0 \le \varphi \le 1$, and choose an annulus $A = A_{\rho,R}$ containing $\supp \varphi$. For any radial $w \in C^\infty_c(\R^N)$, $\nabla (\varphi w) = \varphi \nabla w + w \nabla \varphi$ and $\sup_A |w| \le C \norm{w}$ by Lemma \ref{Lemma Strauss} \textnormal{(\hyperlink{Lemma Strauss.i}{i})}. Thus, $\norm{\varphi w} \le C_\varphi \norm{w}$ since the weight $|x|^{-b}$ is bounded on $A$. By Lemma \ref{Lemma averaging} \textnormal{(\hyperlink{Lemma averaging.iii}{iii})}, this extends to all $w \in W_\rad$, with $\varphi w \in W_\rad$. In particular, $\seq{\varphi w_j}$ is bounded in $W_\rad$ and hence $E'(u_j)(\varphi w_j) \to 0$. Writing $a(x,\xi) = (|\xi|^{p-2} + |x|^{-b}\, |\xi|^{q-2})\, \xi$, we have
	\begin{align*}
		E'(u_j)(\varphi w_j)
		&= \int_{\R^N} \varphi\, a(x,\nabla u_j) \cdot \nabla w_j\, \mathrm{d}x + \int_{\R^N} w_j\, a(x,\nabla u_j) \cdot \nabla \varphi\, \mathrm{d}x\\
		&\quad- \int_{\R^N} f(|x|,u_j)\, \varphi w_j\, \mathrm{d}x.
	\end{align*}
	The last two integrals go to zero since $\sup_A |w_j| \to 0$, while $a(x,\nabla u_j)$ and $f(|x|,u_j)$ are bounded in $L^1(A)$ by H\"{o}lder's inequality, the boundedness of the weights on $A$ and of $\seq{u_j}$ in $W_\rad$, and Lemma \ref{Lemma Strauss} \textnormal{(\hyperlink{Lemma Strauss.i}{i})}. Since $v \mapsto \int_{\R^N} \varphi\, a(x,\nabla u) \cdot \nabla v\, \mathrm{d}x$ belongs to $W_\rad^\ast$ and $w_j \wto 0$, it follows that
	\begin{align*}
		\int_{\R^N} \varphi \left(D_{1,j} + D_{2,j}\right)\, \mathrm{d}x \to 0,
	\end{align*}
	where
	\begin{align*}
		D_{1,j} & = \big(|\nabla u_j|^{p-2}\, \nabla u_j - |\nabla u|^{p-2}\, \nabla u\big) \cdot \nabla w_j, \\
		D_{2,j} & = |x|^{-b}\, \big(|\nabla u_j|^{q-2}\, \nabla u_j - |\nabla u|^{q-2}\, \nabla u\big) \cdot \nabla w_j
	\end{align*}
	are nonnegative. So $D_{1,j} \to 0$ in $L^1(\set{x \in \R^N \colon\varphi(x) = 1})$, and hence a.e.\,in $\set{x \in \R^N \colon\varphi(x) = 1}$ for a renamed subsequence. Therefore,
	\begin{align*}
		\nabla u_j \to \nabla u\quad \text{a.e.\,in }\left\{x \in \R^N \colon\varphi(x) = 1\right\}.
	\end{align*}
	Taking a sequence of such $\varphi$ with $\set{x \in \R^N \colon\varphi (x)= 1}$ exhausting $\R^N \setminus \set{0}$ and passing to a diagonal subsequence gives $\nabla u_j \to \nabla u$ a.e.\,in $\R^N$.

	Finally, the sequences
	\begin{align*}
		|\nabla u_j|^{p-2}\, \nabla u_j,\quad |x|^{-b/q'}\, |\nabla u_j|^{q-2}\, \nabla u_j,\quad\text{and}\quad |x|^{-c_i/s_i'}\, |u_j|^{s_i-2}\, u_j
	\end{align*}
	for $i=1,\ldots, n$, are bounded in
	\begin{align*}
		L^{p'}(\R^N),\quad L^{q'}(\R^N), \quad\text{and}\quad L^{s_i'}(\R^N),
	\end{align*}
	respectively, where the primes denote H\"{o}lder conjugates. Since $u_j \to u$ and $\nabla u_j \to \nabla u$ a.e., these sequences converge pointwise a.e.\,to the corresponding expressions involving $u$. Their boundedness in the respective reflexive Lebesgue spaces implies weak convergence to these limits. Pairing them with $\nabla v$, $|x|^{-b/q}\, \nabla v$, and $|x|^{-c_i/s_i}\, v$, respectively, we find $E'(u)\, v = \lim E'(u_j)\, v = 0$ for all $v \in W_\rad$. Thus, $u$ is a weak radial solution of equation \eqref{1}.
\end{proof}

\section{Scaling-based methods in critical point theory} \label{Section abstract}

In this section, we recall the construction of eigenvalues for scaled eigenvalue problems using the cohomological index, together with several abstract critical point theorems from Mercuri and Perera \cite{Mercuri-Perera-2026} and Perera \cite{Perera-2026}.

\subsection{Scalings}

Let $W$ be an infinite-dimensional reflexive Banach space. The following notion of a scaling in $W$ was introduced by Mercuri and Perera \cite[Definition 1.1]{Mercuri-Perera-2026}.

\begin{definition}\label{Definition scaling}
	A scaling in $W$ is a continuous mapping
	\begin{align*}
		W \times [0,\infty) \to W, \quad (u,t) \mapsto u_t
	\end{align*}
	satisfying
	\begin{enumerate}
		\item[\textnormal{(A$_1$)}]\hypertarget{A1}{}
			$(u_{t_1})_{t_2} = u_{t_1 t_2}$ for all $u \in W$ and $t_1, t_2 \ge 0$.
		\item[\textnormal{(A$_2$)}]\hypertarget{A2}{}
			$(\tau u)_t = \tau u_t$ for all $u \in W$, $\tau \in \R$, and $t \ge 0$.
		\item[\textnormal{(A$_3$)}]\hypertarget{A3}{}
			$u_0 = 0$ and $u_1 = u$ for all $u \in W$.
		\item[\textnormal{(A$_4$)}]\hypertarget{A4}{}
			$u_t$ is bounded on bounded sets in $W \times [0,\infty)$.
		\item[\textnormal{(A$_5$)}]\hypertarget{A5}{}
			There exists $\sigma > 0$ such that $\norm{u_t} = \O(t^\sigma)$ as $t \to \infty$, uniformly in $u$ on bounded sets of $W$.
	\end{enumerate}
\end{definition}

\subsection{Scaled operators}

Denote by $W^\ast$ the dual of $W$. Recall that $q \in C(W,W^\ast)$ is a potential operator if there exists a functional $Q \in C^1(W,\R)$, called a potential for $q$, whose Fr\'{e}chet derivative satisfies $Q' = q$. Replacing $Q$ with $Q - Q(0)$ if necessary, we may assume that $Q(0) = 0$. The next definition is taken from Mercuri and Perera \cite[Definition 2.1]{Mercuri-Perera-2026}.

\begin{definition}\label{Definition scaled operator}
	A scaled operator is an odd potential operator $\As \in C(W,W^\ast)$ that maps bounded sets into bounded sets and satisfies
	\begin{align*}
		\As(u_t) v_t = t^\sigma \As(u) v \quad \text{for all }u, v \in W,\, t \ge 0.
	\end{align*}
\end{definition}

\subsection{Scaled eigenvalue problems} \label{Scaled eigenvalue problems}

Let $\As$ and $\Bs$ be scaled operators satisfying:
\begin{enumerate}
	\item[\textnormal{(A$_6$)}]\hypertarget{A6}{}
		$\As(u) u > 0$ for all $u \in W \setminus \set{0}$.
	\item[\textnormal{(A$_7$)}]\hypertarget{A7}{}
		Every sequence $\seq{u_j}$ in $W$ such that $u_j \wto u$ and $\As(u_j)(u_j - u) \to 0$ has a subsequence that converges strongly to $u$.
	\item[\textnormal{(A$_8$)}]\hypertarget{A8}{}
		$\Bs(u) u > 0$ for all $u \in W \setminus \set{0}$.
	\item[\textnormal{(A$_9$)}]\hypertarget{A9}{}
		if $u_j \wto u$ in $W$, then $\Bs(u_j) \to \Bs(u)$ in $W^\ast$.
\end{enumerate}
The scaled eigenvalue problem
\begin{equation} \label{8}
	\As(u) = \lambda \Bs(u) \quad \text{in } W^\ast
\end{equation}
was studied by Mercuri and Perera \cite{Mercuri-Perera-2026}. We call $\lambda \in \R$ an eigenvalue of this problem if there is a $u \in W \setminus \set{0}$, called an eigenfunction associated with $\lambda$, satisfying \eqref{8}. By the scaling property, $u_t$ is also an eigenfunction associated with $\lambda$ for any $t > 0$ since
\begin{align*}
	\As(u_t) v
	& = \As(u_t) (v_{t^{-1}})_t = t^\sigma \As(u) v_{t^{-1}} = \lambda t^\sigma \Bs(u) v_{t^{-1}} \\
	& = \lambda \Bs(u_t) (v_{t^{-1}})_t = \lambda \Bs(u_t) v
\end{align*}
for all $v \in W$.

The potentials
\begin{equation} \label{9}
	I_\sigma(u) = \int_0^1 \As(\tau u) u\,\mathrm{d}\tau, \quad J_\sigma(u) = \int_0^1 \Bs(\tau u) u\,\mathrm{d}\tau, \quad u \in W
\end{equation}
of $\As$ and $\Bs$, respectively, are even, bounded on bounded sets, and have the scaling property
\begin{align*}
	I_\sigma(u_t) = t^\sigma I_\sigma(u), \quad J_\sigma(u_t) = t^\sigma J_\sigma(u) \quad \text{for all }u \in W,\, t \ge 0,
\end{align*}
as shown in Mercuri and Perera \cite[Proposition 2.2]{Mercuri-Perera-2026}. By \eqref{9}, \textnormal{(\hyperlink{A6}{A$_6$})}, and \textnormal{(\hyperlink{A8}{A$_8$})},
\begin{equation} \label{16}
	I_\sigma(u) > 0, \quad J_\sigma(u) > 0 \quad \text{for all }u \in W \setminus \set{0}.
\end{equation}
It follows from \textnormal{(\hyperlink{A9}{A$_9$})} that if $u_j \wto u$ in $W$, then $J_\sigma(u_j) \to J_\sigma(u)$ (see Mercuri and Perera \cite[Proposition 2.4]{Mercuri-Perera-2026}). We further assume that $I_\sigma$ and $J_\sigma$ satisfy:
\begin{enumerate}
	\item[\textnormal{(A$_{10}$)}]\hypertarget{A10}{}
		$I_\sigma$ is coercive, i.e., $I_\sigma(u) \to \infty$ as $\norm{u} \to \infty$.
	\item[\textnormal{(A$_{11}$)}]\hypertarget{A11}{}
		Every solution of problem \eqref{8} satisfies the Poho\v{z}aev-type identity $I_\sigma(u) = \lambda J_\sigma(u)$.
\end{enumerate}

Since $\As$ is odd, we have $I_\sigma'(0) = \As(0) = 0$, so the origin is a critical point of $I_\sigma$. It is the only critical point of $I_\sigma$, because
\begin{align*}
	I_\sigma'(u) u = \As(u) u > 0 \quad \text{for all }u \in W \setminus \set{0}
\end{align*}
by \textnormal{(\hyperlink{A6}{A$_6$})}. Thus, $I_\sigma(0) = 0$ is the only critical value of $I_\sigma$, and the implicit function theorem shows that
\begin{align*}
	\M = \{u \in W \colon I_\sigma(u) = 1\}
\end{align*}
is a $C^1$-Finsler manifold. As $I_\sigma$ is continuous, even, and coercive, $\M$ is complete, symmetric, and bounded. Define
\begin{align*}
	\Psi(u) = \frac{1}{J_\sigma(u)}, \quad u \in W \setminus \set{0}
\end{align*}
and let $\widetilde{\Psi} = \restr{\Psi}{\M}$. Then the eigenvalues of problem \eqref{8} coincide with the critical values of $\widetilde{\Psi}$ (Mercuri and Perera \cite[Proposition 2.5]{Mercuri-Perera-2026}).

Let $\F$ denote the class of symmetric subsets of $\M$, and $i(M)$ be the cohomological index of $M \in \F$. For $k \ge 1$, define $\F_k = \{M \in \F \colon i(M) \ge k\}$ and set
\begin{align*}
	\lambda_k := \inf_{M \in \F_k}\, \sup_{u \in M}\, \widetilde{\Psi}(u).
\end{align*}
The following theorem was proved by Mercuri and Perera \cite[Theorem 2.10]{Mercuri-Perera-2026}.

\begin{theorem}\label{Theorem 3}
	Assume \textnormal{(\hyperlink{A1}{A$_1$})}--\textnormal{(\hyperlink{A11}{A$_{11}$})}. Then $\seq{\lambda_k}$ is a sequence of eigenvalues of problem \eqref{8} with $\lambda_k \nearrow \infty$. Moreover, the following statements hold:
	\begin{enumerate}
		\item[\textnormal{(i)}]\hypertarget{Theorem 3.i}{}
			The first eigenvalue is given by
			\begin{align*}
				\lambda_1 = \min_{u \in \M}\, \widetilde{\Psi}(u) > 0.
			\end{align*}
		\item[\textnormal{(ii)}]\hypertarget{Theorem 3.ii}{}
			If $\lambda_k = \ldots = \lambda_{k+m-1} = \lambda$ and $E_\lambda$ denotes the set of eigenfunctions associated with $\lambda$ that belong to $\M$, then $i(E_\lambda) \ge m$.
		\item[\textnormal{(iii)}]\hypertarget{Theorem 3.iii}{}
			If $\lambda_k < \lambda < \lambda_{k+1}$, then
			\begin{align*}
				i(\widetilde{\Psi}^{\lambda_k}) = i(\M \setminus \widetilde{\Psi}_\lambda) = i(\widetilde{\Psi}^\lambda) = i(\M \setminus \widetilde{\Psi}_{\lambda_{k+1}}) = k,
			\end{align*}
			where $\widetilde{\Psi}^a = \{u \in \M \colon \widetilde{\Psi}(u) \le a\}$ and $\widetilde{\Psi}_a = \{u \in \M \colon \widetilde{\Psi}(u) \ge a\}$ for $a \in \R$.
	\end{enumerate}
\end{theorem}

\subsection{Scaling-based multiplicity results for even functionals}

\subsubsection{Unconstrained case}

Let $\M$ be as in the preceding subsection. Making essential use of the piercing property of the cohomological index, the following theorem was proved by Mercuri and Perera \cite[Theorem 2.33]{Mercuri-Perera-2026}.

\begin{theorem}
	Let $E \in C^1(W,\R)$ be an even functional, and assume that there exists $c^\ast > 0$ such that $E$ satisfies the \PS{c} condition for all $c \in (0,c^\ast)$. Let $A_0$ and $B_0$ be symmetric subsets of $\M$ such that $A_0$ is compact, $B_0$ is closed, and
	\begin{align*}
		i(A_0) \ge k + m - 1, \quad i(\M \setminus B_0) \le k - 1
	\end{align*}
	for some $k, m \ge 1$. Assume that there exist $R > \rho > 0$ such that
	\begin{align*}
		\sup_{u \in A}\, E(u) \le 0 < \inf_{u \in B}\, E(u), \quad \sup_{u \in X}\, E(u) < c^\ast,
	\end{align*}
	where $A = \{u_R \colon u \in A_0\}$, $B = \{u_\rho \colon u \in B_0\}$, and $X = \{u_t \colon u \in A_0,\, 0 \le t \le R\}$. Then $E$ has $m$ distinct pairs of critical points at levels in $(0,c^\ast)$.
\end{theorem}

Taking $B_0=\M$ and $k=1$, we obtain the following corollary. It is stated in Mercuri and Perera \cite[Corollary 2.34]{Mercuri-Perera-2026} under the condition $i(A_0)=m$, but the argument only requires $i(A_0)\ge m$. The special case $u_t=tu$ was proved by Perera \cite{Perera-2025}.

\begin{corollary}\label{Corollary 1}
	Let $E \in C^1(W,\R)$ be an even functional, and assume that there exists $c^\ast > 0$ such that $E$ satisfies the \PS{c} condition for all $c \in (0,c^\ast)$. Let $A_0$ be a compact symmetric subset of $\M$ with $i(A_0) \ge m$ for some $m \ge 1$. Assume that there exist $R > \rho > 0$ such that
	\begin{align*}
		\sup_{u \in A}\, E(u) \le 0 < \inf_{u \in \M_\rho}\, E(u), \quad \sup_{u \in X}\, E(u) < c^\ast,
	\end{align*}
	where $A = \{u_R \colon u \in A_0\}$, $\M_\rho = \{u_\rho \colon u \in \M\}$, and $X = \{u_t \colon u \in A_0,\, 0 \le t \le R\}$. Then $E$ has $m$ distinct pairs of critical points at levels in $(0,c^\ast)$.
\end{corollary}

In particular, we obtain the following consequence when the \PS{c} condition holds at every positive level $c\in (0,\infty)$.

\begin{corollary} 
	Let $E \in C^1(W,\R)$ be an even functional that satisfies the \PS{c} condition for all $c \in (0,\infty)$. Assume that, for any $m \ge 1$, there exist a compact symmetric subset $A_0$ of $\M$ with $i(A_0) \ge m$ and $R > \rho > 0$ such that
	\begin{align*}
		\sup_{u \in A}\, E(u) \le 0 < \inf_{u \in \M_\rho}\, E(u),
	\end{align*}
	where $A = \{u_R \colon u \in A_0\}$ and $\M_\rho = \{u_\rho \colon u \in \M\}$. Then $E$ has infinitely many critical points at positive levels.
\end{corollary}

\subsubsection{Constrained case}

Let $\As$ be as defined in Subsection \ref{Scaled eigenvalue problems}. We consider the nonlinear operator equation
\begin{equation} \label{18}
	\As(u) = f(u) \quad \text{in } W^\ast,
\end{equation}
where $f \in C(W,W^\ast)$ is an odd potential operator that maps bounded sets into bounded sets and satisfies:
\begin{enumerate}
	\item[\textnormal{(f$_1$)}]\hypertarget{f1}{}
		There exists $\rho \ne \sigma$ such that $f(u_t) v_t = t^\rho f(u) v$ for all $u, v \in W$ and $t \ge 0$.
	\item[\textnormal{(f$_2$)}]\hypertarget{f2}{}
		$f(u) u > 0$ for all $u \in W \setminus \set{0}$.
	\item[\textnormal{(f$_3$)}]\hypertarget{f3}{}
		If $u_j \wto u$ in $W$, then $f(u_j) \to f(u)$ in $W^\ast$.
	\item[\textnormal{(f$_4$)}]\hypertarget{f4}{}
		For any $\mu \in \R$, every solution of the equation $\As(u) = \mu f(u)$ satisfies the Poho\v{z}aev-type identity $\sigma I_\sigma(u) = \rho \mu F(u)$, where $F(u) = \int_0^1 f(\tau u) u\,\mathrm{d}\tau$ is the potential of $f$.
\end{enumerate}
Solutions of equation \eqref{18} are the critical points of the $C^1$-functional
\begin{align*}
	E(u) = I_\sigma(u) - F(u), \quad u \in W,
\end{align*}
where $I_\sigma$ is defined in \eqref{9}. However, verifying the boundedness of \PS{} sequences for $E$ appears difficult.

Let $\M$ be as in Subsection \ref{Scaled eigenvalue problems}. As noted in Mercuri and Perera \cite[Section 2.2.1]{Mercuri-Perera-2026}, the scaling induces a continuous projection
\begin{equation} \label{projection}
	\pi \colon W \setminus \set{0} \to \M, \quad \pi(u) = \widetilde{u} := u_{t_u},
\end{equation}
where $t_u = I_\sigma(u)^{- 1/\sigma}$, see \eqref{16}. Every $u \in W \setminus \set{0}$ lies on the fiber $\{\widetilde{u}_t \colon t > 0\}$ through $\pi(u)$. More precisely, $u = \widetilde{u}_{I_\sigma(u)^{1/\sigma}}$, and $\pi(u_t) = \pi(u)$ for all $t > 0$. This projection allows us to reduce equation \eqref{18} to the nonlinear eigenvalue problem
\begin{align*} 
	\As(u) = \mu f(u)
\end{align*}
on $\M$. Using this constrained approach, Perera \cite[Theorem 3.8]{Perera-2026} proved the following theorem.

\begin{theorem}\label{Theorem 4}
	Assume \textnormal{(\hyperlink{A1}{A$_1$})}--\textnormal{(\hyperlink{A7}{A$_7$})}, \textnormal{(\hyperlink{A10}{A$_{10}$})}, and \textnormal{(\hyperlink{f1}{f$_1$})}--\textnormal{(\hyperlink{f4}{f$_4$})}. Then equation \eqref{18} admits an infinite sequence of nontrivial solutions
	\begin{align*} 
		\overline{u}_k = u^k_{\mu_k^{- 1/(\sigma - \rho)}},
	\end{align*}
	where $u^k \in \M$, $\mu_k > 0$, and $\mu_k \nearrow \infty$. Moreover,
	\begin{align*}
		E(\overline{u}_k) = \left(1 - \frac{\sigma}{\rho}\right) \mu_k^{- \sigma/(\sigma - \rho)}.
	\end{align*}
	In particular, $E(\overline{u}_k) < 0$ and $E(\overline{u}_k) \nearrow 0$ if $\rho < \sigma$, whereas $E(\overline{u}_k) > 0$ and $E(\overline{u}_k) \nearrow \infty$ if $\rho > \sigma$.
\end{theorem}

\subsubsection{Negative critical levels}

The following standard result in critical point theory provides a criterion for the existence of infinitely many critical points at negative levels (Perera, Agarwal, and O'Regan \cite[Proposition 3.36]{Perera-Agarwal-ORegan-2010}).

\begin{proposition} \label{Proposition negative levels}
	Let $\Phi$ be an even $C^1$-functional on a Banach space $W$ such that $\Phi(0) = 0$ and $\Phi$ satisfies the \PS{c} condition for all $c < 0$. Let $\G$ denote the class of symmetric subsets of $W \setminus \set{0}$. For $k \ge 1$, define $\G_k = \{M \in \G \colon i(M) \ge k\}$ and set
	\begin{align*}
		c_k := \inf_{M \in \G_k}\, \sup_{u \in M}\, \Phi(u).
	\end{align*}
	If there exists a $k_0 \ge 1$ such that $- \infty < c_k < 0$ for all $k \ge k_0$, then $c_{k_0} \le c_{k_0 + 1} \le \dotsb \to 0$ is a sequence of critical values of $\Phi$.
\end{proposition}

\section{Existence and multiplicity results} \label{Section results}

In this section, we establish the general results underlying Theorems \ref{thm:intro-eigenvalues}, \ref{thm:intro-pure-power}, \ref{thm:intro-nonexistence}, \ref{thm:intro-large-lambda}, and \ref{thm:intro-concave-convex}. We work in the space $W_\rad$ with the scaling \eqref{scaling}, the scaled operator $\As$ in \eqref{As}, and its potential $I_\sigma$ in \eqref{I_sigma}. Recall that
\begin{align*}
	\M = \{u \in W_\rad \colon I_\sigma(u) = 1\}
\end{align*}
and that $\pi \colon W_\rad \setminus \{0\} \to \M$ is the projection along fibers defined in \eqref{projection}. For $(s,c) \in \Te$, we write
\begin{align*}
	J(u) = \int_{\R^N} |x|^{-c}\, |u|^s\, \mathrm{d}x, \quad u \in W_\rad
\end{align*}
and recall that $J(u_t) = t^{\rho}\, J(u)$, where $\rho = \rho(s,c)$ is given by \eqref{rho sigma}.

\subsection{Scaled eigenvalue problem}

Let $(s,c)$ lie on the open scaling line $(Q_1,Q_2)$, i.e., $p < s < q$ and $c = (\alpha + 1)\, s - \alpha p$. Then  $\rho(s,c) = \sigma$, and the equation
\begin{equation} \label{7}
	- \operatorname{div} \left(|\nabla u|^{p-2}\, \nabla u + |x|^{-b}\, |\nabla u|^{q-2}\, \nabla u\right) = \lambda\, |x|^{-c}\, |u|^{s-2}\, u \quad \text{in } \R^N
\end{equation}
is a scaled eigenvalue problem in the sense of Subsection \ref{Scaled eigenvalue problems}. Thus, if $u$ is a solution, every function in the one-parameter family $\{u_t \colon t > 0\}$ is also a solution. The scaled operator that corresponds to the right-hand side of equation \eqref{7} and its potential are given by
\begin{align*}
	\Bs(u) v = \int_{\R^N} |x|^{-c}\, |u|^{s-2}\, uv\, \mathrm{d}x, \quad u, v \in W_\rad
\end{align*}
and
\begin{align*}
	J_\sigma(u) = \frac{1}{s} \int_{\R^N} |x|^{-c}\, |u|^s\, \mathrm{d}x, \quad u \in W_\rad,
\end{align*}
respectively. Let $\widetilde{\Psi}$ be as in Subsection \ref{Scaled eigenvalue problems}.

\begin{theorem} \label{Theorem 1}
	If $p < s < q$ and $c = (\alpha + 1)\, s - \alpha p$, then the scaled eigenvalue problem \eqref{7} has a sequence of positive eigenvalues $\lambda_k(s) \nearrow \infty$ with radial eigenfunctions, and the conclusions \textnormal{(\hyperlink{Theorem 3.i}{i})}--\textnormal{(\hyperlink{Theorem 3.iii}{iii})} of Theorem \ref{Theorem 3} hold.
\end{theorem}

\begin{proof}
	Assumption \textnormal{(\hyperlink{A8}{A$_8$})} is clearly satisfied, while \textnormal{(\hyperlink{A9}{A$_9$})} follows from the compactness of the embedding $W_\rad \hookrightarrow L^s(\R^N;|x|^{-c})$ since the open segment $(Q_1,Q_2)$ is contained in $\Tc$ (see Theorem \ref{Theorem embedding} \textnormal{(\hyperlink{Theorem embedding.ii}{ii})}). Assumption \textnormal{(\hyperlink{A11}{A$_{11}$})} is the special case $n = 1$ and $\rho_1 = \sigma$ of Theorem \ref{thm:Pohozaev}. So the result follows from Theorem \ref{Theorem 3}.
\end{proof}

\subsection{Pure power equation}

Next, we consider the pure power equation \eqref{7} for $(s,c)$ off the scaling line. Its associated energy functional is
\begin{align*}
	E(u) = \frac{1}{p} \int_{\R^N} |\nabla u|^p\, \mathrm{d}x + \frac{1}{q} \int_{\R^N} |x|^{-b}\, |\nabla u|^q\, \mathrm{d}x - \frac{\lambda}{s} \int_{\R^N} |x|^{-c}\, |u|^s\, \mathrm{d}x,
\end{align*}
for $u \in W_\rad$.

\begin{theorem} \label{Theorem 2}
	Let $\lambda > 0$ and $(s,c) \in \Tc$ with $s > 1$ and $\rho(s,c) \neq \sigma$. Then equation \eqref{7} has an infinite sequence of nontrivial weak radial solutions $\seq{u_k}$.
	\begin{enumerate}
		\item[\textnormal{(i)}]\hypertarget{Theorem 2.i}{}
			If $(s,c)$ lies above the scaling line, i.e., $c > (\alpha + 1)\, s - \alpha p$, then equation \eqref{7} is subscaled, $E(u_k) < 0$, and $E(u_k) \nearrow 0$.
		\item[\textnormal{(ii)}]\hypertarget{Theorem 2.ii}{}
			If $(s,c)$ lies below the scaling line, i.e., $c < (\alpha + 1)\, s - \alpha p$, then equation \eqref{7} is superscaled, $E(u_k) > 0$, and $E(u_k) \nearrow \infty$.
	\end{enumerate}
\end{theorem}

\begin{proof}
	We apply Theorem \ref{Theorem 4} with
	\begin{align*}
		f(u) v &= \lambda \int_{\R^N} |x|^{-c}\, |u|^{s-2}\, u\,v\, \mathrm{d}x,\\
		F(u) &= \frac{\lambda}{s} \int_{\R^N} |x|^{-c}\, |u|^s\, \mathrm{d}x,
	\end{align*}
	for $u, v \in W_\rad$. Assumptions \textnormal{(\hyperlink{f1}{f$_1$})} and \textnormal{(\hyperlink{f2}{f$_2$})} are clearly satisfied, \textnormal{(\hyperlink{f3}{f$_3$})} follows from the compactness of the embedding $W_\rad \hookrightarrow L^s(\R^N;|x|^{-c})$ (see Theorem \ref{Theorem embedding} \textnormal{(\hyperlink{Theorem embedding.ii}{ii})}), and \textnormal{(\hyperlink{f4}{f$_4$})} is the special case $n = 1$ and $\lambda_1 = \lambda \mu$ of Theorem \ref{thm:Pohozaev}. The conclusions now follow from Theorem \ref{Theorem 4} since $\rho < \sigma$ in case \textnormal{(\hyperlink{Theorem 2.i}{i})} and $\rho > \sigma$ in case \textnormal{(\hyperlink{Theorem 2.ii}{ii})}.
\end{proof}

In particular, we have the following corollary for the unweighted pure power equation
\begin{equation} \label{27}
	- \operatorname{div} \left(|\nabla u|^{p-2}\, \nabla u + |x|^{-b}\, |\nabla u|^{q-2}\, \nabla u\right) = \lambda\, |u|^{s-2}\, u \quad \text{in } \R^N.
\end{equation}

\begin{corollary} \label{cor:pure-power}
	If $\lambda > 0$ and $p^\ast < s < q^\ast_b$, then equation \eqref{27} has an infinite sequence of nontrivial weak radial solutions $\seq{u_k}$ such that $E(u_k) > 0$ and $E(u_k) \nearrow \infty$.
\end{corollary}

\begin{proof}
	The point $(s,0)$ lies on the open base of $\T$ and hence belongs to $\Tc$ below the scaling line.
\end{proof}

\subsection{Critical equations} 

In this subsection, we consider the critical equation
\begin{equation} \label{19}
	\begin{aligned}
		& - \operatorname{div} \left(|\nabla u|^{p-2}\, \nabla u + |x|^{-b}\, |\nabla u|^{q-2}\, \nabla u\right) \\
		& = \lambda\, |x|^{-c}\, |u|^{s-2}\, u + \mu_1\, |x|^{- c_1}\, |u|^{s_1 - 2}\, u + \mu_2\, |x|^{- c_2}\, |u|^{s_2 - 2}\, u \quad \text{in } \R^N,
	\end{aligned}
\end{equation}
where the following assumptions hold:
\begin{enumerate}
	\item[\textnormal{(H$_1$)}]\hypertarget{H1}{}
		$s > 1$ and $(s,c) \in \Tc$.
	\item[\textnormal{(H$_2$)}]\hypertarget{H2}{}
		$(s_1,c_1) \in [Q_1,P_1]$ and $(s_2,c_2) \in [Q_2,P_3]$.
	\item[\textnormal{(H$_3$)}]\hypertarget{H3}{}
		$\mu_1, \mu_2 \ge 0$ with $\mu_1 + \mu_2 > 0$.
\end{enumerate}
The energy functional associated with this equation is
\begin{align*}
	E(u) = I_\sigma(u) - \frac{\lambda}{s} J(u) - \frac{\mu_1}{s_1}\, K_1(u) - \frac{\mu_2}{s_2}\, K_2(u), \quad u \in W_\rad,
\end{align*}
where
\begin{align*}
	K_i(u) = \int_{\R^N} |x|^{-c_i}\, |u|^{s_i}\, \mathrm{d}x, \quad i = 1, 2.
\end{align*}
We write $\rho = \rho(s,c)$ and $\rho_i = \rho(s_i,c_i)$, and recall that
\begin{equation} \label{rho1-rho2}
	\rho_1 = \frac{\sigma s_1}{p}, \quad \rho_2 = \frac{\sigma s_2}{q},
\end{equation}
by \eqref{crit-segments}. We also set
\begin{align*}
	\kappa = \frac{N}{p} - \frac{N - b}{q} > 0
\end{align*}
and write
\begin{align*}
	X(u) = \int_{\R^N} |\nabla u|^p\, \mathrm{d}x, \quad Y(u) = \int_{\R^N} |x|^{-b}\, |\nabla u|^q\, \mathrm{d}x.
\end{align*}

\subsubsection{Nonexistence results}

\begin{lemma} \label{Lemma Pohozaev critical}
	Assume \textnormal{(\hyperlink{H1}{H$_1$})}--\textnormal{(\hyperlink{H3}{H$_3$})}. Then every weak radial solution $u \in W_\rad$ of equation \eqref{19} satisfies
	\begin{align}
		\kappa\, \big[X(u) - \mu_1 K_1(u)\big] & = \frac{\lambda}{s} \left(N - c - \frac{Ns}{q^\ast_b}\right) J(u), \label{Pohozaev X} \\
		\kappa\, \big[Y(u) - \mu_2\,  K_2(u)\big] & = \frac{\lambda}{s} \left(c - N + \frac{Ns}{p^\ast}\right) J(u), \label{Pohozaev Y}
	\end{align}
	where both brackets on the right are positive.
\end{lemma}

\begin{proof}
	The expressions in parentheses are positive because $(s,c) \in \Tc$. By the dilation identity \eqref{dilation},
	\begin{align*}
		& \frac{N - p}{p}\, X(u) + \frac{N - q - b}{q}\, Y(u) \\
		& = \frac{\lambda\, (N - c)}{s}\, J(u) + \frac{N - p}{p}\, \mu_1 K_1(u) + \frac{N - q - b}{q}\, \mu_2\,  K_2(u)
	\end{align*}
	since
	\begin{align*}
		\frac{N-c_1}{s_1}=\frac N{p^\ast}=\frac{N-p}{p}, \quad \frac{N-c_2}{s_2}=\frac N{q^\ast_b} =\frac{N-q-b}{q}.
	\end{align*}
	Subtracting $(N - q - b)/q$ (resp.\,$(N - p)/p$) times the Nehari identity
	\begin{align*}
		X(u) + Y(u) = \lambda J(u) + \mu_1 K_1(u) + \mu_2\,  K_2(u)
	\end{align*}
	(see \eqref{Nehari}) gives \eqref{Pohozaev X} (resp.\,\eqref{Pohozaev Y}).
\end{proof}

The preceding identities immediately imply a nonexistence result when at least one critical term is absent and the perturbation parameter is nonpositive.

\begin{theorem} \label{Theorem nonexistence}
	Assume \textnormal{(\hyperlink{H1}{H$_1$})}--\textnormal{(\hyperlink{H3}{H$_3$})} and let $\lambda \le 0$. If $\mu_1 = 0$ or $\mu_2 = 0$, then equation \eqref{19} has no nontrivial weak radial solution $u \in W_\rad$. In particular, the pure critical power equations
	\begin{align*}
		- \operatorname{div} \left(|\nabla u|^{p-2}\, \nabla u + |x|^{-b}\, |\nabla u|^{q-2}\, \nabla u\right) = |x|^{- c_i}\, |u|^{s_i - 2}\, u \quad \text{in } \R^N, \quad i = 1, 2
	\end{align*}
	have no nontrivial weak radial solutions.
\end{theorem}

\begin{proof}
	If $\mu_1 = 0$, then \eqref{Pohozaev X} gives $\kappa X(u) \le 0$, so $X(u) = 0$ and hence $u = 0$. If $\mu_2 = 0$, \eqref{Pohozaev Y} similarly implies $\kappa Y(u) \le 0$, and we again conclude that $u = 0$.
\end{proof}

In particular, we obtain the following corollary for the unweighted pure critical power equations
\begin{equation} \label{p-ast}
	- \operatorname{div} \left(|\nabla u|^{p-2}\, \nabla u + |x|^{-b}\, |\nabla u|^{q-2}\, \nabla u\right) = |u|^{p^\ast - 2}\, u \quad \text{in } \R^N
\end{equation}
and
\begin{equation} \label{q-ast-b}
	- \operatorname{div} \left(|\nabla u|^{p-2}\, \nabla u + |x|^{-b}\, |\nabla u|^{q-2}\, \nabla u\right) = |u|^{q^\ast_b - 2}\, u \quad \text{in } \R^N.
\end{equation}

\begin{corollary} \label{cor:nonexistence}
	Equations \eqref{p-ast} and \eqref{q-ast-b} have no nontrivial weak radial solutions.
\end{corollary}

When $\lambda \le 0$ and both critical terms are present, Lemma \ref{Lemma Pohozaev critical} only gives $X(u) \le \mu_1 K_1(u)$ and $Y(u) \le \mu_2\,  K_2(u)$. The case $\lambda = 0$ is discussed further in Remark \ref{Remark lambda zero}.

\subsubsection{Threshold for compactness} 

To obtain a compactness threshold, we impose the additional assumption:
\begin{enumerate}
	\item[(H$_4$)]\hypertarget{H4}{}
		$s_1 > q$ if $\mu_1 > 0$, and $s_2 > q$ if $\mu_2 > 0$.
\end{enumerate}
In particular, $q < p^\ast$ if $\mu_1 > 0$. Let
\begin{align*}
	S_1 = \inf_{u \in W_\rad \setminus \set{0}}\, \frac{X(u)}{K_1(u)^{p/s_1}}, \quad S_2 = \inf_{u \in W_\rad \setminus \set{0}}\, \frac{Y(u)}{K_2(u)^{q/s_2}},
\end{align*}
which are positive by Lemma \ref{Lemma critical segments}. Set
\begin{align*}
	\beta_1 &= \left(\frac{1}{p} - \frac{1}{s_1}\right) \mu_1^{- p/(s_1 - p)}\, S_1^{s_1/(s_1 - p)},\\
	\beta_2 &= \left(\frac{1}{q} - \frac{1}{s_2}\right) \mu_2^{- q/(s_2 - q)}\, S_2^{s_2/(s_2 - q)},
\end{align*}
with the convention that $\beta_i = \infty$ if $\mu_i = 0$, and define
\begin{align*}
	\beta^\ast = \min \set{\beta_1,\beta_2}.
\end{align*}

We first establish an algebraic estimate that will determine the threshold for loss of compactness.

\begin{lemma} \label{Lemma threshold}
	Assume \textnormal{(\hyperlink{H2}{H$_2$})}--\textnormal{(\hyperlink{H4}{H$_4$})}. Let $X, Y, Z, W \ge 0$ satisfy $(X,Y) \ne (0,0)$ and
	\begin{equation} \label{XYZW}
		X + Y = Z + W, \quad Z \le \mu_1\, S_1^{- s_1/p}\, X^{s_1/p}, \quad W \le \mu_2\, S_2^{- s_2/q}\, Y^{s_2/q}.
	\end{equation}
	Then
	\begin{equation} \label{ineq}
		\frac{1}{p}\, X + \frac{1}{q}\, Y - \frac{1}{s_1}\, Z - \frac{1}{s_2}\, W \ge \beta^\ast.
	\end{equation}
	This inequality is sharp, with equality for $(X,Y,Z,W) = (X_1,0,X_1,0)$ if $\beta^\ast = \beta_1$, and for $(0,Y_2,0,Y_2)$ if $\beta^\ast = \beta_2$, where, for $\mu_1>0$ and $\mu_2>0$, respectively, set
	\begin{equation} \label{X1-Y2}
		X_1 = \mu_1^{- p/(s_1 - p)}\, S_1^{s_1/(s_1 - p)}, \quad Y_2 = \mu_2^{- q/(s_2 - q)}\, S_2^{s_2/(s_2 - q)}.
	\end{equation}
\end{lemma}

\begin{proof}
	Denote the left-hand side of \eqref{ineq} by $I$. Since $X + Y = Z + W$, either $Z \ge X$ or $W \ge Y$.

	\emph{Case 1: $Z \ge X$.} If $X = 0$, then $Z = 0$ by \eqref{XYZW}, so $W = Y$ and we are in Case 2. Thus, assume $X > 0$. Then $Z > 0$, so $\mu_1 > 0$, and
	\begin{align*}
		X \le Z \le \mu_1 S_1^{- s_1/p}\, X^{s_1/p}
	\end{align*}
	implies $X \ge X_1$. Moreover, $W = X + Y - Z \le Y$ and
	\begin{align*}
		I = \left(\frac{1}{p} - \frac{1}{s_1}\right) X + \left(\frac{1}{q} - \frac{1}{s_1}\right) Y + \left(\frac{1}{s_1} - \frac{1}{s_2}\right) W.
	\end{align*}
	If $s_1 \le s_2$, the last two terms are nonnegative because $s_1 > q$ by \textnormal{(\hyperlink{H4}{H$_4$})}. If $s_1 > s_2$, then
	\begin{align*}
		\left(\frac1{s_1}-\frac1{s_2}\right)W \ge \left(\frac1{s_1}-\frac1{s_2}\right)Y,
	\end{align*}
	so the sum of the last two terms is at least $(1/q - 1/s_2)\, Y \ge 0$, since $s_2 \ge q$.
	In either case,
	\begin{align*}
		I\ge \left(\frac1p-\frac1{s_1}\right)X \ge \left(\frac1p-\frac1{s_1}\right)X_1 =\beta_1.
	\end{align*}

	\emph{Case 2: $W \ge Y$.} If $Y = 0$, then $W = 0$ by \eqref{XYZW}, so $Z = X > 0$ and we are in Case 1. We may therefore assume $Y > 0$. Then $W > 0$, so $\mu_2 > 0$, and $Y \le W \le \mu_2 \,S_2^{- s_2/q}\, Y^{s_2/q}$ gives $Y \ge Y_2$. Furthermore, $Z = X + Y - W \le X$ and
	\begin{align*}
		I = \left(\frac{1}{p} - \frac{1}{s_2}\right) X + \left(\frac{1}{q} - \frac{1}{s_2}\right) Y + \left(\frac{1}{s_2} - \frac{1}{s_1}\right) Z.
	\end{align*}
	If $s_2 \le s_1$, the first and the third terms are nonnegative. If $s_2 > s_1$, then
	\begin{align*}
		\left(\frac1{s_2}-\frac1{s_1}\right)Z \ge \left(\frac1{s_2}-\frac1{s_1}\right)X,
	\end{align*}
	so the sum of the first and the third terms is at least $(1/p - 1/s_1)\, X \ge 0$.
	Consequently,
	\begin{align*}
		I\ge \left(\frac1q-\frac1{s_2}\right)Y \ge \left(\frac1q-\frac1{s_2}\right)Y_2 =\beta_2.
	\end{align*}
	The equality cases follow by direct computation.
\end{proof}

We are now ready to prove the compactness result.

\begin{proposition} \label{Proposition PS}
	Assume \textnormal{(\hyperlink{H1}{H$_1$})}--\textnormal{(\hyperlink{H4}{H$_4$})} and let $\lambda \in \R$. If $\seq{u_j}$ is a bounded \PS{c} sequence of $E$, then, after passing to a subsequence, $u_j$ converges weakly to a weak radial solution $u$ of equation \eqref{19}, and either $u_j \to u$ in $W_\rad$ or
	\begin{align*}
		c \ge E(u) + \beta^\ast.
	\end{align*}
	In particular, if $\lambda\, (\rho - \sigma) \ge 0$, then every bounded \PS{c} sequence of $E$ has a strongly convergent subsequence whenever $c<\beta^\ast$.
\end{proposition}

\begin{proof}
	By Lemma \ref{Lemma weak limit}, a subsequence of $\seq{u_j}$, still denoted by $\seq{u_j}$, converges weakly to a weak radial solution $u$ of equation \eqref{19}, and $u_j \to u$ and $\nabla u_j \to \nabla u$ a.e.\,in $\R^N$. Set $v_j = u_j - u$. By the Brezis-Lieb lemma \cite{Brezis-Lieb-1983},
	\begin{align*}
		X(u_j) &= X(u) + X(v_j) + o(1),\\
		Y(u_j) &= Y(u) + Y(v_j) + o(1),\\
		K_i(u_j) &= K_i(u) + K_i(v_j) + o(1), \quad i = 1, 2,
	\end{align*}
	and $J(u_j) \to J(u)$ by Theorem \ref{Theorem embedding} \textnormal{(\hyperlink{Theorem embedding.ii}{ii})}. Since $E'(u_j)\, u_j \to 0$ by the boundedness of $\seq{u_j}$, and $E'(u)\, u = 0$, we arrive at
	\begin{align}
		\frac{1}{p}\, X(v_j) + \frac{1}{q}\, Y(v_j) - \frac{\mu_1}{s_1}\, K_1(v_j) - \frac{\mu_2}{s_2}\, K_2(v_j) &= c - E(u) + o(1), \label{E(v_j)}\\
		X(v_j) + Y(v_j) - \mu_1\, K_1(v_j) - \mu_2\, K_2(v_j) &= o(1). \label{E'(v_j)}
	\end{align}
	Suppose that $u_j \nrightarrow u$. Then, for a renamed subsequence, $\norm{v_j} \ge \delta$ for some $\delta > 0$ and
	\begin{align*}
		X(v_j) \to X, \quad Y(v_j) \to Y, \quad \mu_1 K_1(v_j) \to Z, \qquad \mu_2\,  K_2(v_j) \to W
	\end{align*}
	for some $X, Y, Z, W \ge 0$ with $(X,Y) \ne (0,0)$. Passing to the limit in equations \eqref{E(v_j)} and \eqref{E'(v_j)}, and the inequalities $K_1(v_j) \le S_1^{- s_1/p}\, X(v_j)^{s_1/p}$ and $K_2(v_j) \le S_2^{- s_2/q}\, Y(v_j)^{s_2/q}$, we see that $(X,Y,Z,W)$ satisfies \eqref{XYZW} and
	\begin{align*}
		\frac{1}{p}\, X + \frac{1}{q}\, Y - \frac{1}{s_1}\, Z - \frac{1}{s_2}\, W = c - E(u).
	\end{align*}
	Lemma \ref{Lemma threshold} therefore implies $c - E(u) \ge \beta^\ast$. Finally, by the energy identity \eqref{energy identity},
	\begin{equation} \label{energy critical}
		E(u) = \left(\frac{\rho}{\sigma} - 1\right) \frac{\lambda}{s}\, J(u) + \left(\frac{\rho_1}{\sigma} - 1\right) \frac{\mu_1}{s_1}\, K_1(u) + \left(\frac{\rho_2}{\sigma} - 1\right) \frac{\mu_2}{s_2}\, K_2(u),
	\end{equation}
	where $\rho_1 = \sigma s_1/p \ge \sigma$ and $\rho_2 = \sigma s_2/q \ge \sigma$ by \eqref{rho1-rho2}. Thus, $E(u) \ge 0$ whenever $\lambda\, (\rho - \sigma) \ge 0$, and the last assertion follows.
\end{proof}

\begin{remark} \label{Remark lambda zero}
	When $\lambda = 0$ and $\mu_1, \mu_2 > 0$, Lemma \ref{Lemma Pohozaev critical} gives $X(u) = \mu_1 K_1(u)$ and $Y(u) = \mu_2\,  K_2(u)$ for every weak radial solution $u$. If $u \ne 0$, then the definitions of $S_1$ and $S_2$ imply $X(u) \ge X_1$ and $Y(u) \ge Y_2$ (see \eqref{X1-Y2}). Hence,
	\begin{align*}
		E(u) = \left(\frac{1}{p} - \frac{1}{s_1}\right) X(u) + \left(\frac{1}{q} - \frac{1}{s_2}\right) Y(u) \ge \beta_1 + \beta_2 > \beta^\ast.
	\end{align*}
	Combined with Proposition \ref{Proposition PS}, this shows that when $\lambda = 0$, $E$ has no bounded \PS{c} sequences with $0<c < \beta^\ast$. Whether equation \eqref{19} admits a nontrivial weak radial solution in this case remains an open problem.
\end{remark}

We next show that the threshold $\beta^\ast$ in Proposition \ref{Proposition PS} is optimal when the critical terms are unweighted, i.e., when $(s_1,c_1) = P_1$ and $(s_2,c_2) = P_3$. The two parts of the next proposition realize the two pure profiles in Lemma \ref{Lemma threshold} as Palais-Smale sequences. They correspond to the two mechanisms of noncompactness in Theorem \ref{Theorem embedding} \textnormal{(\hyperlink{Theorem embedding.iii}{iii})}.

\begin{proposition} \label{Proposition sharp}
	Assume \textnormal{(\hyperlink{H1}{H$_1$})}--\textnormal{(\hyperlink{H4}{H$_4$})} and let $\lambda \in \R$.
	\begin{enumerate}
		\item[\textnormal{(i)}]\hypertarget{Proposition sharp.i}{}
			If $\mu_1 > 0$ and $(s_1,c_1) = P_1$, then $E$ has a bounded \PS{\beta_1} sequence that converges weakly to $0$, has no convergent subsequence, and spreads to infinity.
		\item[\textnormal{(ii)}]\hypertarget{Proposition sharp.ii}{}
			If $\mu_2 > 0$ and $(s_2,c_2) = P_3$, then $E$ has a bounded \PS{\beta_2} sequence that converges weakly to $0$, has no convergent subsequence, and concentrates at the origin.
	\end{enumerate}
	In particular, if $(s_i,c_i)$ is the endpoint $P_1$, respectively $P_3$, whenever $\mu_i > 0$, then $E$ has a bounded \PS{\beta^\ast} sequence without convergent subsequences. If, in addition, $\lambda\, (\rho - \sigma) \ge 0$, then $\beta^\ast$ is the largest level below which all bounded Palais-Smale sequences of $E$ are relatively compact.
\end{proposition}

\begin{proof}
	For $v \in W_\rad$, H\"{o}lder's inequality and Lemma \ref{Lemma critical segments} give
	\begin{equation} \label{dual estimates}
		\begin{aligned}
			\left|\int_{\R^N} |\nabla u|^{p-2}\, \nabla u \cdot \nabla v\, \mathrm{d}x\right| & \le X(u)^{(p-1)/p} \norm{v},\\
			\left|\int_{\R^N} |x|^{-b}\, |\nabla u|^{q-2}\, \nabla u \cdot \nabla v\, \mathrm{d}x\right| & \le Y(u)^{(q-1)/q} \norm{v},\\
			\left|\int_{\R^N} |x|^{-c_1}\, |u|^{s_1 - 2}\, u\, v\, \mathrm{d}x\right| & \le C\, X(u)^{(s_1 - 1)/p} \norm{v},\\
			\left|\int_{\R^N} |x|^{-c_2}\, |u|^{s_2 - 2}\, u\, v\, \mathrm{d}x\right| & \le C\, Y(u)^{(s_2 - 1)/q} \norm{v}
		\end{aligned}
	\end{equation}
	for all $u \in W_\rad$, and
	\begin{align*}
		\left|\int_{\R^N} |x|^{-c}\, |u|^{s-2}\, u\, v\, \mathrm{d}x\right| \le C \norm[L^s(\R^N;|x|^{-c})]{u}^{s-1} \norm{v}
	\end{align*}
	by Theorem \ref{Theorem embedding} \textnormal{(\hyperlink{Theorem embedding.i}{i})}. If $w_j \wto 0$ in $W_\rad$, then $w_j \to 0$ in $L^s(\R^N;|x|^{-c})$ by \textnormal{(\hyperlink{H1}{H$_1$})} and Theorem \ref{Theorem embedding} \textnormal{(\hyperlink{Theorem embedding.ii}{ii})}, so the $\lambda$-term in $E$ and in $E'$ is negligible along such sequences.

	\textnormal{(\hyperlink{Proposition sharp.i}{i})} Here $s_1 = p^\ast$ and $c_1 = 0$, and $s_1 > q$ by \textnormal{(\hyperlink{H4}{H$_4$})}. Let
	\begin{align*}
		U(x) = \left(1 + |x|^{p/(p-1)}\right)^{-(N-p)/p}
	\end{align*}
	be the Aubin-Talenti function, which attains the Sobolev constant $S$ of $\D^{1,p}(\R^N)$; see Aubin \cite{Aubin-1976} and Talenti \cite{Talenti-1976}. We have
	\begin{align*}
		|U(x)| \le C\, (1 + |x|)^{-(N-p)/(p-1)}, \quad |\nabla U(x)| \le C\, (1 + |x|)^{-(N-1)/(p-1)}.
	\end{align*}
	Since $b + q < N$ and
	\begin{align*}
		b + q\, \frac{N - 1}{p - 1} > p\, \frac{N - 1}{p - 1} = N + \frac{N - p}{p - 1} > N,
	\end{align*}
	$U$ satisfies \eqref{eq:density}, and hence $U \in W_\rad$ by Lemma \ref{Lemma density}. Since $W_\rad \subset \D^{1,p}(\R^N)$, it follows that $S_1 = S$ and that $S_1$ is attained at $U$. Let $w = tU$, where $t > 0$ is chosen so that $X(w) = \mu_1 K_1(w)$. Since $w$ minimizes $X(u)/K_1(u)^{p/p^\ast}$ over $\D^{1,p}(\R^N) \setminus \set{0}$, the Lagrange multiplier rule gives
	\begin{equation} \label{Talenti equation}
		\int_{\R^N} |\nabla w|^{p-2}\, \nabla w \cdot \nabla v\, \mathrm{d}x = \mu_1 \int_{\R^N} |w|^{p^\ast - 2}\, w\, v\, \mathrm{d}x \quad \text{for all } v \in \D^{1,p}(\R^N).
	\end{equation}
	Moreover, $X(w) = S_1\, K_1(w)^{p/p^\ast} = S_1\, (X(w)/\mu_1)^{p/p^\ast}$, so $X(w) = X_1$ (see \eqref{X1-Y2}).

	Let $\eps_j \searrow 0$ and $w_j(x) = \eps_j^{(N-p)/p}\, w(\eps_j x)$. Then $X(w_j) = X_1$ and $K_1(w_j) = K_1(w)$. Moreover, $w_j$ satisfies \eqref{Talenti equation}, and
	\begin{align*}
		Y(w_j) = \eps_j^{\kappa q}\, Y(w) \to 0.
	\end{align*}
	So $\seq{w_j}$ is bounded in $W_\rad$, and since $|w_j| \le \eps_j^{(N-p)/p} \sup |w| \to 0$ uniformly, $w_j \wto 0$ (as in the proof of Theorem \ref{Theorem embedding} \textnormal{(\hyperlink{Theorem embedding.iii}{iii})}). For $v \in W_\rad$, \eqref{Talenti equation} gives
	\begin{align*}
		E'(w_j)\, v & = \int_{\R^N} |x|^{-b}\, |\nabla w_j|^{q-2}\, \nabla w_j \cdot \nabla v\, \mathrm{d}x - \lambda \int_{\R^N} |x|^{-c}\, |w_j|^{s-2}\, w_j\, v\, \mathrm{d}x\\
		& \quad - \mu_2 \int_{\R^N} |x|^{-c_2}\, |w_j|^{s_2 - 2}\, w_j\, v\, \mathrm{d}x,
	\end{align*}
	and all three terms are bounded by $o(1) \norm{v}$ by \eqref{dual estimates} and the remark following it. Hence $E'(w_j) \to 0$. Similarly,
	\begin{align*}
		E(w_j) = \frac{1}{p}\, X_1 - \frac{\mu_1}{p^\ast}\, K_1(w) + o(1) = \left(\frac{1}{p} - \frac{1}{p^\ast}\right) X_1 + o(1) = \beta_1 + o(1).
	\end{align*}
	Finally, no subsequence of $\seq{w_j}$ converges, since the only possible limit is $0$ while $X(w_j) = X_1 > 0$.

	\textnormal{(\hyperlink{Proposition sharp.ii}{ii})} Here $s_2 = q^\ast_b$ and $c_2 = 0$. We first identify $S_2$. Let $\theta = q/(q + b)$ and $M = Nq/(q + b) > q$, so that $q^\ast_b = Mq/(M - q)$. For radial $u$, the substitution $u(r) = v(t)$, $r = t^\theta$, gives
	\begin{align*}
		Y(u) = \omega_{N-1}\, \theta^{1-q} \int_0^\infty t^{M-1}\, |v'(t)|^q\, \mathrm{d}t, \quad K_2(u) = \omega_{N-1}\, \theta \int_0^\infty t^{M-1}\, |v(t)|^{q^\ast_b}\, \mathrm{d}t.
	\end{align*}
	By Bliss' inequality \cite{Bliss-1930}, the corresponding one-dimensional quotient is minimized by $v(t) = (1 + t^{q/(q-1)})^{-(M-q)/q}$. Thus the infimum of $Y(u)/K_2(u)^{q/q^\ast_b}$ over radial functions with $Y(u) < \infty$ is attained at
	\begin{align*}
		V(x) = \left(1 + |x|^{(q+b)/(q-1)}\right)^{-(N-q-b)/(q+b)}.
	\end{align*}
	Note that $|V(x)| \le C\, (1 + |x|)^{-\gamma}$ and $|\nabla V(x)| \le C\, (1 + |x|)^{-\gamma - 1}$ with $\gamma = (N - q - b)/(q - 1)$. Since the functions $\eta(\cdot/R)\, V$, where $\eta \in C^\infty_c(\R^N)$ is radial with $0 \le \eta \le 1$, $\eta = 1$ on $B_1$, and $\supp \eta \subset B_2$, belong to $W_\rad$ by Lemma \ref{Lemma density} and converge to $V$ in the norms defining $Y$ and $K_2$ as $R \to \infty$, we have $S_2 = Y(V)/K_2(V)^{q/q^\ast_b}$. Let $z = tV$ with $t > 0$ such that $Y(z) = \mu_2 K_2(z)$. As in part \textnormal{(\hyperlink{Proposition sharp.i}{i})}, $Y(z) = Y_2$, and $z$ satisfies
	\begin{equation} \label{Bliss equation}
		\int_{\R^N} |x|^{-b}\, |\nabla z|^{q-2}\, \nabla z \cdot \nabla v\, \mathrm{d}x = \mu_2 \int_{\R^N} |z|^{q^\ast_b - 2}\, z\, v\, \mathrm{d}x
	\end{equation}
	for all radial $v$ with $Y(v) < \infty$, in particular for all $v \in W_\rad$. However, $z$ need not belong to $W_\rad$, since $X(z)$ may be infinite, so we truncate.

	Let $\eps_j \searrow 0$, $a = (N - q - b)/q$, $z_j(x) = \eps_j^{-a}\, z(x/\eps_j)$, and $w_j = \eta\, z_j$. Each $z_j$ satisfies \eqref{Bliss equation}, $Y(z_j) = Y_2$, and $K_2(z_j) = K_2(z)$. By the change of variables $y = x/\eps_j$ and the weighted Hardy inequality,
	\begin{align*}
		Y(z_j - w_j) & \le C \int_{|y| > 1/\eps_j} |y|^{-b} \left(|\nabla z|^q + |y|^{-q}\, |z|^q\right) \mathrm{d}y \to 0, \\
		K_2(z_j - w_j) & \le \int_{|y| > 1/\eps_j} |z|^{q^\ast_b}\, \mathrm{d}y \to 0.
	\end{align*}
	Next,
	\begin{align*}
		X(w_j) \le C\, \eps_j^{N - p\,(a+1)} \int_{B_{2/\eps_j}} |\nabla z|^p\, \mathrm{d}y + C\, \eps_j^{N - pa} \int_{B_{2/\eps_j}} |z|^p\, \mathrm{d}y.
	\end{align*}
	Here $N - p\,(a + 1) = p\kappa > 0$. By the decay of $V$, the two integrals are $O(1 + \eps_j^{-(N - p\,(\gamma + 1))_+} + |\log \eps_j|)$ and $O(1 + \eps_j^{-(N - p\gamma)_+} + |\log \eps_j|)$, respectively. Since $\gamma > a$, it follows that $X(w_j) \le C\, (\eps_j^{p\kappa} |\log \eps_j| + \eps_j^{p\,(\gamma - a)}) \to 0$. So $\seq{w_j}$ is bounded in $W_\rad$. Moreover, $w_j(x) \to 0$ for every $x \ne 0$, since $|z_j(x)| \le C\, \eps_j^{\gamma - a}\, |x|^{-\gamma}$, and hence $w_j \wto 0$.

	For $v \in W_\rad$, we use \eqref{Bliss equation} for $z_j$ and write
	\begin{align*}
		E'(w_j)\, v
		& = \int_{\R^N} |x|^{-b} \left(|\nabla w_j|^{q-2}\, \nabla w_j - |\nabla z_j|^{q-2}\, \nabla z_j\right) \cdot \nabla v\, \mathrm{d}x\\
		& \quad - \mu_2 \int_{\R^N} \left(|w_j|^{q^\ast_b - 2}\, w_j - |z_j|^{q^\ast_b - 2}\, z_j\right) v\, \mathrm{d}x \\
		& \quad + \int_{\R^N} |\nabla w_j|^{p-2}\, \nabla w_j \cdot \nabla v\, \mathrm{d}x - \lambda \int_{\R^N} |x|^{-c}\, |w_j|^{s-2}\, w_j\, v\, \mathrm{d}x\\
		& \quad - \mu_1 \int_{\R^N} |x|^{-c_1}\, |w_j|^{s_1 - 2}\, w_j\, v\, \mathrm{d}x.
	\end{align*}
	The terms in the last two lines are $o(1) \norm{v}$ by \eqref{dual estimates}, since $X(w_j) \to 0$ and $w_j \wto 0$. For the first two lines, note that the rescaling $u \mapsto \eps_j^{a}\, u(\eps_j\, \cdot)$ preserves $Y$, $K_2$, and the corresponding dual norms. It maps $z_j$ to $z$ and $w_j$ to $\eta(\eps_j \cdot)\, z$, and $\eta(\eps_j \cdot)\, z \to z$ in the norms defining $Y$ and $K_2$. By the continuity of the maps $\xi \mapsto |\xi|^{q-2}\, \xi$ from $L^q(\R^N;|x|^{-b})$ to $L^{q'}(\R^N;|x|^{-b})$ and $u \mapsto |u|^{q^\ast_b - 2}\, u$ from $L^{q^\ast_b}(\R^N)$ to $L^{(q^\ast_b)'}(\R^N)$, together with Lemma \ref{Lemma critical segments} \textnormal{(\hyperlink{Lemma critical segments.ii}{ii})} applied to $v$, the first two lines are also $o(1) \norm{v}$. Hence $E'(w_j) \to 0$. Finally,
	\begin{align*}
		E(w_j) = \frac{1}{q}\, Y_2 - \frac{\mu_2}{q^\ast_b}\, K_2(z) + o(1) = \left(\frac{1}{q} - \frac{1}{q^\ast_b}\right) Y_2 + o(1) = \beta_2 + o(1),
	\end{align*}
	and no subsequence of $\seq{w_j}$ converges, since $w_j \wto 0$ and $Y(w_j) \to Y_2 > 0$.

	The final assertions follow from \textnormal{(\hyperlink{Proposition sharp.i}{i})}, \textnormal{(\hyperlink{Proposition sharp.ii}{ii})}, and Proposition \ref{Proposition PS}.
\end{proof}

\begin{remark}
	The proof of Proposition \ref{Proposition sharp} uses radial extremals for the corresponding single-phase inequalities and their decay properties to construct the required sequences in $W_\rad$. The same argument therefore applies at any point of $[Q_1,P_1]$ or $[Q_2,P_3]$ where this holds. For radial functions, the substitution used in part \textnormal{(\hyperlink{Proposition sharp.ii}{ii})} reduces every such point with $s_i > p$, respectively $s_i > q$, to a one-dimensional inequality of Bliss type. We have not pursued the verification of the required integrability of the extremals for the other phase in this generality.
\end{remark}

\subsubsection{Boundedness of Palais-Smale sequences} 

Following the proof of Mercuri and Perera \cite[Lemma 3.6]{Mercuri-Perera-2026}, we use the projection onto $\M$ and an interpolation inequality to show that \PS{} sequences are bounded. The interpolation relies on the subscaled part of the embedding range.

Let $\T_{\rm sub}$ denote the open triangle with the vertices $P_2$, $Q_1$, and $Q_2$. It lies above the scaling line, so that
\begin{align*}
	\rho(s',c') < \sigma\quad \text{for all }(s',c') \in \T_{\rm sub},
\end{align*}
and every point in this triangle satisfies the hypotheses of Lemma \ref{Lemma wedge}.

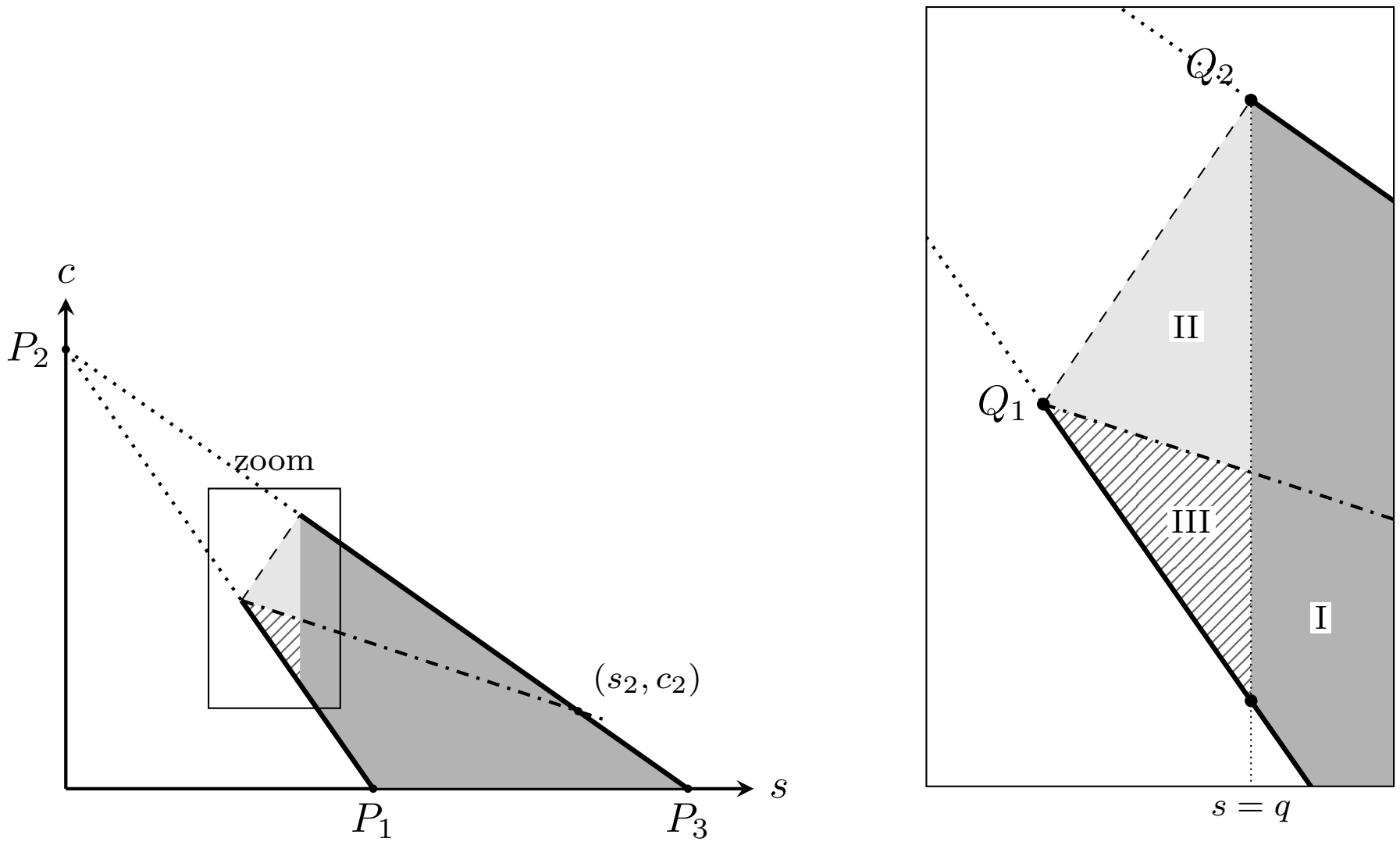
\begin{figure}[ht]
	\centering
	\newcommand{\Bcoords}{%
		\coordinate (P1) at (4.2,0);
		\coordinate (P2) at (0,6);
		\coordinate (P3) at (8.5,0);
		\coordinate (Q1) at (2.4,2.571);
		\coordinate (Q2) at (3.2,3.741);
		\coordinate (Q3) at (3.2,1.429);
		\coordinate (S2) at (7.0,1.059);
		\coordinate (R)  at (3.2,2.308);
	}
	\begin{tikzpicture}
		\begin{scope}[scale=0.62]
			\Bcoords
			\draw[-stealth, thick] (0,0) -- (9.4,0) node[right] {$s$};
			\draw[-stealth, thick] (0,0) -- (0,6.7) node[above] {$c$};
			\fill[black!30] (Q2) -- (Q3) -- (P1) -- (P3) -- cycle;
			\fill[black!10] (Q1) -- (Q2) -- (R) -- cycle;
			\fill[pattern=north east lines, pattern color=black!60] (Q1) -- (R) -- (Q3) -- cycle;
			\draw[thick, dotted] (P2) -- (Q1);
			\draw[thick, dotted] (P2) -- (Q2);
			\draw[very thick] (Q1) -- (P1);
			\draw[very thick] (Q2) -- (P3);
			\draw (P1) -- (P3);
			\draw[dashed] (Q1) -- (Q2);
			\draw[thick, dash dot] (Q1) -- ($(Q1)!1.08!(S2)$);
			\fill (P1) circle (1.6pt) node[below] {$P_1$};
			\fill (P2) circle (1.6pt) node[left] {$P_2$};
			\fill (P3) circle (1.6pt) node[below] {$P_3$};
			\fill (S2) circle (1.6pt) node[above right] {\footnotesize $(s_2,c_2)$};
			\draw[thin] (1.95,1.1) rectangle (3.75,4.1);
			\node[font=\footnotesize] at (2.85,4.45) {zoom};
		\end{scope}
		\begin{scope}[xshift=3cm, yshift=-2.4cm, scale=2.2]
			\Bcoords
			\begin{scope}
				\clip (1.95,1.1) rectangle (3.75,4.1);
				\fill[black!30] (Q2) -- (Q3) -- (P1) -- (P3) -- cycle;
				\fill[black!10] (Q1) -- (Q2) -- (R) -- cycle;
				\fill[pattern=north east lines, pattern color=black!60] (Q1) -- (R) -- (Q3) -- cycle;
				\draw[thick, dotted] (P2) -- (Q1);
				\draw[thick, dotted] (P2) -- (Q2);
				\draw[very thick] (Q1) -- (P1);
				\draw[very thick] (Q2) -- (P3);
				\draw[dashed] (Q1) -- (Q2);
				\draw[densely dotted] (3.2,0) -- (Q2);
				\draw[thick, dash dot] (Q1) -- (S2);
			\end{scope}
			\draw[thin] (1.95,1.1) rectangle (3.75,4.1);
			\fill (Q1) circle (0.7pt) node[left] {$Q_1$};
			\fill (Q2) circle (0.7pt) node[above left] {$Q_2$};
			\fill (Q3) circle (0.7pt);
			\node[below, font=\footnotesize] at (3.2,1.1) {$s = q$};
			\node[font=\footnotesize, fill=white, inner sep=1pt] at (3.47,1.75) {I};
			\node[font=\footnotesize, fill=white, inner sep=1pt] at (2.95,2.87) {II};
			\node[font=\footnotesize, fill=white, inner sep=1pt] at (2.97,2.12) {III};
		\end{scope}
	\end{tikzpicture}
	\caption{The region in condition (B) (schematic, drawn for $q < p^\ast$), with an enlargement near $Q_1$ and $Q_2$ on the right. For $\lambda > 0$ and $\rho \ge \sigma$, the point $(s,c) \in \Tc$ lies on or below the scaling line $Q_1 Q_2$ (dashed). If $\mu_1 > 0$, then (B) holds for all such points. If $\mu_1 = 0$, then (B) holds in the region I (dark), where $s > q$, and in the region II (light), which lies strictly above the line through $Q_1$ and $(s_2,c_2)$ (dash-dotted) and contains the open scaling line $Q_1 Q_2$. Region III (hatched), which lies in the triangle with vertices $Q_1$, $Q_2$, and $(q,N\, (1 - q/p^\ast))$ on or below the dash-dotted line, is not covered when $\mu_1 = 0$.}
	\label{fig:B}
\end{figure}

\begin{lemma} \label{Lemma bounded}
	Assume \textnormal{(\hyperlink{H1}{H$_1$})}--\textnormal{(\hyperlink{H4}{H$_4$})}, and let $\lambda > 0$ and $\rho \ge \sigma$. If
	\begin{enumerate}
		\item[\textnormal{(B)}]\hypertarget{B}{}
			$s > q$, or $\mu_1 > 0$, or $(s,c)$ lies strictly above the line through $Q_1$ and $(s_2,c_2)$ (see Figure \ref{fig:B}),
	\end{enumerate}
	then every \PS{c} sequence of $E$ is bounded.
\end{lemma}

\begin{proof}
	Let $\seq{u_j}$ be a \PS{c} sequence.

	\emph{Case 1: $s > q$.} Let
	\begin{align*}
		\vartheta = \min \left(\set{s} \cup \set{s_i \colon \mu_i > 0}\right) > q.
	\end{align*}
	Then
	\begin{align*}
		&E(u_j) - \frac{1}{\vartheta}\, E'(u_j)\, u_j\\
		& = \left(\frac{1}{p} - \frac{1}{\vartheta}\right) X(u_j) + \left(\frac{1}{q} - \frac{1}{\vartheta}\right) Y(u_j) + \lambda \left(\frac{1}{\vartheta} - \frac{1}{s}\right) J(u_j) \\
		& \quad + \sum_{\mu_i > 0} \mu_i \left(\frac{1}{\vartheta} - \frac{1}{s_i}\right) K_i(u_j) \ge \left(\frac{1}{q} - \frac{1}{\vartheta}\right) \big[X(u_j) + Y(u_j)\big],
	\end{align*}
	while the left-hand side is $c + \o(1) + \o(\norm{u_j})$. Since $1 < p < q$, we conclude that $\seq{u_j}$ is bounded.

	\emph{Case 2: $s \le q$.} We first show that there exist $i \in \set{1,2}$ with $\mu_i > 0$, $\eta \in (0,1)$, and $(s',c') \in \T_{\rm sub}$ such that
	\begin{equation} \label{interpolation point}
		(s,c) = (1 - \eta)(s',c') + \eta\,  (s_i,c_i).
	\end{equation}
	If $\mu_1 > 0$, take $i = 1$. The points $P_2$, $Q_1$, and $(s_1,c_1)$ lie on the line $c = N\, (1 - s/p^\ast)$, with $Q_1$ between the other two. The set of points of $\Tc$ with $\rho \ge \sigma$ and $s \le q$ is contained in the triangle with the vertices $Q_1$, $Q_2$, and $(q,N\, (1 - q/p^\ast))$, the last vertex lying between $Q_1$ and $(s_1,c_1)$ since $p < q < s_1$, and it avoids the lines $c = N\, (1 - s/p^\ast)$ and $c = N\, (1 - s/q^\ast_b)$ as well as the segment joining $Q_2$ and $(s_1,c_1)$. So $(s,c)$ lies in the open triangle with the vertices $P_2$, $(s_1,c_1)$, and $Q_2$. Hence the ray from $(s_1,c_1)$ through $(s,c)$ leaves this triangle through an interior point of the side $P_2 Q_2$, which is also a side of $\T_{\rm sub}$, and the points of this ray just before the exit point lie in $\T_{\rm sub}$. This gives \eqref{interpolation point}. If $\mu_1 = 0$, then $\mu_2 > 0$ and $(s,c)$ lies strictly above the line through $Q_1$ and $(s_2,c_2)$ by $(B)$. Take $i = 2$. The points $P_2$, $Q_2$, and $(s_2,c_2)$ lie on the line $c = N\, (1 - s/q^\ast_b)$, so $(s,c)$ lies in the open triangle with the vertices $P_2$, $Q_1$, and $(s_2,c_2)$, and the ray from $(s_2,c_2)$ through $(s,c)$ leaves this triangle through an interior point of the side $P_2 Q_1$, which is also a side of $\T_{\rm sub}$. This again gives \eqref{interpolation point}.

	Now suppose that $\norm{u_j} \to \infty$ along a subsequence. Then $t_j := I_\sigma(u_j)^{1/\sigma} \to \infty$. Set $\widetilde{u}_j = \pi(u_j) \in \M$, so that $u_j = (\widetilde{u}_j)_{t_j}$ and hence $\norm{u_j} \le t_j^{\sigma/p} \norm{\widetilde{u}_j} = \O(t_j^{\sigma/p})$ by \eqref{eq.1.prop1}, since $\M$ is bounded. The scaling relations \eqref{scaling relations} give
	\begin{equation}\label{proj E}
		\begin{aligned}
			E(u_j)
			& = t_j^\sigma \left[\frac{1}{p}\, X(\widetilde{u}_j) + \frac{1}{q}\, Y(\widetilde{u}_j)\right] - \frac{\lambda}{s}\, t_j^\rho\, J(\widetilde{u}_j) - \sum_{k=1}^2 \frac{\mu_k}{s_k}\, t_j^{\rho_k}\, K_k(\widetilde{u}_j)\\
			&= c + \o(1),
		\end{aligned}
	\end{equation}
	\begin{equation}\label{proj E'}
		\begin{aligned}
			E'(u_j)\, u_j
			 & = t_j^\sigma\, \big[X(\widetilde{u}_j) + Y(\widetilde{u}_j)\big] - \lambda\, t_j^\rho\, J(\widetilde{u}_j) - \sum_{k=1}^2 \mu_k\, t_j^{\rho_k}\, K_k(\widetilde{u}_j)\\
			 &= \o(t_j^{\sigma/p}).
		\end{aligned}
	\end{equation}
	Since $X$ and $Y$ are bounded on $\M$ and all terms involving $\lambda$, $\mu_1$, and $\mu_2$ are nonnegative, we infer from \eqref{proj E'} that, for the index $i$ chosen above,
	\begin{align*}
		K_i(\widetilde{u}_j) \le Ct_j^{- (\rho_i - \sigma)}
	\end{align*}
	for some constant $C > 0$. By \eqref{interpolation point} and H\"{o}lder's inequality with exponents $1/(1 - \eta)$ and $1/\eta$,
	\begin{equation} \label{J interpolation}
		\begin{aligned}
			J(\widetilde{u}_j)
			& = \int_{\R^N} \big(|x|^{-c'}\, |\widetilde{u}_j|^{s'}\big)^{1 - \eta} \left(|x|^{-c_i}\, |\widetilde{u}_j|^{s_i}\right)^\eta\, \mathrm{d}x \\
			& \le \left(\int_{\R^N} |x|^{-c'}\, |\widetilde{u}_j|^{s'}\, \mathrm{d}x\right)^{1 - \eta} K_i(\widetilde{u}_j)^\eta \\
			& \le Ct_j^{- \eta\, (\rho_i - \sigma)},
		\end{aligned}
	\end{equation}
	because the first integral is bounded by Lemma \ref{Lemma wedge}. As $\rho$ is affine in $(s,c)$, equation \eqref{interpolation point} implies
	\begin{align*}
		\rho = (1 - \eta)\, \rho' + \eta\,  \rho_i, \quad \rho' = \rho(s',c') < \sigma.
	\end{align*}
	Combining this with \eqref{J interpolation}, we arrive at
	\begin{equation} \label{J-bd}
		t_j^{\rho - \sigma}\, J(\widetilde{u}_j) \le Ct_j^{\rho - \sigma - \eta\, (\rho_i - \sigma)} = Ct_j^{- (1 - \eta)(\sigma - \rho')} \to 0.
	\end{equation}
	Let $\theta = \min \{s_k \colon \mu_k > 0\}$. Multiplying \eqref{proj E} by $\theta$, subtracting \eqref{proj E'}, and dividing by $t_j^\sigma$, we obtain
	\begin{align*}
		& \left(\frac{\theta}{p} - 1\right) X(\widetilde{u}_j) + \left(\frac{\theta}{q} - 1\right) Y(\widetilde{u}_j) \\
		& = \lambda \left(\frac{\theta}{s} - 1\right) t_j^{\rho - \sigma}\, J(\widetilde{u}_j) + \sum_{k=1}^2 \mu_k \left(\frac{\theta}{s_k} - 1\right) t_j^{\rho_k - \sigma}\, K_k(\widetilde{u}_j) + o(1) \le \o(1)
	\end{align*}
	by \eqref{J-bd}, since $\theta \le s_k$ whenever $\mu_k > 0$. This is a contradiction, because $\theta > q > p$ by \textnormal{(\hyperlink{H4}{H$_4$})}, whereas
	\begin{align*}
		\frac1pX(\widetilde{u}_j) +\frac1qY(\widetilde{u}_j)=1
	\end{align*}
	on $\M$. Hence $\seq{u_j}$ is bounded.
\end{proof}

\subsubsection{Superscaled perturbations}

\begin{theorem} \label{Theorem large lambda}
	Assume \textnormal{(\hyperlink{H1}{H$_1$})}--\textnormal{(\hyperlink{H4}{H$_4$})} and $(B)$, and let $\rho > \sigma$. Then, for any $m \ge 1$, there exists $\lambda^\ast_m > 0$ such that equation \eqref{19} has $m$ distinct pairs of nontrivial weak radial solutions at energy levels in $(0,\beta^\ast)$ for all $\lambda > \lambda^\ast_m$. In particular, the number of solutions tends to infinity as $\lambda \to \infty$.
\end{theorem}

\begin{proof}
	For $\lambda > 0$, Lemma \ref{Lemma bounded} and Proposition \ref{Proposition PS} show that $E$ satisfies the \PS{c} condition for all $c < \beta^\ast$. We apply Corollary \ref{Corollary 1} with $c^\ast = \beta^\ast$. Let $V$ be an $m$-dimensional subspace of $W_\rad$, let $S_V$ denote its unit sphere, and let $A_0 = \pi(S_V)$. Then $A_0$ is a compact symmetric subset of $\M$, and $i(A_0) \ge i(S_V) = m$ by the monotonicity of the cohomological index. As $J > 0$ on the compact set $A_0$, we have $j_0 := \min_{A_0} J > 0$, and for $u \in A_0$ and $t \ge 0$,
	\begin{align*}
		E(u_t) \le t^\sigma - \frac{\lambda\,  j_0}{s}\, t^\rho.
	\end{align*}
	Since $\rho > \sigma$, it follows that
	\begin{align*}
		\sup_{u \in A_0,\, t \ge 0}\, E(u_t) \le \left(1 - \frac{\sigma}{\rho}\right)\!\left(\frac{\sigma s}{\rho\,  \lambda\,  j_0}\right)^{\sigma/(\rho - \sigma)} \to 0 \quad \text{as } \lambda \to \infty.
	\end{align*}
	Hence there exists $\lambda^\ast_m > 0$ such that $\sup_{u \in X}\, E(u) < \beta^\ast$ for all $\lambda > \lambda^\ast_m$ and all $R > 0$, where $X = \{u_t \colon u \in A_0,\, 0 \le t \le R\}$. Fix $\lambda > \lambda^\ast_m$. Choosing $R$ sufficiently large that $R^\sigma - (\lambda\,  j_0/s)\, R^\rho \le 0$, we find $\sup_{u \in A}\, E(u) \le 0$, where $A = \set{u_R \colon u \in A_0}$. On the other hand, since $J$, $K_1$, and $K_2$ are bounded on the bounded set $\M$, for $u \in \M$ and $t \ge 0$,
	\begin{align*}
		E(u_t) \ge t^\sigma - C\left(t^\rho +\sum_{\mu_i>0} t^{\rho_i}\right)
	\end{align*}
	for some constant $C > 0$. Since $\rho>\sigma$ and $\rho_i>\sigma$ whenever $\mu_i>0$, we conclude that
	\begin{align*}
		\inf_{u \in \M_r}\, E(u) > 0
	\end{align*}
	for all sufficiently small $r \in (0,R)$, where $\M_r = \{u_r \colon u \in \M\}$. Corollary \ref{Corollary 1} now provides $m$ distinct pairs of critical points of $E$ at levels in $(0,\beta^\ast)$.
\end{proof}

\subsubsection{Subscaled perturbations: a concave-convex type result}

\begin{theorem} \label{Theorem concave convex}
	Assume \textnormal{(\hyperlink{H1}{H$_1$})}--\textnormal{(\hyperlink{H4}{H$_4$})} and let $\rho < \sigma$. Then there exists $\lambda_\ast > 0$ such that for all $\lambda \in (0,\lambda_\ast)$, equation \eqref{19} has a sequence of weak radial solutions $\seq{u_k}$ such that $E(u_k) < 0$ and $E(u_k) \nearrow 0$.
\end{theorem}

\begin{proof}
	We use a truncation of $E$ inspired by Garc\'{\i}a Azorero and Peral Alonso \cite{GarciaAzorero-PeralAlonso-1991} (see also Mercuri and Perera \cite[Section 3.4.4]{Mercuri-Perera-2026}). For $u \ne 0$, write $u = \widetilde{u}_t$ with $\widetilde{u} = \pi(u) \in \M$ and $t = I_\sigma(u)^{1/\sigma}$. Since $J$, $K_1$, and $K_2$ are bounded on $\M$, the scaling relations \eqref{scaling relations} provide constants $C_0, C_1, C_2 > 0$ such that
	\begin{align*}
		J(u) \le C_0\, I_\sigma(u)^{\rho/\sigma}, \quad K_i(u) \le C_i\, I_\sigma(u)^{\rho_i/\sigma} \quad \text{for all }u \in W_\rad,\, i = 1, 2.
	\end{align*}
	Consequently,
	\begin{align*}
		E(u) \ge g_\lambda(I_\sigma(u)),
	\end{align*}
	where
	\begin{align*}
		g_\lambda(\tau) = \tau - \frac{\lambda\,  C_0}{s}\, \tau^{\rho/\sigma} - \frac{\mu_1 C_1}{s_1}\, \tau^{\rho_1/\sigma} - \frac{\mu_2\,  C_2}{s_2}\, \tau^{\rho_2/\sigma}, \quad \tau \ge 0.
	\end{align*}
	We may write
	\begin{align*}
		g_\lambda(\tau) = \tau^{\rho/\sigma} \left(k(\tau)-\frac{\lambda C_0}{s}\right),
	\end{align*}
	where
	\begin{align*}
		k(\tau) = \tau^{1 - \rho/\sigma} - \frac{\mu_1 C_1}{s_1}\, \tau^{(\rho_1 - \rho)/\sigma} - \frac{\mu_2\,  C_2}{s_2}\, \tau^{(\rho_2 - \rho)/\sigma}.
	\end{align*}
	For every $i$ with $\mu_i>0$, assumption \textnormal{(\hyperlink{H4}{H$_4$})} implies $\rho_i>\sigma$. Thus,
	\begin{align*}
		0<1-\frac{\rho}{\sigma} < \frac{\rho_i-\rho}{\sigma} \quad\text{whenever }\mu_i>0.
	\end{align*}
	It follows that $k(0) = 0$ and that $k$ is strictly increasing on $[0,\tau_1]$ for some $\tau_1 > 0$. Let $\lambda_0 = s\, k(\tau_1)/C_0$. For $\lambda \in (0,\lambda_0)$, let $R_1(\lambda) \in (0,\tau_1)$ be the unique point satisfying $k(R_1(\lambda)) = \lambda\,  C_0/s$ and set $R_2 = \tau_1$. Then
	\begin{align*}
		g_\lambda < 0 \quad\text{on } (0,R_1(\lambda)), \qquad g_\lambda > 0 \quad\text{on } (R_1(\lambda),R_2].
	\end{align*}
	Choose a smooth function $\xi_\lambda \colon [0,\infty) \to [0,1]$ such that $\xi_\lambda = 1$ on $[0,R_1(\lambda)]$ and $\xi_\lambda = 0$ on $[R_2,\infty)$, and define
	\begin{align*}
		\widetilde{E}(u) = \xi_\lambda(I_\sigma(u))\, E(u), \quad u \in W_\rad.
	\end{align*}
	Then $\widetilde{E}$ is an even $C^1$-functional with $\widetilde{E}(0) = 0$ and it is bounded from below since $\widetilde{E} = 0$ on $\{u \in W_\rad \colon I_\sigma(u) \ge R_2\}$ and $E$ is bounded on the bounded set $\{u \in W_\rad \colon I_\sigma(u) < R_2\}$. If $\widetilde{E}(u) < 0$, then $I_\sigma(u) < R_1(\lambda)$ since $\widetilde{E}(u) \ge \xi_\lambda(I_\sigma(u))\, g_\lambda(I_\sigma(u)) \ge 0$ if $I_\sigma(u) \in [R_1(\lambda),R_2]$ and $\widetilde{E}(u) = 0$ if $I_\sigma(u) \ge R_2$. Hence $I_\sigma < R_1(\lambda)$ in a neighborhood of $u$ by continuity, where $\widetilde{E} = E$, so
	\begin{align*}
		\widetilde{E}(u) = E(u), \qquad \widetilde{E}'(u) = E'(u).
	\end{align*}
	In particular, critical points of $\widetilde{E}$ at negative levels are critical points of $E$ at the same levels.

	Next we show that $\widetilde{E}$ satisfies the \PS{c} condition for all $c < 0$, provided $\lambda$ is sufficiently small. Let $c < 0$ and suppose that $\seq{u_j}$ is a \PS{c} sequence of $\widetilde{E}$. For all sufficiently large $j$, $\widetilde{E}(u_j) < 0$, so $\seq{u_j}$ is a \PS{c} sequence of $E$ with $I_\sigma(u_j) < R_1(\lambda) < \tau_1$. In particular, $\seq{u_j}$ is bounded independently of $\lambda \in (0,\lambda_0)$, say $\norm{u_j} \le C_3$. By Proposition \ref{Proposition PS}, after passing to a subsequence, $u_j$ converges weakly to a weak radial solution $u$ of equation \eqref{19}, with $\norm{u} \le C_3$. If $u_j \nrightarrow u$, then $c \ge E(u) + \beta^\ast$. Using \eqref{energy critical}, we derive
	\begin{align*}
		E(u) \ge - \left(1 - \frac{\rho}{\sigma}\right) \frac{\lambda}{s}\, J(u) \ge - C_4 \lambda
	\end{align*}
	for some constant $C_4 > 0$ independent of $\lambda$. Therefore, if $u_j \nrightarrow u$, then $0 > c \ge \beta^\ast - C_4 \lambda$, which is impossible when $\lambda \le \beta^\ast/C_4$.

	Set
	\begin{align*}
		\lambda_\ast = \min\left\{\lambda_0,\frac{\beta^\ast}{C_4}\right\}.
	\end{align*}
	Now let $\lambda \in (0,\lambda_\ast)$ and apply Proposition \ref{Proposition negative levels} to $\widetilde{E}$. Let $c_k$ be as in that proposition and note that $c_k > - \infty$ since $\widetilde{E}$ is bounded from below. For $k \ge 1$, as in the proof of Theorem \ref{Theorem large lambda}, let $V$ be a $k$-dimensional subspace of $W_\rad$, let $S_V$ be its unit sphere, let $A = \pi(S_V)$, so that $A$ is a compact symmetric subset of $\M$ with $i(A) \ge k$, and let $j_k = \min_A J > 0$. For $t > 0$, $M_t = \set{u_t \colon u \in A}$ is a symmetric subset of $W_\rad \setminus \set{0}$ with $i(M_t) = i(A) \ge k$ since $A \to M_t,\, u \mapsto u_t$ is an odd homeomorphism, and
	\begin{align*}
		E(u_t) \le t^\sigma - \frac{\lambda\,  j_k}{s}\, t^\rho \quad \text{for all }u \in A.
	\end{align*}
	As $\rho < \sigma$, taking $t$ so small that $t^\sigma - (\lambda\,  j_k/s)\, t^\rho < 0$ and $t^\sigma \le R_1(\lambda)$ gives
	\begin{align*}
		c_k \le \sup_{u \in M_t}\, \widetilde{E}(u) = \sup_{u \in M_t}\, E(u) < 0.
	\end{align*}
	So $c_1 \le c_2 \le \dotsb \to 0$ is a sequence of negative critical values of $\widetilde{E}$, and hence of $E$, by Proposition \ref{Proposition negative levels}.
\end{proof}

\subsection{Concluding remarks and open problems} \label{Section open problems}

We conclude with some open problems.
\begin{enumerate}
	\item[\textnormal{(i)}]
		\emph{Two critical terms without perturbation.}
		When $\lambda = 0$ and $\mu_1, \mu_2 > 0$, the Poho\v{z}aev identity does not rule out nontrivial solutions of equation \eqref{19}. However, every nontrivial solution has energy at least $\beta_1 + \beta_2$ (see Remark \ref{Remark lambda zero}). Does equation \eqref{19} admit a nontrivial weak radial solution in this case?
	\item[\textnormal{(ii)}]
		\emph{Boundedness of \PS{} sequences.}
		When $\mu_1 = 0$, Lemma \ref{Lemma bounded} does not cover the points $(s,c)$ below the scaling line with $s \le q$ that lie on or below the line through $Q_1$ and $(s_2,c_2)$.
	\item[\textnormal{(iii)}]
		\emph{Nonradial solutions.}
		The embeddings of the whole space $W$ cover a smaller range (see Remark \ref{Remark 1}). It would be interesting to obtain nonradial solutions in this smaller range.
\end{enumerate}


\end{document}